\documentclass{agtart_a}
\pdfoutput=1

\usepackage{pinlabel}

\title[Knots with unknotting number 1 and essential Conway spheres]{Knots with unknotting number 1\\ and essential Conway spheres}

\author{C\,McA Gordon}
\givenname{Cameron McA}
\surname{Gordon}
\address{Department of Mathematics\\The University of Texas at 
Austin\\\newline
1 University Station C1200\\Austin, TX 78712-0257\\USA}
\email{gordon@math.utexas.edu}
\urladdr{}

\author{John Luecke}
\givenname{John}
\surname{Luecke}
\email{luecke@math.utexas.edu}
\urladdr{}

\volumenumber{6}
\issuenumber{}
\publicationyear{2006}
\papernumber{72}
\startpage{2051}
\endpage{2116}

\doi{}
\MR{}
\Zbl{}

\keyword{unknotting number 1}
\keyword{Conway spheres}
\keyword{algebraic knots}
\keyword{mutation}
\subject{primary}{msc2000}{57N10}
\subject{secondary}{msc2000}{57M25}

\received{9 January 2006}
\revised{}
\accepted{29 September 2006}
\published{19 November 2006}
\publishedonline{19 November 2006}
\proposed{}
\seconded{}
\corresponding{}
\editor{CPR}
\version{}

\arxivreference{math.GT/0601265}

\AtBeginDocument{\let\bar\wbar\let\tilde\wtilde\let\hat\what%
\def\notin{\not\in}}
\numberwithin{figure}{section}

\makeatletter
\def\cnewtheorem#1[#2]#3{\newtheorem{#1}{#3}[section]
\expandafter\let\csname c@#1\endcsname\c@thm}
\makeatother

\newtheorem{thm}{Theorem}[section]
\cnewtheorem{lem}[thm]{Lemma}
\cnewtheorem{quest}[thm]{Question}
\newtheorem*{Quest}{Question}
\newtheorem*{conjecture}{Conjecture}
\newtheorem*{thmJknot}{\fullref{thm2.1}} 
\newtheorem*{thmMain-bis}{\fullref{thm:Main}} 
\newtheorem*{cor6.3bis}{\fullref{cor4.2}} 
\newtheorem*{thm7.1bis}{\fullref{thm:M}} 
\newtheorem*{thm8.2bis}{\fullref{thm8.3}} 
\newtheorem*{thm11.1bis}{\fullref{thm11.1}} 
\newtheorem*{thm11.5bis}{\fullref{thm8.4}} 
 
\cnewtheorem{cor}[thm]{Corollary}
\cnewtheorem{prop}[thm]{Proposition}
\cnewtheorem{adden}[thm]{Addendum}
\newtheorem*{Addendum}{Addendum}
\newtheorem{claim}{Claim}
\newtheorem*{Claim}{Claim}

\theoremstyle{definition}
\newtheorem*{Defin}{Definition}
\cnewtheorem{Definition}[thm]{Definition}
\makeautorefname{Definition}{Definition}
\newtheorem*{remark}{Remark}
\newtheorem*{remarks}{Remarks}

\def\Fix{\operatorname{Fix}}
\def\id{\operatorname{id}}
\def\Int{\operatorname{int}}
\def\Mut{\operatorname{Mut}}
\def\rel{\operatorname{rel}}
\def\seif{\operatorname{Seif}}
\def\que{{\mathbb Q}}
\def\real{{\mathbb R}}
\def\zed{{\mathbb Z}}
\def\A{{\mathcal A}}
\def\B{{\mathcal B}}
\def\C{{\mathcal C}}

\def\M{{\mathcal M}}
\def\P{{\mathcal P}}
\def\R{{\mathcal R}}
\def\S{{\mathcal S}}
\def\T{{\mathcal T}}
\def\re{{\mathcal R}}
\def\bT{{\mathbf T}}

\def\wt#1{\widetilde#1}
\def\ep{\varepsilon}

\begin{document}
\begin{asciiabstract}
For a knot K in S^3, let T(K) be the characteristic toric sub-orbifold
of the orbifold (S^3,K) as defined by Bonahon and Siebenmann.  If K
has unknotting number one, we show that an unknotting arc for K can
always be found which is disjoint from T(K), unless either K is an
EM-knot (of Eudave-Munoz) or (S^3,K) contains an EM-tangle after
cutting along T(K).  As a consequence, we describe exactly which large
algebraic knots (ie algebraic in the sense of Conway and containing
an essential Conway sphere) have unknotting number one and give a
practical procedure for deciding this (as well as determining an
unknotting crossing).  Among the knots up to 11 crossings in Conway's
table which are obviously large algebraic by virtue of their
description in the Conway notation, we determine which have unknotting
number one.  Combined with the work of Ozsvath-Szabo, this
determines the knots with 10 or fewer crossings that have unknotting
number one.  We show that an alternating, large algebraic knot with
unknotting number one can always be unknotted in an alternating
diagram.

As part of the above work, we determine the hyperbolic knots in a solid 
torus which admit a non-integral, toroidal Dehn surgery. 
Finally, we show that having unknotting number one is invariant under 
mutation.
\end{asciiabstract}

\begin{htmlabstract}
<p class="noindent">
For a knot K in S<sup>3</sup>, let T(K) be the characteristic toric
sub-orbifold of the orbifold (S<sup>3</sup>,K) as defined by Bonahon&ndash;Siebenmann.
If K has unknotting number one, we show that an unknotting arc for K
can always be found which is disjoint from  T(K), unless either K is
an EM&ndash;knot (of Eudave-Munoz) or (S<sup>3</sup>,K) contains an EM&ndash;tangle after
cutting along T(K).
As a consequence, we describe exactly which large algebraic knots (ie,
algebraic in the sense of Conway and containing an essential Conway sphere)
have unknotting number one and give a practical procedure for deciding
this (as well as determining an unknotting crossing).
Among the knots up to 11&nbsp;crossings in Conway's table which are
obviously large algebraic by virtue of their description in the Conway
notation,  we determine which have unknotting number one.
Combined with the work of Ozsv&aacute;th&ndash;Szab&oacute;, this determines the knots with
10 or fewer crossings that have unknotting number one.
We show that an alternating, large algebraic knot with unknotting number
one can always be unknotted in an alternating diagram.
</p>
<p class="noindent">
As part of the above work, we determine the hyperbolic knots in a solid
torus which admit a non-integral, toroidal Dehn surgery.
Finally, we show that having unknotting number one is invariant under
mutation.
</p>
\end{htmlabstract}

\begin{abstract}
For a knot $K$ in $S^3$, let $\mathbf T(K)$ be the characteristic toric 
sub-orbifold of the orbifold $(S^3,K)$ as defined by Bonahon--Siebenmann. 
If $K$ has unknotting number one, we show that an unknotting arc for $K$ 
can always be found which is disjoint from $\mathbf T(K)$, unless either $K$ is
an EM--knot (of Eudave-Mu\~noz) or $(S^3,K)$ contains an EM--tangle after 
cutting along $\mathbf T(K)$.
As a consequence, we describe exactly which large algebraic knots (ie, 
algebraic in the sense of Conway and containing an essential Conway sphere) 
have unknotting number one and give a practical procedure for deciding 
this (as well as determining an unknotting crossing). 
Among the knots up to 11~crossings in Conway's table which are 
obviously large algebraic by virtue of their description in the Conway 
notation,  we determine which have unknotting number one.
Combined with the work of Ozsv\'ath--Szab\'o, this determines the knots with 
10 or fewer crossings that have unknotting number one. 
We show that an alternating, large algebraic knot with unknotting number
one can always be unknotted in an alternating diagram. 

As part of the above work, we determine the hyperbolic knots in a solid 
torus which admit a non-integral, toroidal Dehn surgery. 
Finally, we show that having unknotting number one is invariant under 
mutation.
\end{abstract}

\maketitle

\section{Introduction}

Montesinos showed \cite{Mon} that if a knot $K$ has unknotting
number~1 then its double branched cover $M$ can be obtained by a
half-integral Dehn surgery on some knot $K^*$ in $S^3$.  Consequently,
theorems about Dehn surgery can sometimes be used to give necessary
conditions for a knot $K$ to have unknotting number~1.  For instance,
$H_1(M)$ must be cyclic, and the $\que/\zed$--valued linking form on
$H_1(M)$ must have a particular form (Lickorish \cite{L}).  If $K$ is
a 2--bridge knot, then $M$ is a lens space, and hence, by the Cyclic
Surgery Theorem (Culler, Gordon, Luecke and Shalen \cite{CGLS}), $K^*$
must be a torus knot.  In this way the 2--bridge knots with unknotting
number~1 have been completely determined (Kanenobu and Murakami
\cite{KM}).  Another example is Scharlemann's theorem that unknotting
number~1 knots are prime \cite{S1}; this can be deduced from the fact
(proved later, however by Gordon and Luecke \cite{GL1}) that only
integral Dehn surgeries can give reducible manifolds (see Zhang
\cite{Z}).  Finally, we mention the recent work of Ozsv\'ath and
Szab\'o \cite{OS}, in which the Heegaard Floer homology of $M$ is used
to give strong restrictions on when $K$ can have unknotting number~1,
especially if $K$ is alternating.

The present paper explores another example of this connection.  Here,
the Dehn surgery theorem is the result of Gordon and Luecke \cite{GL3}
that the hyperbolic knots with non-integral toroidal Dehn surgeries
are precisely the Eudave-Mu\~noz knots $k(\ell,m,n,p)$ \cite{E1}; this
gives information about when a knot $K$ whose double branched cover
$M$ is toroidal can have unknotting number~1.

First we clarify the extension of the main result of \cite{GL3} to 
knots in solid tori that is described in the Appendix of \cite{GL3}. 
In Section~3, we define a family of hyperbolic knots 
$J_\ep (\ell,m)$ in a solid 
torus, $\ep \in \{1,2\}$ and $\ell,m$ integers, each of which admits
a half-integer surgery yielding a toroidal manifold. 
(The knots $J_\ep (\ell,m)$ in the solid torus are the 
analogs of the knots $k(\ell,m,n,p)$ in the 3--sphere.)  
We then use \cite{GL3} to show that these are the only such:

\begin{thmJknot}  		
Let $J$ be a knot in a solid torus whose exterior is irreducible and atoroidal.
Let $\mu$ be the meridian of $J$ and suppose that $J(\gamma)$ contains an 
essential  torus for some $\gamma$ with $\Delta (\gamma,\mu)\ge2$. 
Then $\Delta (\gamma,\mu) =2$ and $J= J_\ep (\ell,m)$ for some $\ep,\ell,m$.
\end{thmJknot}

This theorem along with the main result of \cite{GL3} then allows us to 
describe in \fullref{thm:JSJ} the relationship between the torus 
decomposition of the exterior of a knot $K$ in $S^3$ and the torus 
decompositon of any non-integral surgery on $K$. 
In particular, \fullref{thm:JSJ} says that the canonical tori of the 
exterior of $K$ and of the Dehn surgery will be the same unless 
$K$ is a cable knot 
(in which case an essential torus of the knot exterior can become 
compressible), or $K$ is  a $k(\ell,m,n,p)$, or $K$ is a satellite with 
pattern 
$J_\ep (\ell,m)$ (in the latter cases a new essential torus is created).

We apply these theorems about non-integral Dehn surgeries to address
questions about unknotting number.  The knots $k(\ell,\!m,\!n,\!p)$
($J_\ep (\ell,m)$) are strongly invertible.  Their quotients under the
involutions give rise to {\em EM--knots}, $\!K(\ell,\!m,\!n,\!p)$
({\em EM--tangles}, $\A_\ep (\ell,m)$ resp.) which have essential
Conway spheres and yet can be unknotted (trivialized, resp.) by a
single crossing change.  \fullref{thm:Main} descibes when a knot with an
essential Conway sphere or 2--torus can have unknotting number~1.
This is naturally stated in the context of the characteristic
decomposition of a knot along toric 2--suborbifolds given by Bonahon
and Siebenmann in \cite{BS}.  The characteristic torus decomposition
of the double branched cover of a knot $K$ corresponds to the
characteristic decomposition of the orbifold $O(K)$, where $O(K)$
refers to $S^3$ thought of as an orbifold with singular set $K$ and
cone angle $\pi$ (see \cite{BS}).  This decomposition of $O(K)$ is
along Conway spheres and along tori disjoint from $K$, the collection
of which is denoted $\bT (K)$.  When $O(K)$ is cut along $\bT(K)$,
$\seif(K)$ denotes the components corresponding to Seifert-fibered
components of the canonical torus decomposition in the double branched
cover.  An {\em unknotting arc}, $(a,\partial a)$, for $K$ is an arc
such that $a\cap K=\partial a$ that guides a crossing move that
unknots $K$.  Under the correspondence between crossing changes and
Dehn surgeries in the double branched cover, \fullref{thm:JSJ} becomes

\begin{thmMain-bis}		
Let $K$ be a knot with unknotting number 1. 
Then one of the following three possibilities holds.
\begin{itemize}
\item[(1)] {\rm (a)}\qua
Any unknotting arc $(a,\partial a)$ for $K$ can be isotoped 
in $(S^3,K)$ so that $a\cap \bT(K) = \emptyset$.

\item[]{\rm (b)}\qua
If $\bT(K)\ne\emptyset$ and $K$ has an unknotting arc $(a,\partial a)$ 
in $\seif (K)$ then $(a,\partial a)$ 
is isotopic to an $(r,s)$--cable of an exceptional fiber of $\seif (K)$, 
for some $s\ge1$.

\item[(2)] {\rm (a)}\qua
$K$ is an EM--knot $K(\ell,m,n,p)$.

\item[]{\rm (b)}\qua $O(K)$ has a unique connected incompressible
2--sided toric 2--suborbi\-fold $S$, a Conway sphere, $K$ has an
unknotting arc $(a,\partial a)$ with $|a\cap S|=1$ (the standard
unknotting arc for $K(\ell,m,n,p)$) , and $K$ has no unknotting arc
disjoint from $S$.

\item[(3)] 
$K$ is the union of essential tangles $\P\cup \P_0$, where 
$\P_0$ is an EM--tangle $\A_\ep (\ell,m)$ and $\partial\P_0$ is in $\bT(K)$.
Any unknotting arc for $K$ can be isotoped into $\P_0$.
The standard unknotting arc for $\A_\ep (\ell,m)$ is an unknotting arc for $K$.
\end{itemize}
\end{thmMain-bis} 

Scharlemann and Thompson proved \cite{ST1,ST2} 
that if a satellite knot has unknotting number one, then an unknotting arc 
can be isotoped off any companion 2--torus.
(This follows from Corollary~3.2 of \cite{ST1} when the genus of $K$, $g(K)$, 
is $\ge 2$.
When $g(K) =1$, it follows from the proof of Corollary~3.2 of \cite{ST2}, 
or from Corollary~1 of Kobayashi \cite{Kob}, 
which say that a knot $K$ has $u(K) = g(K)=1$ if and only if it is a 
Whitehead double.)  
The following corollary of \fullref{thm:Main} can be thought of as 
a generalization of this result.

\begin{cor6.3bis}
Let $K$ be a knot with unknotting number 1, that is neither an EM--knot 
nor a knot with an EM--tangle summand with essential boundary. 
Let $F$ be an incompressible 2--sided toric 2--suborbifold of $O(K)$. 
Then any unknotting arc $(a,\partial a)$ for $K$ can be isotoped in 
$(S^3,K)$ so that $a\cap F=\emptyset$.
\end{cor6.3bis}

When a knot or link contains an essential Conway sphere, one can perform 
mutations along that sphere. 
Boileau asked \cite[Problem 1.69(c)]{Ki} if the unknotting number of a 
link is a mutation invariant. 
We prove that it is at least true for knots with unknotting number one.

\begin{thm7.1bis}
Having unknotting number 1 is invariant under mutation.
\end{thm7.1bis}

We would like to thank Alan Reid for suggesting that we consider mutation.

We apply our results to the knots that are algebraic in the sense of
Conway \cite{C} (see also Thistlethwaite \cite{Th1}), and which have
an essential Conway sphere.  We call such a knot $K$ a {\em large
algebraic} knot.  Note that the double branched cover of $K$ is a
graph manifold (ie, the union of Seifert fiber spaces identified along
their boundaries).  \fullref{thm:Main} gives particularly strong
constraints on unknotting arcs for knots in this class.

\begin{thm8.2bis}
Let $K$ be a large algebraic knot with unknotting number $1$. 
Then either

\begin{itemize}
\item[(1)] any unknotting arc for $K$ can be isotoped into either

\item[] {\rm (a)}\qua one of the rational tangles $\R(p/q)$ in an elementary tangle of 
type I; or 

\item[] {\rm (b)}\qua the rational tangle $\R(p/q)$ in an elementary tangle of type II.
\item[] 

In case (a), the crossing move transforms $\R(p/q)$ to $\R(k/1)$ for some $k$, 
and $p/q = \frac{2s^2}{2rs\pm1} +k$, where $s\ge1$ and $(r,s)=1$.

\item[] In case (b), the crossing move transforms $\R(p/q)$ to $\R(1/0)$, and 
$p/q = \frac{2rs\pm1}{2s^2}$, where $s\ge1$ and $(r,s)=1$.

\item[(2)] {\rm (a)}\qua $K$ is an EM--knot $K(\ell,m,n,p)$.

\item[]  {\rm (b)}\qua
$O(K)$ has a unique connected incompressible 2--sided toric 2--suborbifold $S$,
a Conway sphere, $K$ has an unknotting arc $a$ with $|a\cap S|=1$ (the
standard unknotting arc for $K(\ell,m,n,p)$) , and
$K$ has no unknotting arc disjoint from $S$.

\item[(3)]
$K$ is the union of essential tangles $\P\cup \P_0$, where
$\P_0$ is an EM--tangle $\A_\ep (\ell,m)$ and $\partial\P_0$ is in $T(K)$.
Any unknotting arc for $K$ can be isotoped into $\P_0$.
The standard unknotting arc for $\A_\ep (\ell,m)$ is an unknotting arc
for $K$.
\end{itemize}\end{thm8.2bis}

In Section 10 we apply \fullref{thm8.3} to the knots in Conway's 
tables \cite{C} of knots up to 11 crossings that can be immediately seen 
to be large algebraic by virtue of their description in terms of Conway's 
notation. 
There are 174 such knots, and we show that exactly 24 of them have  unknotting 
number~1. 
In particular, combining our results with those of Ozsv\'ath and Szab\'o, the 
knots with 10 or fewer crossings that have unknotting number~1 are now 
completely determined (see \cite{OS}).

It follows from \fullref{thm8.3} that the unknotting number~1 question 
is decidable for large algebraic knots.

\begin{thm11.1bis} 
There is an algorithm to decide whether or not a given large algebraic 
knot $K$, described as a union of elementary marked tangles (\fullref{alg-tangle}) 
and 4--braids in $S^2\times [0,1]$, has unknotting number~1, and, if so, 
to identify an unknotting crossing move.
\end{thm11.1bis}

We remark that the algorithm in \fullref{thm11.1} is straightforward 
to carry out in practice.

Finally, in Section 12, we consider large algebraic knots which are alternating
and show

\begin{thm11.5bis}
Let $K$ be an alternating large algebraic knot with unknotting number~$1$. 
Then $K$ can be unknotted by a crossing change in any alternating diagram 
of $K$.
\end{thm11.5bis}

The authors would like to thank Mario Eudave-Mun\~oz 
for pointing out a gap in the original proof of \fullref{D8} in the case 
that the double branched cover of $K$ is a Seifert fiber space. 
Also, the first named author wishes to acknowledge partial support 
for this work by the National Science Foundation (grant DMS-0305846).

\section{Preliminaries}

For us, a {\em tangle} will be a pair $(B,A)$ where $B$ is $S^3$ with the 
interiors of a finite number $(\ge1)$ of disjoint 3--balls removed, and $A$ 
is a disjoint union of properly embedded arcs in $B$ such that $A$ meets 
each component of $\partial B$ in four points. 
Two tangles $(B_1,A_1)$ and $(B_2,A_2)$ are {\em homeomorphic} if there 
is a homeomorphism of pairs $h\co (B_1,A_1)\to (B_2,A_2)$.

A {\em marking} of a tangle $(B,A)$ is an identification of each pair 
$(S,S\cap A)$, where $S$ is a component of $\partial B$, with 
$(S^2 , Q= \{NE,NW,SW,SE\})$. 
A {\em marked} tangle is a tangle together with a marking.
Two marked tangles are {\em equivalent} if they are homeomorphic by 
an orientation-preserving homeomorphism that preserves the markings.

A tangle $(B^3,A)$ in the 3--ball is {\em essential} if $S^2-A$ is 
incompressible in $B^3-A$.

Let $\T = (B,A)$ be a knot in $S^3$ or a tangle. 
A {\em Conway sphere} in $\T$ is a 2--sphere $S\subset \Int B$ such that 
$S$ meets $A$ transversely in four points. 
$S$ is {\em essential} if $S-A$ is incompressible in $B-A$ and $(S,S\cap A)$ 
is not pairwise parallel in $(B,A)$ to $(S_0,S_0\cap A)$ for any component 
$S_0$ of $\partial B$.

A {\em rational} tangle is a marked tangle that is homeomorphic to the
trivial tangle in the 3--ball, $(D^2,2\text{ points})\times I$.  As
marked tangles, rational tangles are parametrized by $\que \cup
\{1/0\}$.  We denote the rational tangle corresponding to $p/q \in\que
\cup \{1/0\}$ by $\re (p/q)$.  We will adopt the convention of
Eudave-Mu{\~n}oz \cite{E2} for continued fractions.  Thus
$[a_1,a_2,\ldots,a_n]$ will denote the rational number $\frac{p}q =
a_n + \frac1{a_{n-1} +\frac1{\cdots + \frac1{a_1}}}$.  We will
sometimes write $\re (p/q)$ as $\re (a_1,\ldots,a_n)$.

Let $\T = (B^3,A)$ be a tangle in the 3--ball. 
A {\em slope\/} of $\T$ is the isotopy class $(\rel\partial)$ of an 
embedded arc 
$\tau$ in $\partial B^3$ such that $\partial\tau \subset A\cap \partial B^3$. 
Given a marking on $\T$, the slopes of $\T$ are in 1--1 correspondence with 
$\que \cup \{1/0\}$  (via the double branched cover, $S^1\times S^1$, 
of $\partial B^3$ along $A\cap \partial B^3$). 
If $\T$ is rational, then $A$ defines a slope on $\partial B^3$. 
The rational number corresponding to this slope is that assigned to $\T$ 
in the preceding paragraph, $p/q$. 
If $\frac{p_1}{q_1}$, $\frac{p_2}{q_2}$ are slopes on some tangle $\T$, then 
the {\em distance\/} between these slopes, denoted $\Delta (\frac{p_1}{q_1},
\frac{p_2}{q_2})$, is $|p_1q_2 - p_2q_1|$, and 
is the minimal intersection number 
between the corresponding isotopy classes in the double branched 
cover of $\partial B^3$ along $A\cap \partial B^3$.

\begin{figure}[ht!]
\centerline{\includegraphics[height=1.0truein]{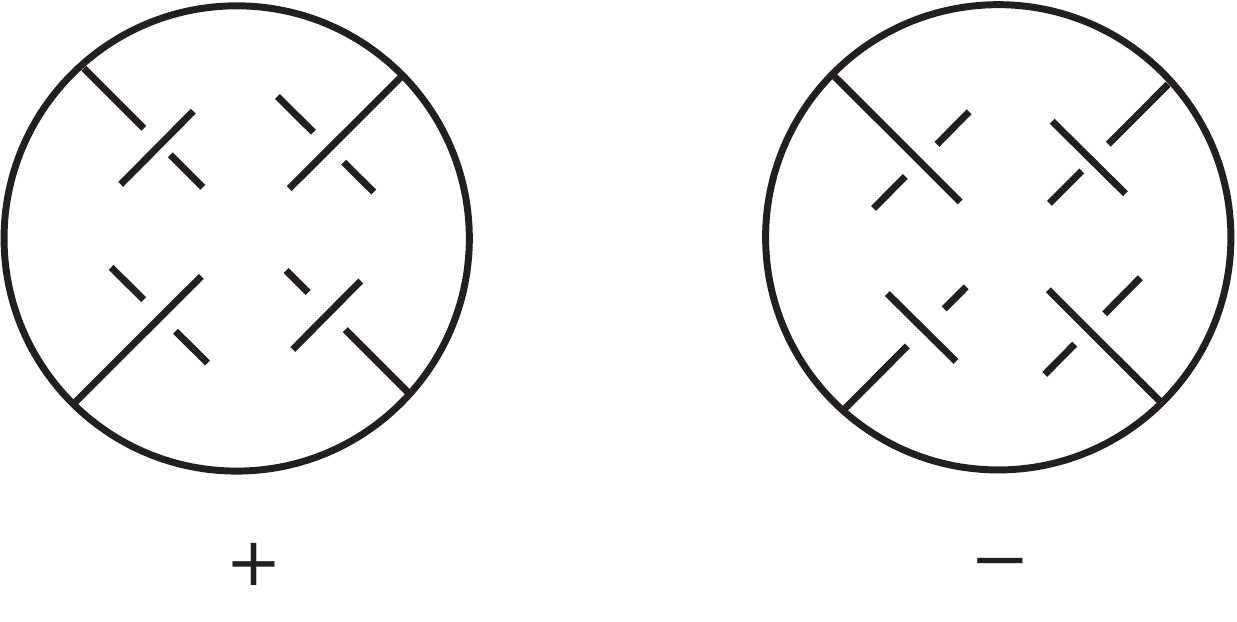}}
\caption{}
\label{Altplus}
\end{figure}

\begin{Defin}
An alternating diagram of a marked tangle 
in $B^3$ is said to be {\em positive} ({\em negative}, resp.) 
if the first crossings encountered from the boundary (with pictured marking) 
are as shown in \fullref{Altplus}. 
An alternating diagram of a marked tangle in $S^2\times I$ is said to be 
{\em positive} ({\em negative}, resp.) if filling it with $\re (1/0)$ 
gives a positive (negative, resp.) diagram of a tangle in $B^3$.
\end{Defin}

By our conventions, then, $\re (p/q)$ has a positive alternating diagram 
when $p/q >0$.

\begin{figure}[ht!]

\centerline{\includegraphics[height=.75truein]{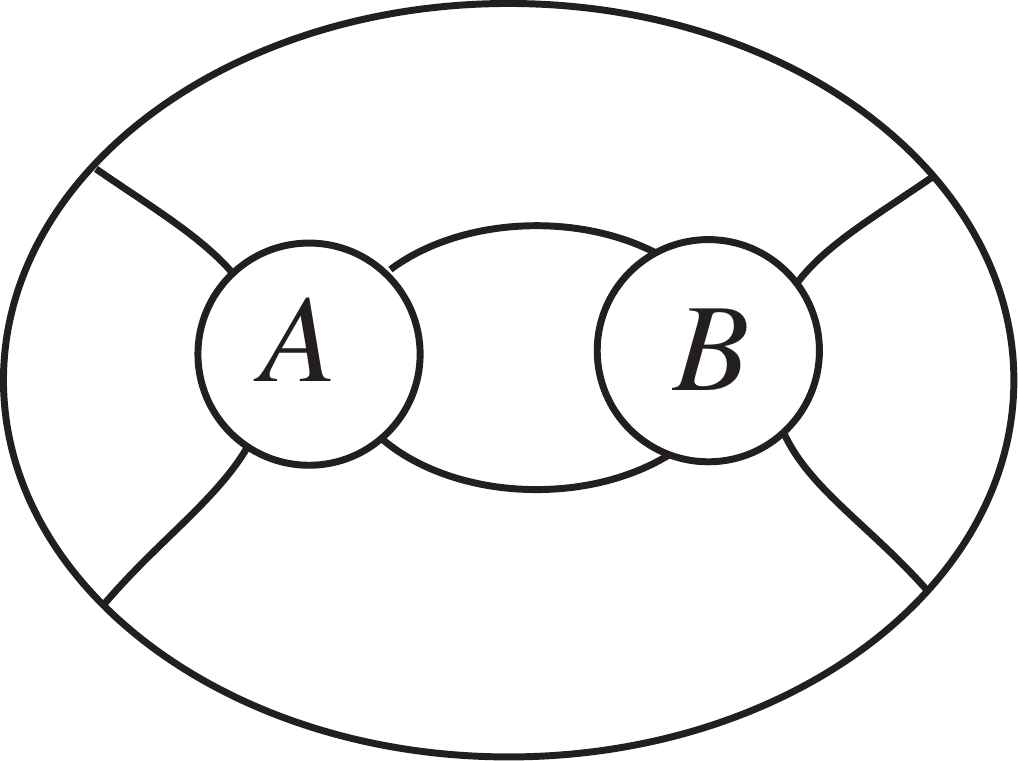}}
\caption{} 
\label{thrice-punc2}
\end{figure}

Let $\M(*,*)$ be the tangle in the thrice-punctured 3--sphere 
illustrated in \fullref{thrice-punc2}.
If $\alpha,\beta\in\que \cup \{1/0\}$, then $\M(\alpha,\beta)$ will denote 
the tangle in the 3--ball obtained by inserting rational tangles $\re (\alpha)$,
$\re(\beta)$ into $A$, $B$ respectively (with respect to the markings of 
$\partial A$ and $\partial B$ given by \fullref{thrice-punc2}). 
Similarly, $\M(\alpha,*)$ (resp.\ $\M(*,\beta)$) will denote the tangle 
in $S^2\times I$ obtained by inserting $\re (\alpha)$ (resp.\ $\re(\beta)$)
into $A$ (resp.\ $B$).

If $\alpha,\beta \in\que-\zed$ then $\M(\alpha,\beta)$ is a 
{\em Montesinos tangle of length~2}. 
Note that transferring horizontal twists between $A$ and $B$ shows that 
$\M (\alpha +m,\beta-m) = \M (\alpha,\beta)$ for all $m\in\zed$.

In general we will denote the double branched cover of a tangle 
$\T$ by $\wt{\T}$. 
However, we will denote $\wt{\M} (*,*)$ by $D^2(*,*)$;
it is homeomorphic to $P\times S^1$, where $P$ is a pair of pants. 
Similarly, denote the double branched cover of $\M(p/q,*)$ by $D^2 (p/q,*)$.
If $q>1$ this is a Seifert fiber space over the annulus with one 
exceptional fiber of multiplicity $q$. 
Finally, the double branched cover of $\M(p_1/q_1,p_2/q_2)$ is 
$D^2 (p_1/q_1,p_2/q_2)$; if $q_1,q_2>1$ this is a Seifert fiber space 
over the disk with two exceptional fibers of multiplicities $q_1$ and $q_2$.

\begin{figure}[ht!]
\centerline{\includegraphics[height=0.9truein]{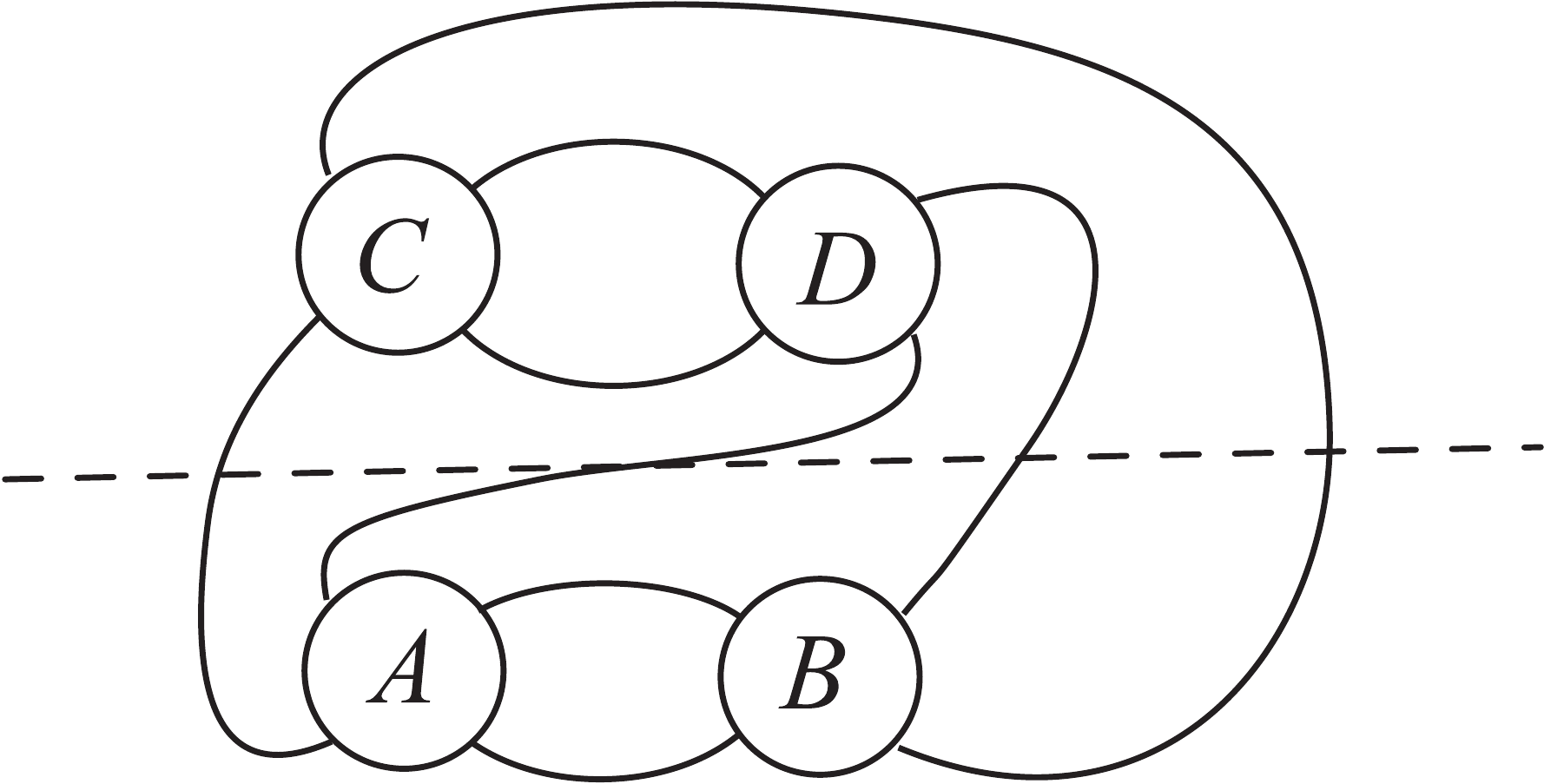}}
\caption{$\S(*,*;*,*)$}
\label{union}
\end{figure}

Let $\S (*,*;*,*)$, the {\em square} tangle, 
be the marked tangle shown in \fullref{union}; it is 
the union of two copies of $\M(*,*)$. 
If $\alpha,\beta,\gamma,\delta\in \que \cup \{\infty\}$, then 
$\S(\alpha,\beta;\gamma,\delta)$ is the knot or link obtained by 
inserting the corresponding rational tangle into $A,B,C,D$ respectively.

\begin{lem}\label{lem:rotate}  
$\S(\alpha,\beta;*,*) \cong \S(\beta,\alpha;*,*)$ by a 
homeomorphism whose restriction to $\partial C$ $(\partial D)$ is rotation 
through $180^\circ$ about the horizontal axis.
\end{lem}

\begin{proof} 
This follows by rotating \fullref{union-flipped}
through $180^\circ$ about the vertical
axis shown, using the fact that a rational tangle is unchanged by 
rotation through $180^\circ$ about the vertical axis.
\begin{figure}[ht!]
\centerline{\includegraphics[height=1.3truein]{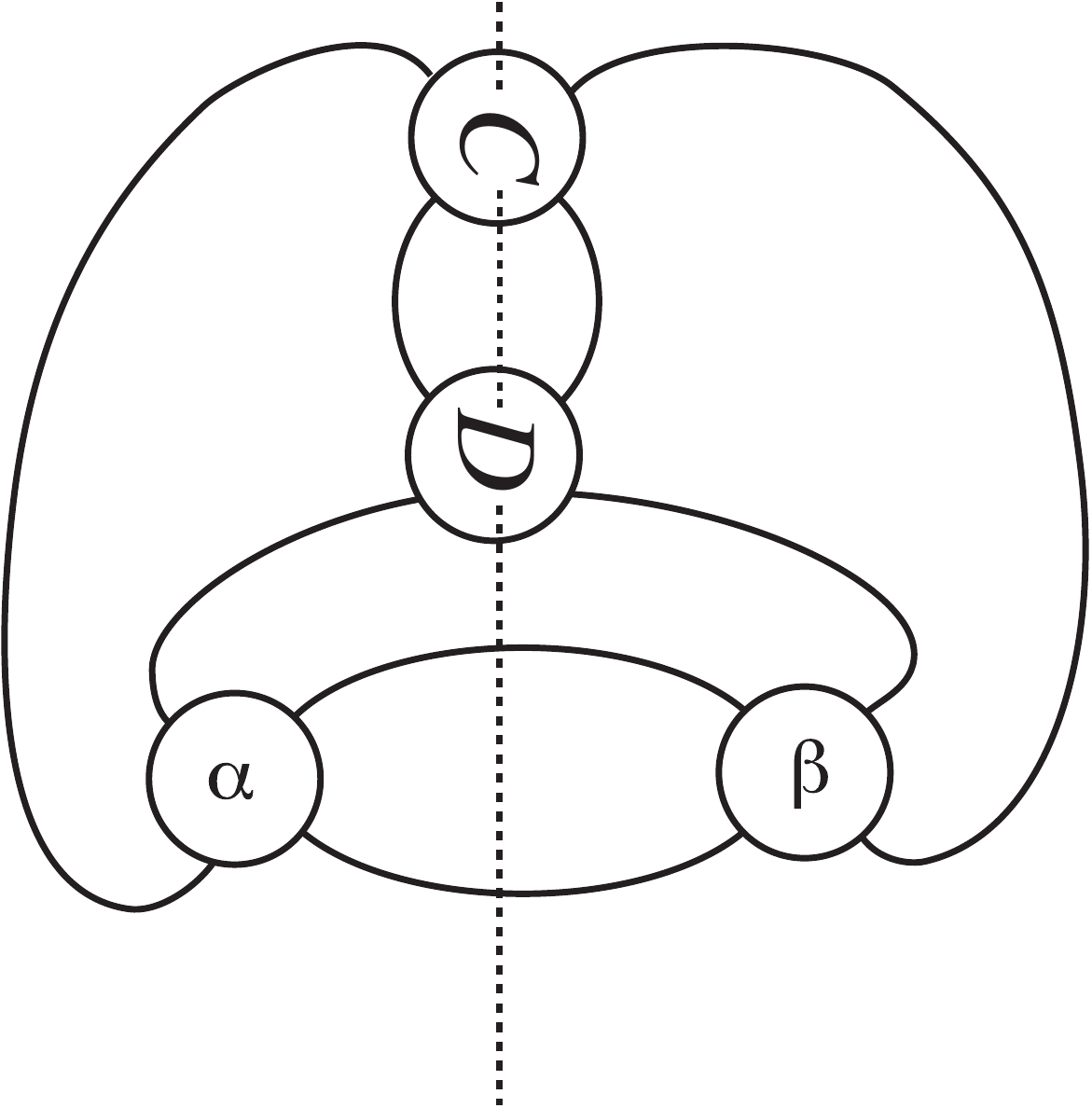}}
\caption{$\S(\alpha,\beta;*,*)$}
\label{union-flipped}
\end{figure}
\end{proof}

To state the next lemma, let $D_8$ be the order 8 dihedral group of all 
permutations of $\{\alpha,\beta,\gamma,\delta\}$ that preserve the 
partition $\{\{\alpha,\beta\},\{\gamma,\delta\}\}$. 

\begin{lem}\label{D8}\
\begin{itemize}
\item[(1)] $\S(\alpha,\beta;\gamma,\delta) = \S (\alpha +m,\beta-m; \gamma+n,
\delta -n)$ for all $m,n\in\zed$.

\item[(2)] $\S(\alpha,\beta;\gamma,\delta) = \S (\pi (\alpha),\pi(\beta);
\pi(\gamma),\pi(\delta))$ for all $\pi\in D_8$.

\item[(3)] $\S (-\alpha,-\beta;-\gamma,-\delta) = -\S(\alpha,\beta;\gamma,
\delta)$.
\end{itemize}
If $\alpha,\beta,\gamma,\delta,\alpha',\beta',\gamma',\delta' \in\que-\zed$ 
then 
\begin{itemize}
\item[(4)] $\S (\alpha,\beta;\gamma,\delta) =\S(\alpha',\beta';\gamma',
\delta')$ (resp.\ $\pm\S(\alpha',\beta';\gamma',\delta')$) 
if and only if $(\alpha,\beta$; $\gamma,\delta)$ and $(\alpha',\beta';
\gamma',\delta')$ are related by a composition of the transformations 
in (1) and (2) (resp.\ the transformations in (1), (2) and (3)). 
\end{itemize}
\end{lem}

\begin{proof}
(1) follows from the property of $\M(\alpha,\beta)$ noted earlier.

To prove (2), observe that rotating \fullref{union} 
through $180^\circ$ about 
an axis perpendicular to the plane of the paper shows that 
$\S(\alpha,\beta;\gamma,\delta) = \S(\delta,\gamma;\beta,\alpha)$.
(A rational tangle is unchanged by rotation through $180^\circ$ about any 
of the three co-ordinate axes.)
Also, by \fullref{lem:rotate}, $\S(\alpha,\beta;\gamma,\delta) 
= \S (\beta,\alpha;\gamma,\delta)$. 
The group generated by these two permutations is the dihedral group $D_8$.

(3) follows by changing all the crossings in the diagrams of 
$\R(\alpha)$, $\R(\beta)$, $\R(\gamma)$, $\R(\delta)$.

To prove (4), let $K = \S(\alpha,\beta;\gamma,\delta)$, 
$K' = \S(\alpha',\beta';\gamma',\delta')$, and let $M,M'$  be the 
double branched covers of $K$ and $K'$ respectively; thus 
$M = D^2 (\alpha,\beta) \cup D^2 (\gamma,\delta)$, and similarly for $M'$.
Parametrize slopes on the torus $T$, the double branched cover of
$S= \partial\M(\alpha,\beta)$, by the parametrization of slopes on $S$ 
coming from the marking of $\M(\alpha,\beta)$ in \fullref{union-flipped}.
Thus $1/0$ is the slope of the Seifert fiber $\phi$ of $D^2(\alpha,\beta)$,
$0/1$ is the slope of the Seifert fiber $\psi$ of $D^2 (\gamma,\delta)$, 
and similarly for $\phi',\psi'$.

Suppose $K=K'$. 
Then there is an orientation-preserving homeomorphism $h\co M\to M'$. 
Since $T$ is, up to orientation-preserving homeomorphism, 
the unique separating, incompressible torus in $M$, 
and similarly for $T'$, we may suppose that $h(T) =T'$. 
We may assume further, by interchanging $\{\alpha,\beta\}$ and 
$\{\gamma,\delta\}$ if necessary, that $h (D^2 (\alpha,\beta)) = 
D^2 (\alpha',\beta')$, and, since the Seifert fiberings of $D^2(\alpha,\beta)$ 
etc.\ are  unique, that $h(\phi)=\phi'$, $h(\psi) = \psi'$.

Recall that if $N$ is a Seifert fiber space over $D^2$ with two exceptional 
fibers, then to describe $N$ as $D^2 (\mu,\nu)$ $(\mu,\nu\in\que-\zed)$, 
we remove disjoint Seifert fibered neighborhoods of the exceptional fibers, 
getting $P\times S^1$, where $P$ is a pair of pants, and choose a section
$s\co P\to P\times S^1$. 
In identifying $\widetilde{\M} (\alpha,\beta)$ with $D^2(\alpha,\beta)$ 
we use the section that takes the boundary components of $P$ to curves 
of slope $0/1$ with respect to the markings in \fullref{union-flipped}
of $S$ and 
the boundaries of the rational tangles $\re(\alpha)$ and $\re(\beta)$.

Since $h\co D^2(\alpha,\beta)\to D^2 (\alpha',\beta')$ is an 
orientation-preserving homeomorphism which preserves the slopes $1/0$ 
and $0/1$ on $S$, 
the descriptions $D^2 (\alpha,\beta)$ and $D^2 (\alpha',\beta')$ 
differ only in the possible re-ordering of the two exceptional fibers and 
the choice of section~$s$, subject to $s(\partial_0P)$ having slope 
$0/1$, where $\partial_0P$ is the boundary component of $P$ that 
corresponds to the outer (unfilled) boundary component of $\M(*,*)$ 
in \fullref{thrice-punc2}. 
This choice corresponds to twisting a given section along an annulus
$a\times S^1\subset P\times S^1$, where $a$ is an arc in $P$ with one
endpoint in each component of $\partial P-\partial_0P$. 
This in turn corresponds to replacing $(\alpha,\beta)$ by 
$(\alpha+m,\beta-m)$ for some $m\in\zed$. 
Applying the same considerations to $D^2 (\gamma,\delta)$ and 
$D^2 (\gamma',\delta')$ gives the desired conclusion. 

The parenthetical statement in (4) now follows from (3).
\end{proof}

\begin{figure}[ht!]

\centerline{\includegraphics[height=1.0truein]{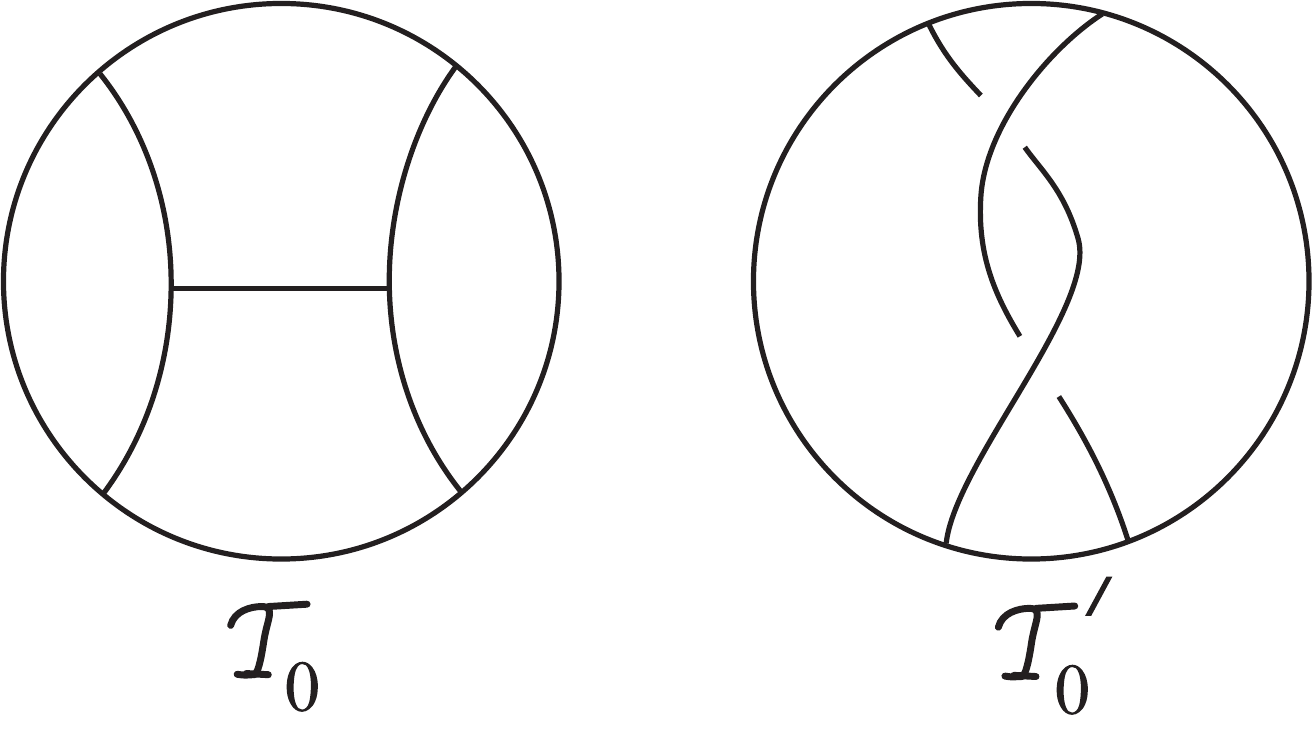}}
\caption{}		
\label{fig-main1}
\end{figure}

By a {\em crossing move} on a knot $K$ we mean the operation of passing 
one strand of $K$ through another. 
More precisely, we take a 3--ball 
$B_0$ in $S^3$ such that $(B_0,B_0\cap K)= \T_0$
is a trivial  tangle, and replace it by the trivial  tangle
$\T'_0$ shown in \fullref{fig-main1}.
This determines an arc $(a,\partial a) \subset (S^3,K)$ as shown in 
\fullref{fig-main1}.
Note that $\T_0$ is a relative regular neighborhood of $(a,\partial a)$ in 
$(S^3,K)$. 
Conversely, the arc $a$, together with a framing of $a$, determines the 
crossing move. 
If the resulting knot  
$K'$ is the unknot, we say that $a$ is an {\em unknotting arc\/} for $K$.

Note that we distinguish between a crossing {\em move} and a crossing 
{\em change}, reserving the latter term for a change of crossing in a 
knot diagram. 

If $K$ is a knot in $S^3$, then $u(K)$ is its {\em unknotting number}. 
That is, $u(K)$ is the smallest number of crossing moves required to unknot $K$.

Write $\T = (S^3,K)-\Int \T_0$, and assume that $K'$ is the unknot. 
Then taking double branched covers gives 
\begin{equation*}
\begin{split}
M &= B_2 (K) = B_2 (\T) \cup B_2 (\T_0) = X\cup V_0\ ,\\
S^3 &= B_2 (K') = B_2 (\T) \cup B_2 (\T'_0) = X\cup V'_0\ ,
\end{split}
\end{equation*}
where $V_0,V'_0$ are solid tori with meridians $\gamma,\mu$, say, on 
$\partial X$, such that $\Delta (\gamma,\mu)=2$. 
Thus the core of $V'_0$ is a knot $K^*$ in $S^3$, with exterior $X$ 
and meridian $\mu$, and $K^* (\gamma ) = X(\gamma) \cong M$.

This connection between crossing moves and Dehn surgery 
in the double branched cover is due to Montesinos \cite{Mon}. 

We now recall the characteristic toric orbifold decomposition of a knot, 
due to Bonahon and Siebenmann \cite{BS}. 

Let $K$ be a prime knot in $S^3$. 
Regard $S^3$ as an orbifold $O(K)$ with singular set $K$, each point of $K$ 
having isotropy group rotation of $\real^3$ about $\real^1$ 
through angle $\pi$. 
Since $K$ is prime, the Characteristic Toric Orbifold Splitting Theorem 
of \cite{BS} asserts the existence of a collection $\bT (O(K)) = \bT(K)$ 
of disjoint incompressible 2--sided toric 2--suborbifolds, unique up 
to orbifold isotopy, such that 
(i)~each component of $O(K)$ cut along $\bT(K)$ is either atoroidal or 
$S^1$--fibered (as an orbifold), and 
(ii)~$\bT(K)$ is minimal with respect to this property. 
(See \cite[Splitting Theorem 1]{BS}.) 
Each component of $\bT(K)$ is either a 2--torus disjoint from $K$ or 
a Conway sphere. 

$\bT(K)$ may be described as follows; see Boileau and Zimmermann
\cite{BZ}.  Let $M$ be the double branched cover of $(S^3,K)$, with
covering involution $h\co M\to M$.  Let $\bT(M)$ be the
JSJ--decomposition of $M$.  By Meeks and Scott \cite{MS}, we may
assume that $\bT(M)$ is $h$--invariant.  For each component $T$ of
$\bT(M)$ such that $h(T) =T$ and $h$ exchanges the sides of $T$,
replace $T$ by two parallel copies that are interchanged by $h$.
Denote this new collection of tori by $\bT^+(M)$.  Then $\bT(K)$ is
the quotient $\bT^+(M)/h$ in $S^3$.

\section{EM--knots and EM--tangles}	

In \cite{E1} Eudave-Mu\~noz constructed an infinite family of knots 
$K = K(\ell,m,n,p)$ such that 
(1)~$K$ has unknotting number~1, 
(2)~$K$ has a (unique) essential Conway sphere $S$, and 
(3)~no unknotting arc for $K$ is disjoint from $S$. 
Passing to double branched covers these give rise (by \cite{Mon}; see 
Section~2) to a family of hyperbolic knots $k(\ell,m,n,p)$ in $S^3$, 
called the {\em Eudave-Mu\~noz} knots in \cite{GL3}, each of which has 
a half-integral toroidal surgery. 
To distinguish the $K(\ell,m,n,p)$'s from the $k(\ell,m,n,p)$'s we 
shall call the former {\em EM--knots}. 

\begin{Defin}
The crossing move described in \cite{E1} that unknots 
$K(\ell,m,n,p)$ will be called the {\em standard crossing move} 
of $K(\ell,m,n,p)$.
\end{Defin}

Recall \cite{E1} that the parameters $\ell,m,n,p$ are restricted as 
follows: 
one of $n,p$ is always $0$; $|\ell| >1$; if $p=0$ then $m\ne0$, 
$(\ell,m)\ne (2,1)$ or $(-2,-1)$, and $(m,n)\ne (1,0)$ or $(-1,1)$; 
if $n=0$ then $m\ne 0$ or 1, and $(\ell,m,p) \ne (-2,-1,0)$ or $(2,2,1)$.

The EM--knots can be conveniently described in terms of the tangle $\S$ 
defined in Section~2
(see \fullref{union} to describe $\S$ as a marked tangle).

\begin{lem}\label{lem:EMknot} 	
The EM--knot $K(\ell,m,n,p) = \S(\alpha,\beta;\gamma,\delta)$ where 
$\alpha,\beta,\gamma,\delta$ are as follows:
\begin{equation*}
\begin{split}
p =0\ :\quad &\alpha = -\frac{1}\ell\ ,\enspace
\beta = \frac{m}{\ell m-1}\ ,\enspace
\gamma = \frac{2mn+1-m-n}{4mn-2m+1}\ ,\enspace
\delta = -\frac12\\
\noalign{\vskip6pt}
n=0\ :\quad &\alpha = -\frac1{\ell}\ ,\enspace
\beta = \frac{2mp-m-p}{\ell(2mp-m-p)-2p+1}\ ,\enspace
\gamma = \frac{m-1}{2m-1}\ ,\enspace
\delta = -\frac12\ .
\end{split}
\end{equation*}
\end{lem}

\begin{proof} 
This follows immediately from \cite[Proposition 5.4]{E2} 
(after allowing for sign errors).      
\end{proof}

\fullref{lem:EMknot}, together with the restrictions on the parameters 
$\ell,m,n,p$, easily implies

\begin{cor}\label{cor:EMknot}
Any EM--knot is of the form $\S(\alpha,\beta;\gamma,\delta)$ with 
$\alpha,\beta,\gamma,\delta\in \que-\zed$, 
$|\alpha|,|\beta|,|\gamma|,|\delta| <1$, and 
$\alpha\beta <0$, $\gamma\delta<0$.
\end{cor}

We will need to consider a collection of tangles in the 3--ball, the 
EM--tangles, closely related to the EM--knots. 
We describe them as two families, corresponding to the cases $p=0$ and 
$n=0$ of the EM--knots. 
More precisely, let $\A_1(\ell,m)$ be the tangle obtained from the 
knot $K(\ell,m,n,0)$ by removing the ``$C$''--tangle, and $\A_2(\ell,m)$ 
be the tangle obtained from $K(\ell,m,0,p)$ by removing the ``$B$''--tangle.
Here $\ell$ and $m$ are subject to the same restrictions as for 
$K(\ell,m,n,p)$, ie, $|\ell|>1$ in both cases, and for $\A_1(\ell,m)$, 
$m\ne 0$, $(\ell,m)\ne (2,1)$ or $(-2,-1)$, while for $\A_2(\ell,m)$, 
$m\ne 0$ or 1.

We therefore have the following:

\begin{Definition}\label{def:EMtangle} 
The {\em EM--tangle} $\A_\ep (\ell,m)$ is given by 
$\A_1 (\ell,m)  = \S (\alpha,\beta;*,\delta)$,\break
$\A_2 (\ell,m)  = \S (\alpha,*;\gamma,\delta)$,
where $\alpha,\beta,\gamma,\delta$ are as follows: 
\begin{equation*}
\begin{split}
\ep =1 \ :\ &\alpha = -\frac1{\ell}\ ,\quad 
\beta = \frac{m}{\ell m-1}\ ,\quad 
\delta = -\frac12\ ;\ |\ell| >1,\ m\ne 0,\\
&\hskip2.8in (\ell,m) \notin  \{(2,1),(-2,-1)\}\\
\ep =2 \ :\ &\alpha = -\frac1{\ell}\ ,\quad 
\gamma = \frac{m-1}{2m-1}\ ,\quad 
\delta = -\frac12\ ;\ |\ell| >1,\ m\notin \{0,1\} .
\end{split}
\end{equation*}
\end{Definition}
\noindent
(The $\A_\ep (\ell,m)$ are pictured in \fullref{Decomp} where the 
twist boxes represent vertical twists.) 

$K(\ell,m,n,p)$ contains an essential Conway sphere $S$, decomposing it into 
two Montesinos tangles: 
$K(\ell,m,n,p) = \M(\alpha,\beta)\cup_S \M(\gamma,\delta)$. 
This gives rise to a decomposition of the double branched cover of 
$K(\ell,m,n,p)$ as $N_1\cup_T N_2$, where $N_i$ is a Seifert fiber space 
over the disk with two exceptional fibers, $i=1,2$, and $T=\partial N_1 
= \partial N_2 = \widetilde S$ is the double branched cover of $S$. 
Similarly, the essential Conway sphere $S$ in $\A_\ep (\ell,m)$ gives 
a decomposition of its double branched cover as $N_1\cup_T N_2$, where 
$N_2$ is as above, and $N_1$ is a Seifert fiber space over the annulus 
with one exceptional fiber.

The remainder of this section is devoted to proving the 
following theorem which says that the EM--knots and the EM--tangles are 
determined by their double branched covers. 

\begin{thm}\label{thm:doublebranched}\
\begin{itemize}
\item[{\rm (1)}] Let $K$ be a knot in $S^3$ whose double branched cover is 
homeomorphic to that of $K(\ell,m,n,p)$. 
Then $K= \pm K(\ell,m,n,p)$.

\item[{\rm (2)}] Let $\T$ be a tangle in $B^3$ whose double branched
cover is homeomorphic to that of $\A_\ep (\ell,m)$.  Then $\T$ and
$\A_\ep (\ell,m)$ are homeomorphic tangles.
\end{itemize}\end{thm}

In order to prove \fullref{thm:doublebranched}, we first study 
involutions on the manifolds $D^2 (*,*)$, $D^2(p/q,*)$ and 
$D^2(p_1/q_1,p_2/q_2)$. 
The definition of equivalence that is appropriate to our purposes is the 
following. 
Two homeomorphisms $f,g\co X\to X$ are {\em strongly conjugate} if there is 
a homeomorphism $h\co X\to X$ isotopic to the identity such that $f=h^{-1}gh$.
If $X_0\subset X$, then $f$ and $g$ are strongly conjugate {\em rel $X_0$}
if $h$ can be chosen to be isotopic to the  identity by an isotopy fixed 
on $X_0$.
The set of fixed points of an involution $\tau$ will be denoted 
by $\Fix (\tau)$.

To define a standard model for a pair of pants $P$, let $D^2$ be the 
unit disk in $\real^2$, let $D_1$ and $D_2$ be disjoint round disks in 
$\Int D^2$ with their centers on the $x$--axis, and let 
$P = D^2-\Int (D_1\cup D_2)$.
The map $(x,y)\mapsto (x,-y)$ defines an orientation-reversing involution 
$\rho_P$ of $P$, which we will call {\em reflection}.

\begin{lem}\label{lem:reflection}
Let $\tau$ be a non-trivial involution on a pair of pants $P$. 
If each boundary component of $P$ is invariant under $\tau$ then $\tau$ 
is strongly conjugate to reflection.
\end{lem}
\vspace{-4pt}

\begin{proof} 
Let $\alpha = P\cap (x\text{--axis})$, a disjoint union of three arcs properly 
embedded in $P$. 
Since the restriction of $\tau$ to each boundary component of $P$ is either
the identity, conjugate to rotation through $\pi$, or conjugate to 
reflection, we may assume that $\partial\alpha$ is invariant under $\tau$. 
By analyzing the intersections of $\alpha$ and $\tau (\alpha)$, one can 
show that after conjugating $\tau$, $\alpha$ can be taken to be invariant 
under $\tau$. 
This implies that $\alpha$ is fixed by $\tau$. 
The two disks of $P-\alpha$ are either exchanged or invariant. 
In the first case, $\tau$ is conjugate to reflection, and in the second 
$\tau$ must be the identity.
Finally, since any homeomorphism of $P$ is isotopic to one that commutes 
with $\rho_P$, in the first case $\tau$ is strongly conjugate to $\rho_P$.
\end{proof}
\vspace{-4pt}

Recall (Section 2) that $D^2 (*,*)$ (resp.\ $D^2 (p/q,*)$, resp.\ 
$D^2 (p_1/q_1,p_2/q_2)$) is the double branched cover of the tangle 
$\M(*,*)$ (resp.\ $\M(p/q,*)$, resp. $\M(p_1/q_1,\allowbreak p_2/q_2)$). 
The {\em standard involution} on $D^2(*,*)$, $D^2 (p/q,*)$ or 
$D^2 (p_1/q_1,p_2/q_2)$ is the non-trivial covering transformation 
corresponding to this double branched cover. 
Note that in particular, identifying $D^2(*,*)$ with $P\times S^1$, 
the standard involution on $D^2(*,*)$ is the map $(x,\theta) \mapsto 
(\rho_P(x),-\theta)$.
\vspace{-4pt}

\begin{lem}\label{lem:non-trivial}
Let $\tau$ be a non-trivial orientation-preserving involution on $P\times S^1$.
If each component of the boundary is invariant under $\tau$ then $\tau$ is 
strongly conjugate to either the standard involution or a free involution 
that leaves each $S^1$--fiber invariant.
\end{lem}
\vspace{-4pt}

\begin{proof}
By Tollefson \cite{To}, there is a Seifert fibration of $P\times S^1$ that is 
invariant under $\tau$. 
As the Seifert fibration is unique up to isotopy, we may therefore assume, 
after strongly conjugating $\tau$, that $\tau$ preserves the product 
$S^1$--fibration of $P\times S^1$. 
Thus $\tau$ induces an involution $\tau_P$ on $P$. 
By \fullref{lem:reflection}, $\tau_P$ is either the identity or 
strongly conjugate to $\rho_P$. 
\vspace{-4pt}

If $\tau_P$ is the identity, then $\tau$ takes each $S^1$--fiber to 
itself by an orientation-preserving involution, hence by either the 
identity or a map conjugate to rotation through $\pi$. 
By continuity, the action is the same on each fiber. 
Therefore $\tau$ is either the identity or free.
\vspace{-4pt}

We may suppose, then, that $\tau_P$ is strongly conjugate to $\rho_P$, 
and hence, by strongly conjugating $\tau$, that $\tau_P = \rho_P$. 
\vspace{-4pt}

Let $C$ be a boundary component of $P$.  Note that $\tau_P$ fixes two
points in $C$.  The restriction of $\tau$ to each of the two
corresponding $S^1$--fibers is therefore conjugate to reflection
$\theta\mapsto -\theta$.  It follows that the restriction of $\tau$ to
$C\times S^1$ is strongly conjugate to $-I$, given by
$(\varphi,\theta)\mapsto (-\varphi,-\theta)$ (see Hartley \cite{Har}).
So we may assume that $\tau =-I$ on each boundary component of
$P\times S^1$.  Since any two $S^1$--fibrations of $P\times S^1$ that
agree on the boundary are isotopic $\rel \partial$, we can still
assume that $\tau$ preserves the product fibration and that $\tau_P =
\rho_P$.

Let $\alpha_1,\alpha_2$ be the two arc components of $\Fix(\rho_P)$ 
shown in \fullref{Fix-rho}, and let $C_1,C_2$ be the two boundary components 
of $P$ indicated in the same figure. 

\begin{figure}[ht!]
\centerline{\includegraphics[height=0.9truein]{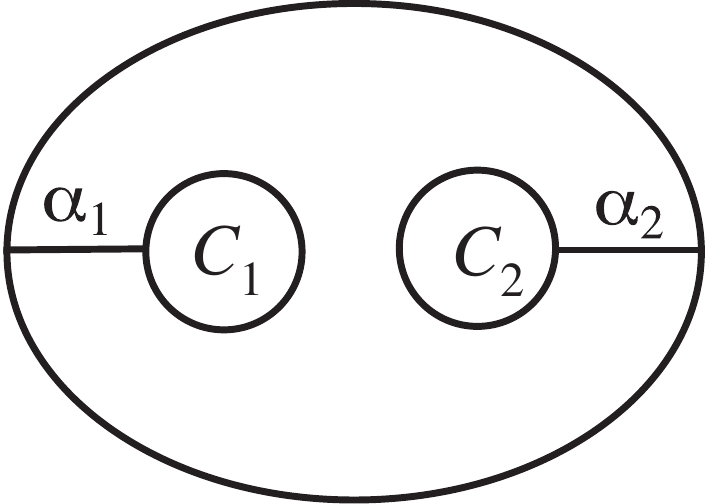}}
\caption{}			
\label{Fix-rho}
\end{figure}

Let $A_i$ be the vertical annulus $\alpha_i\times S^1$, and let $T_i$ 
be the boundary torus $C_i\times S^1$, $i=1,2$. 
Note that $A_i$ is invariant under $\tau$, $i=1,2$. 
Then $\tau | (A_i= \alpha_i\times S^1)$ is conjugate to the involution 
$(x,\theta)\mapsto (x,-\theta)$, by a homeomorphism $g_i\co A_i\to A_i$ 
that is isotopic $\rel\partial$ to a power of a Dehn twist along the core 
of $A_i$.
Hence, conjugating $\tau$ by the corresponding power of a vertical Dehn 
twist $h_i$ along a torus in a collar neighborhood of $T_i$, we may assume 
that $\tau | A_i$ is $(x,\theta) \mapsto (x,-\theta)$. 
Since $h_i$ is isotopic to the identity, the strong conjugacy class 
of $\tau$ is unchanged. 

By a further isotopy $\rel\partial$, we may assume that $A_i = A_i\times\{0\}$ 
has a neighborhood $A_i\times [-1,1]$ on which $\tau$ acts by 
$(x,\theta,t)\mapsto (x,-\theta,-t)$. 
Removing $(A_1\cup A_2)\times (-1,1)$ from $P\times S^1$, $\tau$ induces 
an involution $\tau_0$ on $D^2\times S^1$ which is equal to $-I$ on the 
boundary. 
Hence $\tau_0$ is strongly conjugate $\rel\partial$ to the involution 
$((x,y),\theta)\mapsto ((x,-y),-\theta)$ \cite{Har}. 
Reattaching $(A_1\cup A_2)\times [-1,1]$ we get that $\tau$ is strongly 
conjugate to the standard involution.
\end{proof}

\begin{lem}\label{lem:Fix-tau}
Let $\tau$ be a non-trivial orientation-preserving involution 
on\break $D^2 (p_1/q_1,p_2/q_2)$, where $q_1\ne q_2$. 
\begin{itemize}
\item[(1)] If $\Fix(\tau)$ has non-empty intersection with the boundary 
then $\tau$ is strongly conjugate to the standard involution.
\item[(2)] If $\tau$ acts freely on the boundary then the Seifert fibration 
of $D^2 (p_1/q_1,p_2/q_2)$ may be isotoped so that $ \tau$ leaves each 
Seifert fiber on the boundary invariant.
\end{itemize}
\end{lem}

\begin{proof}
By \cite{To}, $D^2 (p_1/q_1,p_2/q_2)$ has a Seifert fibration for which 
$\tau$ is fiber-pres\-erving.
Since the Seifert fibration is unique up to isotopy, and since $q_1\ne q_2$, 
the exceptional fibers must be invariant, and hence they have disjoint 
invariant fibered neighborhoods $V_1$ and $V_2$, say. 
Thus $\tau$ restricts to an involution $\tau_0$ on the complement 
of these neighborhoods, $D^2(*,*) = P\times S^1$. 
Note that $\tau_0$ leaves each boundary component of $D^2(*,*)$ invariant.

(1)\qua Here $\Fix (\tau_0) \ne\emptyset$, so by \fullref{lem:non-trivial}
$\tau_0$ is strongly conjugate to the standard involution. 
Extending over $V_1$ and $V_2$, and using 
\cite{Har}, we get that $\tau$ is strongly conjugate to the standard 
involution.

(2)\qua Since $\tau_0$ acts freely on at least one of the boundary components 
of $D^2 (*,*)$, by \fullref{lem:non-trivial} the (product) Seifert 
fibration of $D^2(*,*)$ may be isotoped so that $\tau$ leaves each fiber 
invariant. 
Now $\tau|\partial V_i$ can be extended to an involution $\tau_i$ of 
$V_i$ that leaves each Seifert fiber invariant. 
Since $\tau_i$ and $\tau |V_i$ agree on $\partial V_i$, they are strongly 
conjugate $\rel \partial V_i$ \cite{Har}. 
The corresponding isotopy of $V_i$ ($\rel \partial V_i$) takes the Seifert 
fibration of $V_i$ to one such that each fiber is invariant under $\tau$.
\end{proof}

\begin{lem}\label{lem:Fix-tau2}
Let $\tau$ be a non-trivial orientation-preserving involution 
on $D^2 (p/q,*)$ such that $\Fix (\tau)$ has non-empty intersection with the 
boundary. 
Then $\tau$ is strongly conjugate to the standard involution.
\end{lem}

\begin{proof}
This is the same as the proof of Part (1) of \fullref{lem:Fix-tau}.
\end{proof}

\begin{proof}[Proof of \fullref{thm:doublebranched} Case (1)]
Write $K_0 = K(\ell,m,n,p)$. 
Then $(S^3,K_0)\allowbreak 
= (B_1,A_1) \cup_S (B_2,A_2)$, where $S$ is an essential 
Conway sphere and $(B_i,A_i)$ is a Montesinos tangle of length~2, $i=1,2$. 
The double branched cover of $(S^3,K_0)$ is $N= N_1\cup_{\widetilde S}N_2$, 
where $N_i$, the double branched cover of $(B_i,A_i)$, is a Seifert 
fiber space over the disk with two exceptional fibers, $i=1,2$, and 
$\widetilde S = \partial N_1 = \partial N_2$ is the double branched 
cover of $(S,S\cap K_0)$. 
The Seifert fibers of $N_1$ and $N_2$ intersect once on $\widetilde S$.
The covering involution $\sigma \co N\to N$ restricts to the standard 
involution $\sigma_i$ on $N_i$, $i=1,2$.

Now suppose $K$ is a knot in $S^3$ whose double branched cover is 
homeomorphic to $N$. 
Let $\tau \co N\to N$ be the corresponding covering involution. 
Since $\widetilde S$ is the unique incompressible torus in $N$, up to 
isotopy, by \cite[Theorem~8.6]{MS} we may assume that $\widetilde S$ 
is invariant under~$\tau$.

\begin{claim}\label{claim1}	
Each $N_i$ is invariant under $\tau$.
\end{claim}

\begin{proof}
If $\tau$ interchanges $N_1$ and $N_2$ then $\Fix (\tau)$ is contained
in $\widetilde S$.
With respect to some parametrization of $\widetilde S$ as $S^1\times S^1$, 
$\Fix (\tau)$ is a (2,1)--curve and $\tau$ leaves each (0,1)--curve 
$\gamma$ invariant, taking it to itself by reflection in the pair of 
points $\Fix (\tau)\cap\gamma$ (thus the quotient $\widetilde S/\tau$ 
is a M\"obius band). 
Then $S^3  = N_1 /(\tau|\widetilde S)$ is homeomorphic to $N_1$ with 
a solid torus $V$ attached so that $\gamma$ bounds a meridian disk of $V$.
Hence $N_1$ is a knot exterior with meridian $\gamma$. 
Applying the same argument to $N_2$, we see that in $N= N_1\cup N_2$ 
the meridians of $N_1$ and $N_2$ are identified. 
But this is not true: when each side of $N$ is the
exterior of a knot in $S^3$, the argument in Lemma~1.3 of \cite{E1} 
(or Lemmas~\ref{lem:EMknot} and \ref{newlem:9.5})
shows that the meridian of one side is 
identified with the Seifert fiber of the other.
\end{proof}
\vspace{-3pt}

Let $\tau_i$ be the restriction of $\tau$ to $N_i$, $i=1,2$.
\vspace{-3pt}

\begin{claim}\label{claim2}
$\tau_i$ is strongly conjugate to the standard involution $\sigma_i$ 
on $N_i$, $i=1,2$. 
\end{claim}
\vspace{-3pt}

\begin{proof}
If $\Fix (\tau)$ meets $\widetilde S$ then the result follows 
from \fullref{lem:Fix-tau}(1) and \fullref{lem:doublebranched} below.
\vspace{-3pt}

If $\Fix (\tau)$ is disjoint from $\widetilde S$, then by 
\fullref{lem:Fix-tau}(2) and \fullref{lem:doublebranched} the Seifert 
fibrations of $N_1$ and $N_2$ can be isotoped so that, on $\widetilde S$, 
each $S^1$--fiber of each fibration is invariant under $\tau$. 
But since the fibers of the two fibrations intersect once on $\widetilde S$,
this is clearly impossible.
\end{proof}
\vspace{-3pt}

Write $S_i = \partial B_i$, $i=1,2$, and let $f\co (S_1,S_1\cap A_1)\to 
(S_2,S_2\cap A_2)$ be the gluing homeomorphism that defines 
$(S^3,K_0) = (B_1,A_1)\cup_f (B_2,A_2)$. 
To compare $K_0$ and $K$, we need the notion of a mutation involution, which 
is defined at the beginning of Section~7.
\vspace{-3pt}

\begin{claim}\label{claim3}
$(S^3,K)$ is homeomorphic to $(B_1,A_1)\cup_{\mu f}(B_2,A_2)$ for some 
mutation involution $\mu$ of $(S_2,S_2\cap A_2)$. 
\end{claim}
\vspace{-3pt}

\begin{proof}
Let $\tilde f\co \partial N_1\to \partial N_2$ be a lift of $f$, giving 
$N= N_1 \cup_{\tilde f} N_2$.
Note that $\tilde f\sigma_1 = \sigma_2 \tilde f$ and $\tilde f\tau_1 = 
\tau_2\tilde f$. 
Also, $N/\tau = (N_1/\tau_1)\cup_g (N_2/\tau_2)$ for some $g\co  
\partial (N_1/\tau_1) \to \partial (N_2/\tau_2)$ such that $\tilde f$ 
is a lift of $g$.
\vspace{-3pt}

By Claim~\ref{claim2}, there is a homeomorphism $\tilde h_i \co N_i\to N_i$, 
isotopic to the identity, such that $\tau_i = \tilde h_i^{-1}\sigma_i
\tilde h_i$, $i=1,2$. 
Then $\tilde h_i$ induces a homeomorphism $h_i\co N_i/\tau_i \to N_i/\sigma_i
= (B_i,A_i)$, $i=1,2$. 
Let $\partial h_i$ be the restriction of $h_i$ to $\partial (N_i/\tau_i)$. 
Then $h_1\cup h_2$ induces a homeomorphism $h\co N/\tau = (N_1/\tau_1) 
\cup_g (N_2/\tau_2) \to (B_1,A_1) \cup_e (B_2,A_2)$, where 
$e = (\partial h_2)^{-1} g(\partial h_1)$. 
Then $e$ lifts to $\tilde e = (\partial\tilde h_2)^{-1}\tilde f(\partial 
\tilde h_1)$, which is isotopic to $\tilde f$.
\vspace{-3pt}

Let $\mu$ be the mutation involution of $(S_2,S_2\cap A_2)$ such that the 
composition $\mu f$ agrees with $e$ on some point of $S_1\cap A_1$. 
Since $\tilde e$ is isotopic to $\tilde f$, $e$ and $f$ induce the same 
function from the set of (unoriented) isotopy classes of essential 
simple closed curves in $S_1-  (S_1\cap A_1)$ to the set of those in 
$S_2-(S_2\cap A_2)$. 
Hence $e$ and $\mu f$ do also. 
Since $e$ and $\mu f$ agree on a point of $S_1\cap A_1$, they must be 
isotopic as maps of pairs. 
Therefore $(S^3,K) \cong (B_1,A_1)\cup_e (B_2,A_2) \cong 
(B_1,A_1) \cup_{\mu f} (B_2,A_2)$.
\end{proof}
\vspace{-3pt}

By Claim~\ref{claim3}, $K$ is a mutation of $K_0$ along $S$. 
By \cite{E1} such mutations yield $K_0$ again. 
\vspace{-3pt}

This completes the proof of the theorem in Case (1). 

{\bf Case (2)}\qua
We have $\A_\ep(\ell,m) = (B_1,A_1)\cup_S (B_2,A_2)$, where $B_1$ 
is $S^2\times I$, $(B_1,A_1) = \M(p/q,*)$, and $(B_2,A_2)$ is a 
Montesinos tangle of length~2 as in Case~(1); 
see \fullref{def:EMtangle}. 
The double branched cover of $\A_\ep (\ell,m)$ is then 
$N = N_1\cup_{\widetilde S} N_2$, where $N_1$ is now a Seifert fiber 
space over the annulus with one exceptional fiber.
The covering involution $\sigma$ restricts to the standard involution
$\sigma_i$ on $N_i$, $i=1,2$.

Suppose $\T$ is a tangle in $B^3$ whose double branched cover is 
homeomorphic to $N$, and let $\tau$ be the corresponding covering  involution.
Since $\widetilde S$ is the unique essential torus in $N$, we may assume 
$\tau (\widetilde S)=\widetilde S$. 
Let $\tau_i$ be the restriction of $\tau$ to $N_i$, $i=1,2$.
By \fullref{lem:Fix-tau2}, $\tau_1$ is strongly conjugate to the 
standard involution on $N_1$. 
In particular $\Fix (\tau)\cap \widetilde S\ne \emptyset$, and so, by 
Lemmas~\ref{lem:Fix-tau}(1) and \ref{lem:doublebranched}, $\tau_2$ is 
strongly conjugate to the standard involution on $N_2$. 
Now Claim~\ref{claim3} holds, exactly as in Case~(1). 
Since $\A_\ep(\ell,m)$ is unchanged by mutation along $S$, 
(by \fullref{lem:rotate} and the fact that rotating a rational tangle 
through $\pi$ about a co-ordinate axis does 
not change it), $\T$ is homeomorphic to $\A_\ep (\ell,m)$.
\end{proof}

\begin{lem}\label{lem:doublebranched}
Let $N=N_1\cup N_2$ be the double branched cover of $S^3$ over 
$K(\ell,m,n,p)$ or  $\A_\ep (\ell,m)$, as in the discussion before the
statement of \fullref{thm:doublebranched}. 
Let $q_1,q_2$ be the orders of the two exceptional fibers of $N_i$ 
for some $i$. 
Then $q_1\ne q_2$.
\end{lem}

\begin{proof}
When $N$ is the cover of an EM--tangle, the lemma applies to $N_2$. 
By \fullref{def:EMtangle} 
the exceptional fibers are of orders $(|\ell|, |1-\ell m|)$ with 
$|\ell|>1$, $m \ne 0$ or of orders $(2,|2m-1|)$ where $m \ne 0,1$. In either
case, the lemma easily follows.

So assume $N$ is the double branched cover of $K(\ell\!,m\!,n\!,p)$,
and recall \fullref{lem:EMknot}. 

If $p=0$, the orders of the exceptional fibers are 
$(|\ell|,|1-\ell m|)$ for $N_1$ and $(2,|4mn-2m+1|)$ for $N_2$. 
In this case the lemma is clear.

So we assume $n=0$. In this case the exceptional fibers have orders
$(2,|2m-1|)$ in $N_1$ and $(|\ell|,|2\ell mp-\ell m-\ell p-2p+1|)$ in $N_2$. 
We need to consider the solutions to 
\begin{align}
\pm l &= 2\ell mp-\ell m-\ell p-2p+1 \notag\\
\Longleftrightarrow \quad \qquad 2p-1 &= \ell (2mp-m-p \pm 1) \label{star}\\
\Longleftrightarrow  \qquad 2(2p-1) &= \ell(4mp-2m-2p \pm 2) \notag\\
\Longleftrightarrow  \qquad 2(2p-1) &= \ell((2m-1)(2p-1) + a)\label{starstar} 
\end{align}
where $a$ is $-3$ or $1$. When $n=0$, we assume that $m \neq 0,1$. 
Thus $|2m-1| \geq 3$. 
Then
$$ |(2m-1)(2p-1) + a| \geq |(2m-1) (2p-1)| - 3 \geq 2|2p-1| + |2p-1| - 3\ .$$
Since $|\ell| \geq 2$, \eqref{star}  implies that $|2p-1| \geq 3$. 
Thus \eqref{starstar}  becomes
$$ |2(2p-1)| \geq |\ell | |2(2p-1)|\ ,$$ a contradiction.
\end{proof}

\section{Knots in solid tori}	

In this section we describe the hyperbolic knots in a solid torus that have 
a non-integral toroidal Dehn surgery. 
We gave a description of such knots in \cite[Corollary A.2]{GL3}; here we 
sharpen this to a complete characterization.

\begin{figure}[ht!]
\centerline{\includegraphics[height=1truein]{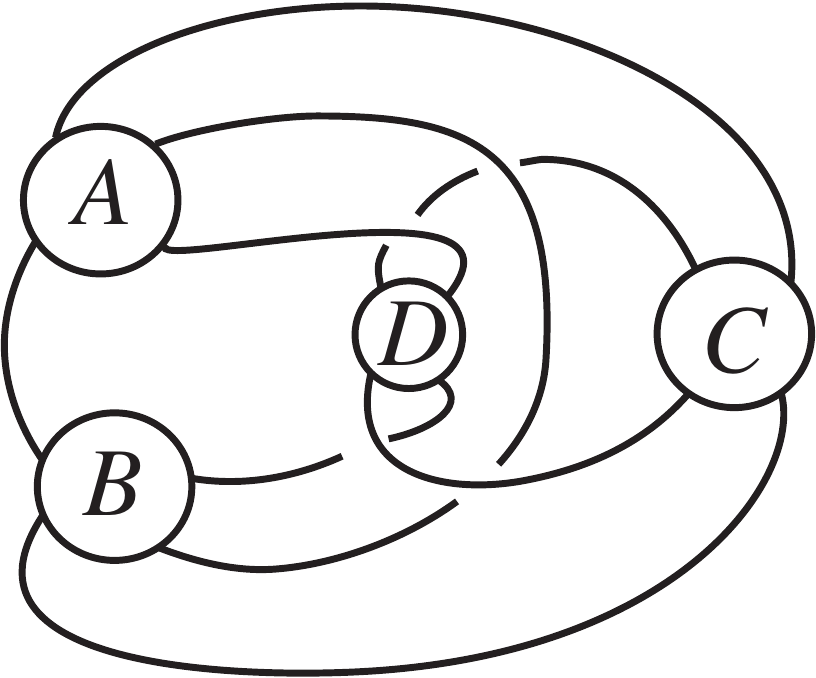}}
\caption{$\C(A,B,C,D)$}
\label{solid-tori1}
\end{figure}

\begin{figure}[ht!]
\centerline{\includegraphics[height=2.5truein]{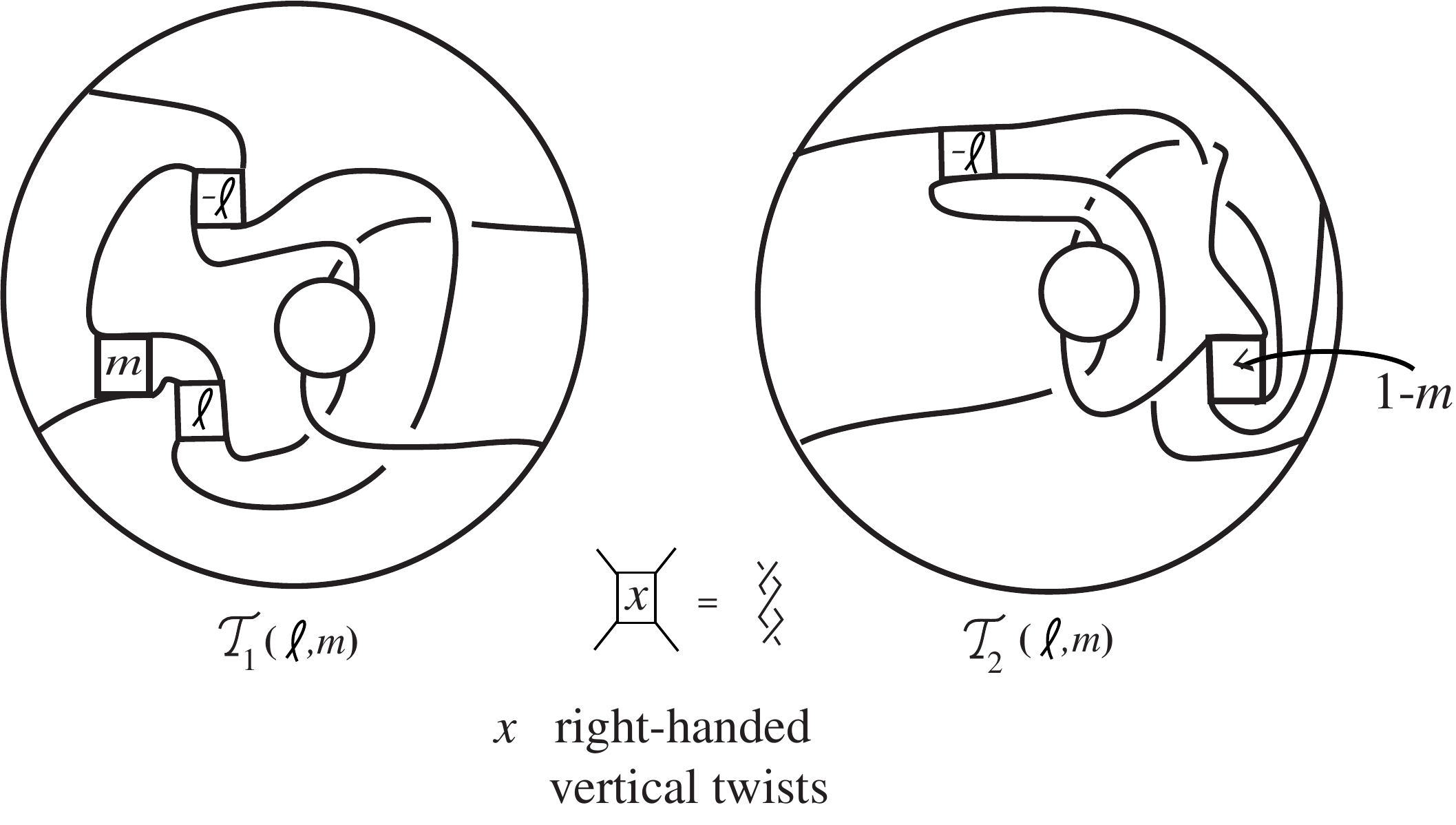}}
\caption{}
\label{ScriptT}
\end{figure}

The exteriors of the knots are the double branched covers of certain 
tangles in $S^2\times I$, which we now describe. 
Let $\C(A,B,C,D)$ be the tangle shown in \fullref{solid-tori1}. 
(To be consistent with the notation of \cite{E2}, we here regard $A,B,C,D$ as
denoting rational tangles; a puncture that is not filled in will as usual be 
indicated by a $*$.)
Note that in the terminology of \cite{GL3}, $\C(A,B,C,D) = \B(A,B,C)+R(D)
= \P(A,B,C,\frac12,D)$. 
Define tangles $\T_1(\ell,m)$ and $\T_2(\ell,m)$ in $S^2\times I$ as follows. 
\begin{itemize}
\item[{}] $\T_1(\ell,m) = \C(A,B,*,*)$, where $A= \R (\ell)$, 
$B= \R (m,-\ell)$, and $\ell,m$ are integers such that $|\ell|>1$, $m\ne 0$, 
and $(\ell,m)\ne (2,1)$ or $(-2,-1)$. 
\item[{}] $\T_2(\ell,m) = \C(A,*,C,*)$, where $A= \R (\ell)$, 
$C= \R (m-1,2,0)$, and $\ell,m$ are integers such that 
$|\ell| >1$, $m\ne 0,1$.
\end{itemize}
See \fullref{ScriptT} (the boxes correspond to vertical twists).

Note that $\T_1(\ell,m)$ is obtained by removing the ``$C$--tangle" from 
$B(\ell,m,n,0)$ of \cite{E1}, and similarly $\T_2(\ell,m)$ is obtained by 
removing the ``$B$--tangle" from $B(\ell,m,0,p)$. 

Filling the ``$D$''--puncture of $\T_1(\ell,m)$ or $\T_2(\ell,m)$ 
with the $\frac10$--tangle, ie, $\C(A,B$, $*,\frac10)$ or $\C(A, *,C,\frac10)$
with $A,B,C$ as indicated above, gives a rational tangle. 
So, the corresponding $\frac10$--Dehn filling of the double branched cover 
of $\T_\ep(\ell,m)$, $\ep\in \{1,2\}$, is a solid torus. 
Let $J_\ep (\ell,m)$ denote the core of this Dehn filling, seen as a 
knot in this solid torus.
\vspace{2pt}

\begin{figure}[ht!]

\centerline{\includegraphics[height=1.7truein]{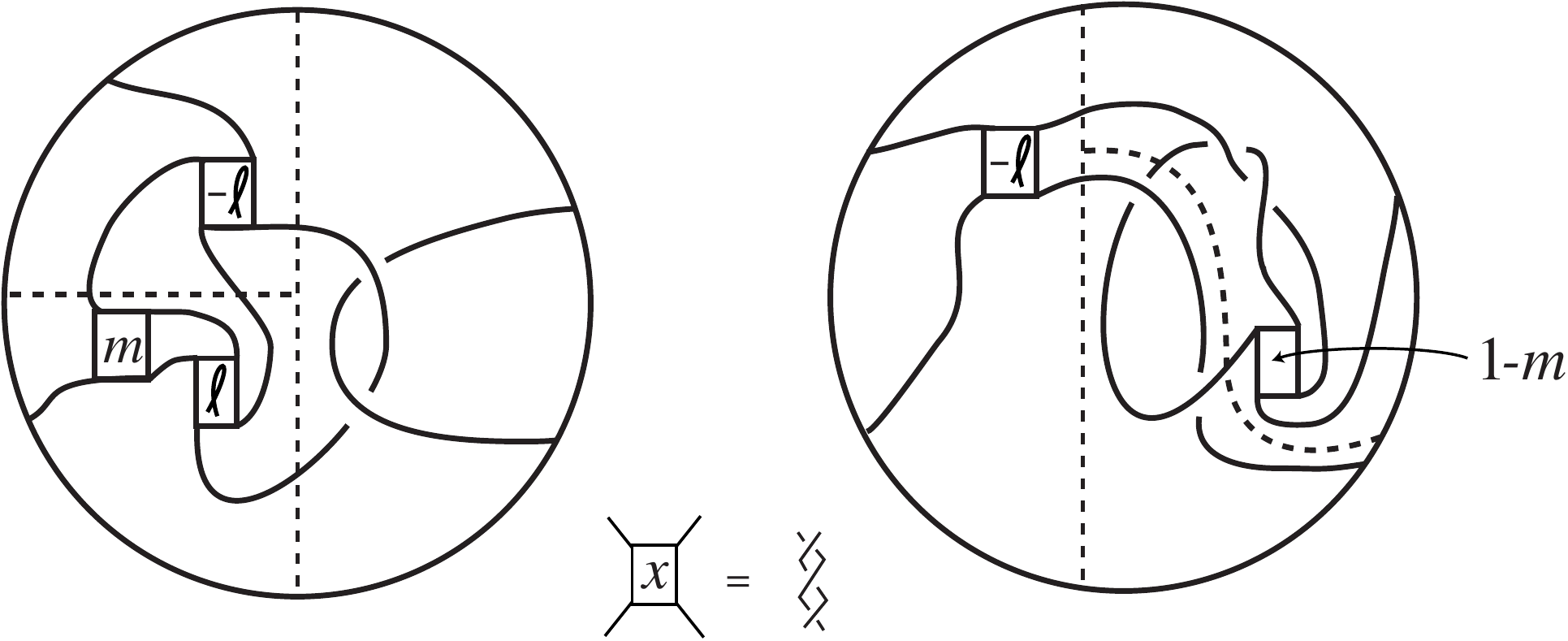}}
\centerline{\small$\A_1(\ell,m)$\hskip1.7truein $\A_2 (\ell,m)$\qquad}	
\caption{} 
\label{Decomp}
\end{figure}

Note that the EM--tangle 
$\A_\ep(\ell,m)$, $\ep =1,2$, defined in Section~3, is 
the tangle in $B^3$ obtained by 
filling the ``$D$''--puncture of $\T_\ep(\ell,m)$ with the $\frac12$--tangle. 
That is 
\begin{equation*}
\begin{split}
\A_1 (\ell,m) & = \C\Big(A,B,*,\tfrac12\Big)\\
\A_2 (\ell,m) & = \C\Big(A,*,C,\tfrac12\Big)
\end{split}
\end{equation*}
(with $A,B$ as for $\T_1(\ell,m)$, $\T_2(\ell,m)$). 
See \fullref{Decomp}, which gives a marking to $\A_\ep (\ell,m)$.
\vspace{2pt}

\begin{Defin}
Denote by $\T_\ep (\ell,m)(\frac{p}q)$ the tangle in the 3--ball gotten by 
filling the $D$--puncture of $\T_\ell (\ell,m)$ with the rational tangle 
$\re (p/q)$.
The change in filling $\T_\ep (\ell,m)(\frac12)$ to $\T(\ell,m)(\frac10)$ 
corresponds to a crossing move taking the toroidal tangle $\A_\ep (\ell,m)$ 
to a rational tangle.
We will refer to this as the {\em standard crossing move} on 
$\A_\ep (\ell,m)$.
The arc guiding this crossing change is the {\em standard unknotting arc} 
for $\A_\ep (\ell,m)$.
\end{Defin}
\vspace{2pt}

$\A_\ep(\ell,m)$ $(\ep =1$ or 2) contains an essential Conway sphere $S$, 
which induces a decomposition of the double branched covering 
$\wt{\A}_\ep(\ell,m) = N_1 \cup_T M_2$, where $N_1$ is a Seifert fiber space
over $A^2$ with one exceptional fiber, and $M_2$ is a Seifert fiber space 
over $D^2$ with two exceptional fibers.
See \fullref{Decomp}. 
Thus if $J= J_\ep (\ell,m)$ for some $\ep, \ell,m$, then $J$ is a knot in 
$S^1\times D^2$, with meridian $\mu$, say, and 
$J(\gamma)= \wt{\A}_\ep (\ell,m) = N_1 \cup_T M_2$ contains an essential 
separating torus $T$, for some $\gamma$ such that 
$\Delta (\gamma,\mu) =2$.
\vspace{2pt}

\vspace{2pt}

\begin{thm}\label{thmJlm}
$J_\ep (\ell,m)$ is a hyperbolic knot in the solid torus.
\end{thm}
\vspace{2pt}

We give the proof of \fullref{thmJlm} at the end of this section. 
\vspace{2pt}

The following theorem says the $J_\ep (\ell,m)$ are exactly the 
hyperbolic knots in solid tori which admit non-integral toroidal surgeries.
\vspace{2pt}

\begin{thm}\label{thm2.1}		
Let $J$ be a knot in a solid torus whose exterior is irreducible and atoroidal.
Let $\mu$ be the meridian of $J$ and suppose that $J(\gamma)$ contains an 
essential  torus for some $\gamma$ with $\Delta (\gamma,\mu)\ge2$. 
Then $\Delta (\gamma,\mu) =2$ and $J= J_\ep (\ell,m)$ for some $\ep,\ell,m$.
\end{thm}
\vspace{2pt}

\begin{remark}
$J = J_\ep (\ell,m)$ means there is a homeomorphism of the solid torus 
(possibly orientation-reversing) taking $J$  to $J_\ep (\ell,m)$.
\end{remark}
\vspace{2pt}

The following is the corresponding statement about tangles.
\fullref{thm:doublebranched}(2) says that $\A_\ep (\ell,m)$ is 
determined by its double branched cover. 
At this point, it is not known if the same is true for $\T_\ep (\ell,m)$. 
This complicates the statement of \fullref{thm2.2}.	
\vspace{2pt}

\begin{thm}\label{thm2.2} 		
Let $\T(*,*)$ be a tangle in $S^2\times I$ which is irreducible and atoroidal 
as a $\zed_2$--orbifold.
If $\T(*,\alpha)$ is rational, and $\T(*,\beta)$ is orbifold-toroidal, where 
$\Delta (\alpha,\beta) \ge2$, then the double branched cover 
of $\T(*,*)$ is homeomorphic to the double branched cover of 
$\T_\ep(\ell,m)$ for some $\ep,\ell,m$.
Under this homeomorphism, 
the slopes $\alpha,\beta$ on $\T$ correspond to slopes 
$1/0$, $1/2$ respectively on $D$ of $\C$. 

Furthermore, there are tangle homeomorphisms 
$h_1\co \T(*,\beta)\to \A_\ep(\ell,m) =$\break 
$\T_\ep (\ell,m)(\frac12)$ and 
$h_2\co \T(*,\alpha)\to \T_\ep(\ell,m)(\frac10)$  such that 
$(h_2|\partial) (h_1|\partial)^{-1}$ is the identity (where 
$\partial\T(*,\beta) = \partial\T(*,\alpha)\subset\partial \T(*,*)$, 
and $\T_\ep (\ell,m)(\frac10)$, $\T_\ep (\ell,m)(\frac12)$ are 
marked from \fullref{ScriptT}).
\end{thm}

\begin{Addendum}
$\T(*,\alpha)$ is rational and thus determines a slope, $\frac{p_1}{q_1}$, 
on its boundary.
$\T(*,\beta)$ is the tangle $\T_\ep (\ell,m)$ whose double 
branched cover contains a unique essential annulus. 
The boundary of this annulus determines a slope on the boundary of the cover,
which determines a tangle slope, $\frac{p_2}{q_2}$, on the boundary of 
$\T(*,\beta)$. 
Then $\Delta (p_1/q_1,p_2/q_2) = |p_1q_2 - p_2q_1| >1$.
\end{Addendum}

\begin{proof}[Proof of \fullref{thm2.2}] 
Let $\T(*,*)$ be as in the theorem. 
Let $X$ be its double branched cover. 
Let $X(\alpha),X(\beta)$ be the Dehn fillings of $X$ corresponding to the 
double branched covers of $\T(*,\alpha)$, $\T(*,\beta)$ (resp.). 
Then by assumption, $X$ is irreducible and atoroidal, $X(\alpha)$ is a solid 
torus and $X(\beta)$ is toroidal. 
Since $\Delta (\alpha,\beta) >1$, Corollary~A.2 of \cite{GL3} proves 
that $X$ is the double branched cover of $\C(A,B,C,*)$ where one of $A,B,C$ 
is the empty tangle and the others are rational (using the fact proven 
there that $\frac12\in \{\alpha',\beta',\gamma'\}$).
Furthermore, under this identification, $\alpha =\frac10$, $\beta=\frac12$. 
That is, $\C(A,B,C,\frac10)$ is a rational tangle, and $\C(A,B,C,\frac12)$ 
is orbifold-toroidal. 
By symmetry we may assume that either $B$ or $C$ is the empty tangle above. 
\fullref{lemmaB} and \fullref{lemmaC} below 
show that $A,B,C$ are as in the definition of 
$\T_\ep (\ell,m)$. 
Now the double branched cover of $\T(*,\beta)$ is the same as the 
double branched cover of $\C(A,B,C,\frac12)= \A_\ep (\ell,m)$. 
By \fullref{thm:doublebranched}, $\T(*,\beta)$ and $\A_\ep (\ell,m)$ are 
homeomorphic tangles.
Such a homeomorphism $h_1$ determines a framing on $\partial X(\beta)$ 
($\A_\ep (\ell,m)$ is a marked tangle), hence on $\partial X(\alpha)$. 
The meridian disk of $X(\alpha)$ determines  a rational number in this
framing which corresponds to the rational tangle $\T_\ep (\ell,m)(\frac10)$. 
This implies that $\T(*,\alpha)$ is the same as the rational tangle 
$\T_\ep (\ell,m)(\frac10)$ under the marking determined by $h_1$. 
This is the second paragraph of \fullref{thm2.2}. 
\end{proof}

\begin{proof}[Proof of Addendum to \fullref{thm2.2}]
In the context of the proof of \fullref{thm2.2}, $X(\alpha)$, $X(\beta)$ 
are the double branched covers of $\C(A,B,C,\frac10)$, 
$\C(A,B$, $C,\frac12)$ with $A,B,C$ as in the definition of $\T_\ep(\ell,m)$. 
The slope of the essential annulus of $X(\beta)$ corresponds to the tangle 
slope $\frac10$ in \fullref{ScriptT} (see \fullref{Decomp}). 
The slope corresponding to the rational tangle $\T(*,\alpha)$ corresponds 
to the slope of the meridian disk in $X(\alpha)$. 
This in turn corresponds in \fullref{ScriptT} 
to the $\frac{p}q$ of the 
rational tangle $\C(A,B,C,\frac10)$. 
We need to show then that $|q|>1$. 
For $\T_1(\ell,m)$, we see this by simply noting that $\C(A,B,C,\frac10)$ 
capped off by strands of slope $\frac10$ is the unlink only when $m=0$. 
For $\T_2(\ell,m)$, $\C(A,*,C,\frac10)$ corresponds to the rational 
number $\frac{2m-1}{\ell(2m-1)-2}$. 
That is, $q=\ell (2m-1)-2$. 
The conditions that $|\ell|>1$, $m\ne 0,1$ imply that $|q|>1$.
\end{proof}

\begin{proof}[Proof of \fullref{thm2.1}]
This is the proof of \fullref{thm2.2} using 
the definition of $J_\ep (\ell,m)$, and without 
\fullref{thm:doublebranched}.
\end{proof}

\begin{lem}\label{lemmaB}
If $\C(A,*,C,\frac10)$ is a rational tangle and $\C(A,*,C,\frac12)$ is 
orbifold-toroidal, then $A= \R (s)$, $C= \R (t,2,0)$ for $s,t\in\zed$ with 
$|s|,|2t+1|>1$.
\end{lem}

\begin{proof}
We follow the argument of Lemma~5.1 of \cite{E2}. 
Rewrite $\C(A,*,C,\frac10)$ as in \fullref{FigB1}.		
{From} this we deduce that $C' = \R (\frac1{t'})$, $t'\in\zed$. 
Thus $\C(A,*,C,\frac10)$ is as in \fullref{FigB2}. 
(Note that  our convention in this paper is that twist boxes represent 
{\em vertical} twists).

\begin{figure}[ht!]
\centerline{\includegraphics[height=0.65truein]{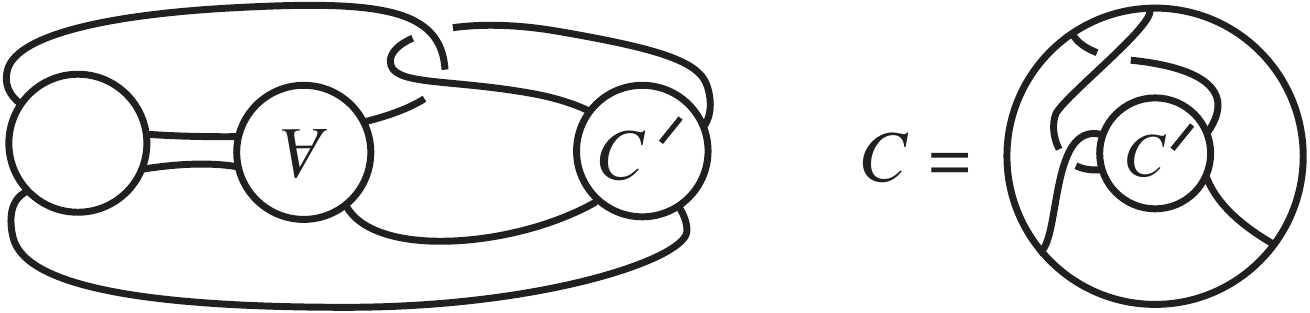}}
\caption{}		
\label{FigB1}
\end{figure}

\begin{figure}[ht!]
\centerline{\includegraphics[height=0.8truein]{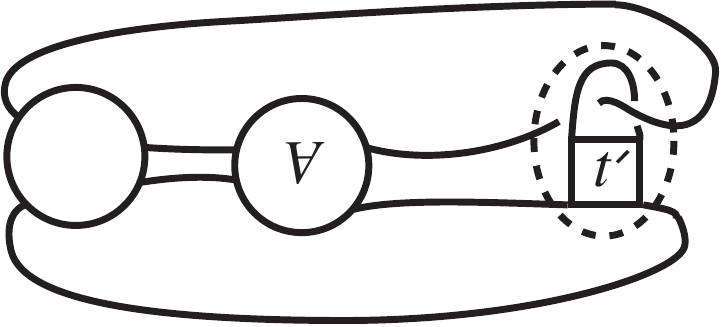}}
\caption{}		
\label{FigB2}
\end{figure}

The tangle encapsulated in the Conway sphere in \fullref{FigB2} 
is $\R (2,t',0)= \R (\frac2{2t'+1})$. 
Thus either $t'=0,-1$ or $A= \R (s)$ for some $s\in\zed$. 
In either case, $C= \R (t',1,1,0) = \R (t,2,0)$ where $t= -(t'+1)$. 

\fullref{FigB3}	
shows $\C(A,*,C,\frac12)$, 
where $C = \R (t,2,0)$ corresponds to the rational number $\frac{t}{2t+1}$. 
Thus $\C(A,*,C,\frac12)$ orbifold-toroidal implies that $\Delta (\frac10,
\frac{t}{2t+1})>1$. 
Thus $|2t+1|>1$. 
Consequently, $t' = -(t+1) \ne -1,0$ and $A= \R (s)$. 
Again, that $\C(A,*,C,\frac12)$ is orbifold-toroidal guarantees that $|s|>1$.
\end{proof} 

\begin{figure}[ht!]		
\centerline{\includegraphics[height=1.2truein]{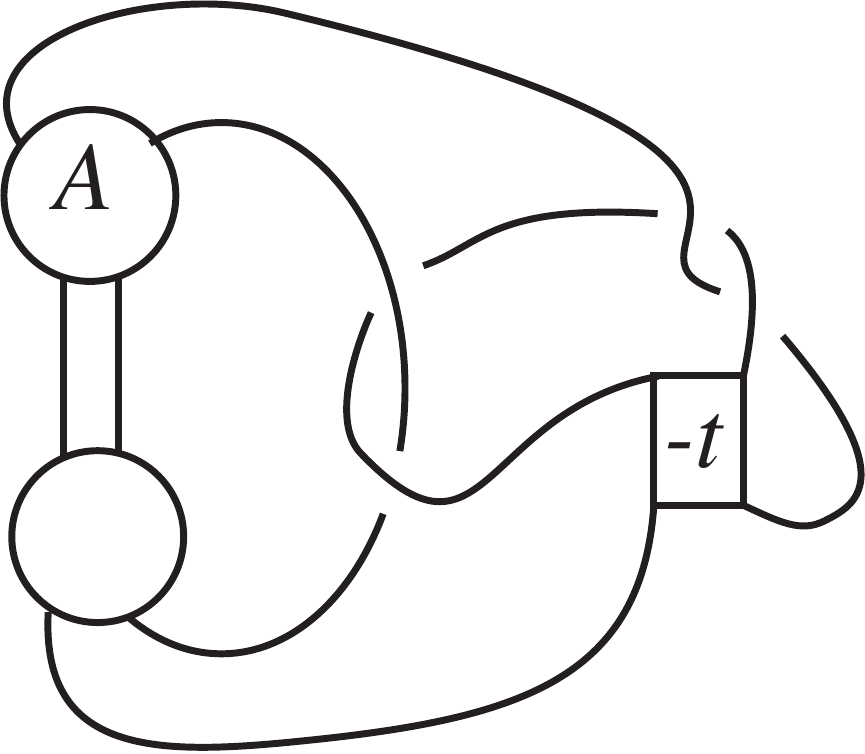}}
\caption{}		
\label{FigB3}
\end{figure}

\begin{lem}\label{lemmaC}
If $\C(A,B,*,\frac10)$ is a rational tangle then, up to symmetry exchanging
$A$ and $B$, $A= \R (s)$ and $B=\R (t,-s)$ for $s,t\in\zed$. 
If $\C(A,B,*,\frac12)$ is orbifold-toroidal then $|s|>1$, $t\ne 0$ and 
$(s,t) \ne (2,1),(-2,-1)$.
\end{lem}

\begin{figure}[ht!]
\centerline{\includegraphics[height=1truein]{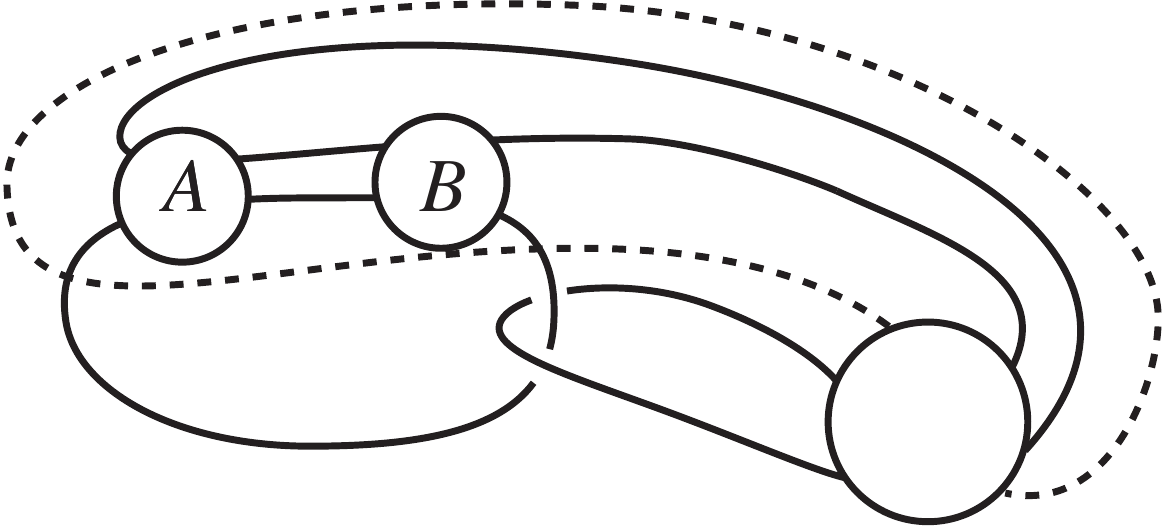}}
\caption{}		
\label{FigC1}
\end{figure}

\begin{figure}[ht!]
\centerline{\includegraphics[height=1.3truein]{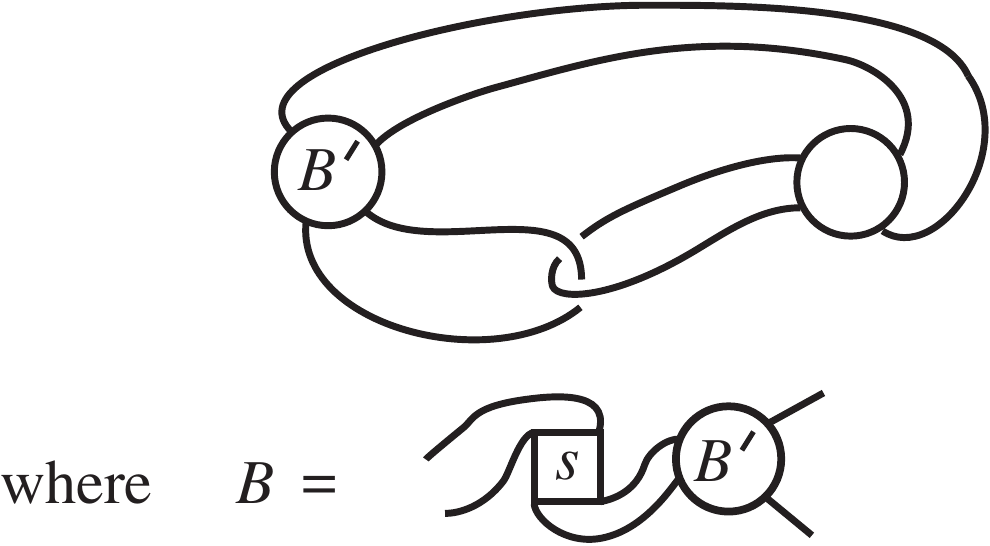}}
\caption{}		
\label{FigC2}
\end{figure}

\begin{figure}[ht!]
\centerline{\includegraphics[height=1.1truein]{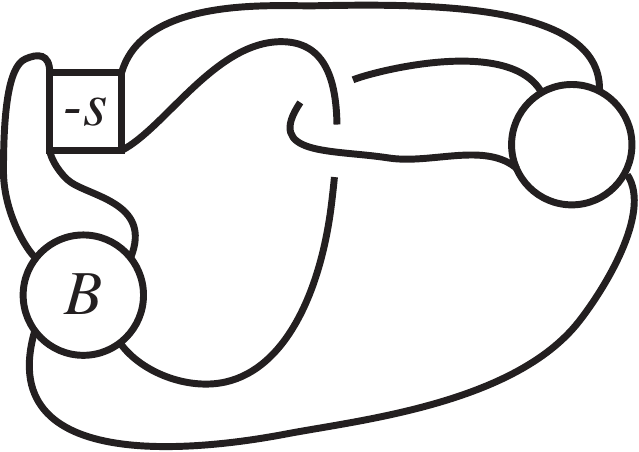}}
\caption{}		
\label{C3}
\end{figure}

\begin{proof}
Isotoping $\C(A,B,*,\frac10)$ to \fullref{FigC1} 
we see that one of $A$ or $B$ must be an integral tangle. 
By symmetry we assume it is $A$, $A= \R (s)$.
Thus we are as in \fullref{FigC2}, 
from which we see that $B'= \R (1/t)$. 
Thus $B = \R (t,-s)$. 
Now assume $\C(A,B,*,\frac12)$ is orbifold-toroidal. 
See \fullref{C3}. 
Then $|s| >1$, and $B= \R (t,-s)$ corresponds to tangle slope $\frac{1-st}t$.
Thus $\Delta (\frac01,\frac{1-st}t) >1$. 
That is, $|1-st| >1$. 
Thus $t\ne 0$ and $(s,t)\ne (2,1),(-2,-1)$.
\end{proof}

\begin{figure}[ht!]
\centerline{\includegraphics[height=1.7truein]{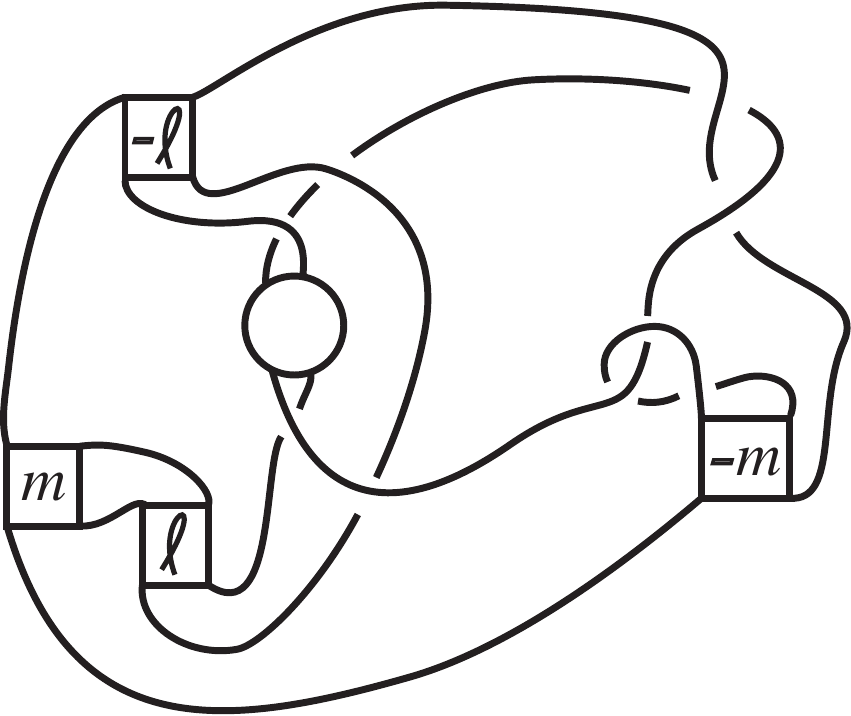}}
\caption{}		
\label{J1gamma}
\end{figure}

\begin{proof}[Proof of \fullref{thmJlm}]
We prove this for $J_1(\ell,m)$. 
The proof for $J_2(\ell,m)$ is similar. 
Recall that by attaching the appropriate ``$C$''--tangle to $\T_1(\ell,m)$
we get the tangle $B(\ell,m,n,0)$ of \cite{E1}. 
There are infinitely many such fillings corresponding to different 
values of $n$. 
Looking at double branched covers, this says that the corresponding Dehn 
fillings of $X$, the exterior of $J_1(\ell,m)$, are the exteriors of the 
hyperbolic knots $k(\ell,m,n,0)$. 
Denote the two components of $\partial X$ as $\partial_1X$, $\partial_2X$, 
where $\partial_1X$ is the component along which these fillings are made 
(corresponding to the boundary of the ambient solid torus of $J_1(\ell,m)$).
Because infinitely many fillings of $X$ are hyperbolic, either $X$ is 
hyperbolic or there is a cable space along $\partial_1X$ 
(Theorem~2.4.4 of \cite{CGLS}). 
We assume the latter for contradiction. 
Then the slope of each of these Dehn fillings is distance~1 from a 
unique slope $\gamma$, and furthermore, Dehn filling $X$ along $\gamma$ 
has a lens space summand. 
The $\zed_2$--orbifold quotient of the $\gamma$--Dehn filling of $X$ is 
pictured in \fullref{J1gamma} 
(the $C$--tangle for this picture corresponds 
to replacing the $n$ twist box of $B(\ell,m,n,0)$ with $\R(1/0)$). 
But inserting $\R (1/0)$ into the ``$D$''--tangle of 
\fullref{J1gamma} 
gives the unlink of two components. 
That is, there is a filling of $\partial_2X$ such that, along with 
the filling of $\partial_1X$ along $\gamma$, gives $S^2\times S^1$.
But this contradicts the fact that first filling $X$ along $\gamma$ yields 
a lens space summand. 
\end{proof}

\section{Non-integral surgery and the JSJ--decomposition} 

In this section we consider a non-integral Dehn filling $X(\gamma)$ on the 
exterior $X$ of a knot in $S^3$, and analyze the relation between the 
JSJ--decompositions of $X$ and $X(\gamma)$. 

If $M$ is an irreducible 3--manifold we shall denote by $\seif (M)$ the 
disjoint union of the Seifert fibered pieces of the JSJ--decomposition of $M$. 
In the case  of a knot exterior $X$, 
the possible components of $\seif (X)$ have been described by Jaco 
and Shalen.

\begin{lem} \label{lem:JS}
{\rm\cite[Lemma VI.3.4]{JS}}\qua
Let $X$ be the exterior of a knot in $S^3$, and let $W$ be a component 
of $\seif (X)$.
Then $W$ is either a torus knot space, a cable space, or a composing space.
\end{lem}

The relation between $\bT(X)$ and $\bT(X(\gamma))$ is 
described in the following theorem and its addendum.

\begin{thm}\label{thm:JSJ}
Let $X$ be the exterior of a knot $k$ in $S^3$, and suppose 
$\Delta (\gamma,\mu)\ge 2$ where $\mu$ is the meridian of $k$. 
Let $W$ be the component of $X$ cut along $\bT(X)$ that contains $\partial X$.
Then exactly one of the following four possibilities holds.
\begin{itemize}
\item[(1)] $\bT(X(\gamma)) = \bT(X)$;
\item[(2)] 
$\bT(X) = \emptyset$, $X$ is hyperbolic, $k$ is an Eudave-Mu\~noz 
knot $k(\ell,m,n,p)$,\break $\Delta (\gamma,\mu) =2$, 
$X(\gamma) = M_1 \cup_T M_2$, 
where $M_i$ is a Seifert fiber space over $D^2$ with two exceptional fibers, 
$i=1,2$, and $\bT(X(\gamma))=T$;
\item[(3)] 
$k$ is contained in a tubular neighborhood $N(k_0)$ of a non-trivial
knot $k_0$ as a $J_\ep (\ell,m)$--satellite of $k_0$, $\partial W = T_0 \cup
\partial X$, where $T_0 = \partial N(k_0)$, $\Delta (\gamma,\mu)=2$, 
$W(\gamma) = N_1 \cup_T M_2$, where $N_1$ is a Seifert fiber space over $A^2$ 
with one exceptional fiber and $M_2$ is a Seifert fiber space over $D^2$ with 
two exceptional fibers, and $\bT(X(\gamma)) = \bT(X) \cup T$; 
\item[(4)] 
$k$ is a $(p,q)$--cable of a non-trivial knot $k_0$ with exterior $X_0$, 
$q\ge 2$, $\partial W = T_0 \cup \partial X$ as in (3), 
$\gamma = \frac{npq\pm 1}{n}$, $n\ge2$, and $\bT(X(\gamma)) = \bT(X) - T_0$.
\end{itemize}
\end{thm}

We spell out more details about cases (1) and (4) of \fullref{thm:JSJ}
in the following addendum, where $k_\gamma$ denotes the core of the 
Dehn filling solid torus in $X(\gamma)$.

\begin{adden}\label{add:JSJ}
In case {\rm(1)} of \fullref{thm:JSJ} we have 
\begin{itemize}
\item[(a)] if $W= X$ then $X(\gamma)$ is atoroidal;
\item[(b)] if $W\ne X$ then $W(\gamma)$ is hyperbolic if and only if $W$ 
is hyperbolic;
\item[(c)] if $W\ne X$ then $W(\gamma)$ is Seifert fibered if and only if 
$W$ is Seifert fibered. 
If $W$ is a cable space then $W(\gamma)$ is a Seifert fiber space over $D^2$ 
with two exceptional fibers, and if $W$ is a composing space with $(n+1)$
boundary components then $W(\gamma)$ is a Seifert fiber space over an 
$n$--punctured sphere with a single exceptional fiber of multiplicity 
$\Delta (\gamma,\mu)$. 
In both cases $k_\gamma$ is isotopic to an exceptional fiber of $W(\gamma)$.
\end{itemize}

In case {\rm(4)} of \fullref{thm:JSJ} we have 
\begin{itemize}
\item[(a)] $W$ is a cable space and $W(\gamma)$ is a solid torus, with 
meridian $\gamma_0$, say, on $T_0$;
\item[(b)] conclusion {\rm(1)} holds for $X_0 (\gamma_0)$;
\item[(c)] $k_\gamma$ is isotopic to an $(r,s)$--cable of the core of 
$W(\gamma)$, for some $s>1$.
\end{itemize}
\end{adden}

\begin{proof}[Proof of \fullref{thm:JSJ}] 

{\bf Case I\qua $|\partial W| =1$}

Here $\bT(X) = \emptyset$ and $W=X$, so $X$ is either hyperbolic or Seifert 
fibered. 
In the first case, either (1) or (2) holds by \cite{GL3}. 
In the second case, $k$ is a torus knot, $X(\gamma)$ is a Seifert fiber space
over $S^2$ with at most three exceptional fibers, and $\bT(X(\gamma)) = 
\bT(X) =\emptyset$.
\medskip

{\bf Case II\qua $|\partial W| =2$}

There are two subcases.

{\bf (A)\qua $W$ hyperbolic}

$W(\mu)$ is a solid torus; therefore $W(\gamma)$ is irreducible and 
$\partial$--irreducible, by \cite{S2} and \cite{CGLS} respectively.
If $W(\gamma)$ is hyperbolic then $\bT(X(\gamma)) =\bT(X)$, and (1) holds. 
By \cite[Proposition 9]{MZ}, $W(\gamma)$ is not Seifert fibered. 
Hence we may assume that $W(\gamma)$ is toroidal. 
Then by \fullref{thm2.1} 
$k$ is a $J_\ep (\ell,m)$--satellite in $W(\mu)\cong S^1\times D^2$. 
Then $W(\gamma) = N_1\cup_T M_2$  as in conclusion~(3). 
Let $\varphi$ be the fiber of $N_1$ on $T_0 = \partial W(\gamma)$.

Consider the component $Z\ne W$  of the JSJ--decomposition of $X$ with 
$T_0 \subset \partial Z$. 
Assume that $Z$ is a Seifert fiber space, with fiber $\psi$ on $T_0$. 
We will show that $\psi \ne \varphi$; hence $\bT(X(\gamma)) = \bT(X)\cup T$ 
and (3) holds.

Let $\mu_0$ be the meridian of $W(\mu) \cong S^1\times D^2$. 
Since $\Delta (\mu_0,\varphi) \geq 2$ by the Addendum to \fullref{thm2.2}, 
it suffices to prove that $\Delta(\mu_0,\psi) \le1$.

By \fullref{lem:JS}, $Z$ is either 
(i)~a torus knot exterior, or
(ii)~a cable space, or 
(iii)~a composing space. 
Let $X_0 = \overline{X-W}$; so $\partial X_0 = T_0$. 
Since $X_0 (\mu_0) \cong S^3$, we must have 
$\Delta (\mu_0,\psi)=1$ in case~(i), 
$\Delta (\mu_0,\psi)=1$ in case~(ii) (since $Z\cup W(\mu)$ is a solid torus),
and, $\mu_0=\psi$ in case~(iii) (since $\partial (Z\cup W(\mu))$ must be 
compressible).

{\bf (B)\qua $W$ Seifert fibered.}

Then $W$ is a cable space. 
Thus $k$ is a $(p,q)$--cable, $q\ge2$, of a non-trivial knot $k_0$. 
Let $\varphi$ be the slope on $\partial X$ of the Seifert fiber of $W$; 
thus $\varphi = pq/1$ with respect to the usual meridian-longitude basis.

If $\Delta (\gamma,\varphi)\ge2$ then $W(\gamma)$ is a Seifert fiber space 
over $D^2$ with two exceptional fibers, which is atoroidal. 
Therefore $\bT(X(\gamma)) = \bT(X)$.

Since $\Delta (\gamma,\mu)\ge2$ ($\mu=1/0=$ meridian of $k$), 
$\gamma\ne \varphi$. 
So assume $\Delta (\gamma,\varphi) =1$; ie, 
$\gamma =\frac{npq\pm1}{n}$, $n\ge2$. 
Then $W(\gamma)$ is a solid torus, with  meridian $\gamma_0$, say, on $T_0$. 
Let $\mu_0$ be the slope on $T_0$ of the meridian of the solid torus $W(\mu)$.
Then $\Delta (\gamma_0,\mu_0) = q^2 \Delta (\gamma,\mu) \ge8$. 
Let $X_0 = \overline{X-W}$, the exterior of $k_0$. 
By induction on the number of components of $\bT(X)$, we may assume that the 
theorem holds  for $X_0$. 
(The start of the induction is Case~(1) above.) 
Since $\Delta (\gamma_0,\mu_0)>2$, conclusions~(2) and (3) of the theorem 
do not hold for $X_0$. 
If (1) holds for $X_0$, then we get 
$\bT(X(\gamma)) = \bT(X_0(\gamma_0)) = \bT(X_0) = \bT(X) - T_0$, 
which is conclusion~(4) for $X$. 
Finally, assume that (4) holds for $X_0$. 
So $k_0$ is a $(p_1,q_1)$--cable, $q_1 \ge2$, of a non-trivial knot $k_1$, 
and $\gamma_0 = \frac{n_1p_1q_1\pm 1}{n_1}$, $n_1\ge2$, with respect to the 
usual basis for $k_0$. 
But we also have $\gamma_0 = \frac{npq\pm 1}{nq^2}$, $n\ge2$, and this is 
a contradiction (see \cite[page 704]{G}).
\medskip

{\bf Case III\qua $|\partial W|\ge 3$}

Let the components of $\partial W-\partial X$ be $T_1,\ldots,T_n$, $n\ge2$. 
Then the components of $\overline{X-W}$ are $Y_1,\ldots,Y_n$, say, 
where $\partial Y_i = T_i$ and $T_i$ is incompressible in $Y_i$, $1\le i\le n$.
Since $T_i$ compresses in $S^3= X(\mu)$, it compresses in $W(\mu)$. 
Hence $W(\mu)$ is reducible. 
(In fact it is easy to show that $W(\mu)$ is a connected sum of $n$ 
solid tori.)

Again we distinguish two subcases.

{\bf (A)\qua $W$ hyperbolic}

Since $\Delta (\gamma,\mu) \ge2$, $W(\gamma)$ is irreducible by \cite{GL2}. 
Also, since $|\partial W| \ge3$, $W(\gamma)$ is atoroidal and anannular 
by \cite{W}. 
Therefore $W(\gamma)$ is hyperbolic, and $\bT(X(\gamma)) = \bT(X)$.  
This is conclusion~(1).

{\bf (B)\qua $W$ Seifert fibered}

Since $|\partial W| = n+1$, $n\ge2$, $W$ is a composing space by 
\fullref{lem:JS}. 
Also, the meridian $\mu$ is the Seifert fiber of $W$, since $W(\mu)$ 
is reducible. 
Hence $W(\gamma)$ is a Seifert fiber space over an $n$--punctured sphere 
with one exceptional fiber, of multiplicity $\Delta (\gamma,\mu)$. 
Since the Seifert fibers of $W$ and $W(\gamma)$ are the same on each $T_i$, 
we have $\bT(X(\gamma)) = \bT(X)$, and (1)~holds. 
\end{proof}

\begin{proof}[Proof$\!$ of$\!\!$ Addendum$\!$ \ref{add:JSJ}] 
$\!\!\!$This follows by examining the proof of \fullref{thm:JSJ}.
\end{proof}

\section{Main theorem}\label{sec:mainthm}	

Recall from Section 2 the definition of the characteristic 2--sided toric 
2--suborbi\-fold $\bT (K)$ of a prime knot $K$ in $S^3$. 
Let $\seif (K) = \seif(O(K))$ be the disjoint union of the $S^1$--fibered 
components of $O(K)$ cut along $\bT(K)$. 

Let $(a,\partial a)$ be an unknotting arc for $K$. 
As described in Section~2, a relative regular neighborhood of $(a,\partial a)$ 
in $(S^3,K)$ determines a marked tangle $\T_0$ which is replaced with a 
tangle $\T'_0$ under the crossing move. 
Let $M$ be the double branched cover of $S^3$ along $K$, $V_0$ be the 
solid torus preimage of $\T_0$ under the branched covering, and $X=M-V_0$.
Then $M= X(\gamma)$ and $S^3 = X(\mu)$ where $\Delta (\mu,\gamma)=2$.
Let $k= k_\mu$ be the knot in $S^3$ of which $X$ is the exterior, and let 
$k_\gamma$ be the core of $V_0$ in $M$.

\begin{Defin}
The unknotting arc $(a,\partial a)$ is said to be an
{\em $(r,s)$--cable of an exceptional fiber of $\seif (K)$} 
iff $k_\gamma$ is an $(r,s)$--cable of an exceptional fiber in 
$\seif (M)$.
\end{Defin}

\begin{remark}
If $(a,\partial a)$ is an $(r,s)$--cable of an exceptional fiber of 
$\seif (K)$, then the corresponding $\T_0$ lies in a rational tangle 
$\R(p/q)$ in $\seif (K)$ which is the quotient of a neighborhood of this
exceptional fiber.
The tangle $\R(p/q) - \T_0$ in $S^2\times I$, has a double branched cover 
which is a cable space. 
By \fullref{lem:Fix-tau2}, 
$\R(p/q)-\T_0$ is homeomorphic to $\M(\frac{v}w,*)$ for some $v,w\in\zed$.
The results of applying the crossing move associated to such an 
$(a,\partial a)$ 
are further discussed in Lemmas~\ref{lemK}, 
\ref{lem:crossingchange}, \ref{lemH}, and \fullref{thm8.3}.
\end{remark}

\begin{lem}\label{lem6.1}
Let $(a,\partial a)$ be an unknotting arc for $K$. 
One of the following holds:
\begin{itemize}
\item[{(1)}] 
$(a,\partial a)$ can be isotoped in $(S^3,K)$ to be disjoint from $\bT(K)$. 
Furthermore, if $\bT(K)\ne\emptyset$ and $(a,\partial a)$ can be isotoped 
into $\seif (K)$, then $a$ is isotopic to an $(r,s)$--cable of an 
exceptional fiber in $\seif (K)$, for some $s\ge1$.
\item[{(2)}] 
{\rm(a)}\qua $K$ is an EM--knot $K(\ell,m,n,p)$.
\item[{}] 
{\rm(b)}\qua $O(K)$ has a unique connected, incompressible, 2--sided, toric 
2--suborbifold $S$, a Conway sphere, $K$ has an unknotting arc $(b,\partial b)$ 
with ${|b\cap S| =1}$ (the standard unknotting arc for $K(\ell,m,n,p)$),  
and no unknotting arc is disjoint from $S$.
\item[{(3)}] 
$K$ is the union of essential tangles $\P \cup \P_0$, where $\P_0$ is the 
EM--tangle $\A_\ep(\ell,m)$, $\partial\P_0\subset \bT(K)$, and $(a,\partial a)$
can be isotoped into $\P_0$. 
If $\T(*,*)$ is the exterior in $\P_0$ of the crossing ball corresponding 
to $(a,\partial a)$, then $\T(*,*)$ is as described in 
\fullref{thm2.2}.
\end{itemize}
\end{lem}

\begin{proof}[Proof of \fullref{lem6.1}] 
Let $M$, $k=k_\mu,k_\gamma$ be as described above.
We are now in the context of \fullref{thm:JSJ}. 
Possibilities (1), (2), and (3) in the conclusion of \fullref{thm:JSJ} 
will lead to conclusions (1), (2), and (3), respectively of 
\fullref{lem6.1}, and possibility (4) will lead to conclusion (1).
\vspace{2pt}

Let $h\co M\to M$ be the covering involution, with quotient orbifold 
$\widehat M = O(K)$. 
Write $V_0 = V_\gamma$, with quotient $\widehat V_\gamma$ the 3--suborbifold 
$\T_0$ of $O(K)$. 
\vspace{2pt}

(1)\qua  Here $\bT(M) = \bT(X)$. 
The covering involution $h\co M\to M$ restricts to $h\co X\to X$, and we can 
isotop $\bT(X)$ in $X$ to be $h$--invariant \cite{MS}. 
Let $T$ be a component of $\bT(X)$. 
Then $T$ separates $X$ into two components, one of which contains 
$\partial X$. 
It follows that if $h(T) =T$ then $h$ preserves the sides of $T$, and 
hence $\bT^+ (M) = \bT^+ (X) = \bT(X) = \bT(M)$.
\vspace{2pt}

By (1)(a) of Addendum~\ref{add:JSJ}, if $\bT(X) = \emptyset$ then 
$M= X(\gamma)$ is atoroidal, and hence $\bT(K)=\emptyset$.
Thus conclusion~(1) holds trivially.
\vspace{2pt}

If $\bT(X) \ne\emptyset$, then $X= X_0 \cup W$, say, $X_0\ne \emptyset$. 
Hence $M= X_0 \cup W(\gamma)$, and, taking quotients, $O(K)= \hat X_0 \cup 
\widehat W(\gamma)$, where $\widehat V_\gamma = \T_0 \subset \widehat W
(\gamma)$. 
By (1)(b) of Addendum~\ref{add:JSJ}, if $W$ is hyperbolic then  $W(\gamma)$ 
is hyperbolic, and so $\widehat W(\gamma)$ is atoroidal and is a 
component of $O(K)$ cut along $\bT(K)$. 
Hence $\bT(K)$ can be orbifold-isotoped off $\widehat W(\gamma)$, in 
particular, off $\T_0$. 
If $W$ is Seifert fibered, then by (1)(c) of Addendum~\ref{add:JSJ}
the Seifert fibering of $W$ extends to a Seifert fibering of $W(\gamma)$. 
Thus $\widehat W(\gamma)$ is $S^1$--fibered, and $\widehat V_\gamma = \T_0$ is 
a neighborhood of an exceptional (orbifold) fiber.
Also, since $\bT(M) = \bT(X)$, $W(\gamma)$ is a component of $\seif (M)$, and 
hence $\widehat W(\gamma)$ is a component of $\seif (K)$.
\vspace{2pt}

(2)\qua Here $k^*$ is an Eudave-Mu\~noz knot $k(\ell,m,n,p)$. 
By \fullref{thm:doublebranched}(1), $K= K(\ell,m,n,p)$, $M= M_1\cup_T M_2$, 
$\bT(M) =T$, and we may isotop $T$ so that $h(M_i) = M_i$, $i=1,2$.
Hence $\bT(K) = S= \widehat T$.
The facts that $K$ has an unknotting arc $b$ with $|b\cap S|=1$, and 
no unknotting arc disjoint from $S$, are proved in \cite{E1}.
\vspace{2pt}

(3)\qua Here $k^*$ is a $J_\ep (\ell,m)$--satellite of $k_0$. 
Thus $X = X_0\cup_{T_0} W$, where $X_0$ is the exterior of $k_0$, 
$T_0 = \partial X_0$, and $W$ is the exterior of $J_\ep (\ell,m)$ in 
$S^1\times D^2$. 
Also, $M= X_0 \cup_{T_0} W(\gamma)$, $W(\gamma)\cong N_1\cup_T M_2$ as in 
\fullref{thm:JSJ}, and $\bT(M) = \bT(X) \cup T$. 
Since $W$ is hyperbolic, $T_0 \subset \bT(X)$.
\vspace{2pt}

Now $h\co M\to M$ leaves $X$ invariant. 
Hence we can isotop $\bT(X)$ in $X$ to be $h$--invariant \cite{MS}. 
In particular, we must clearly have $h(T_0)=T_0$. 
Hence $h$ leaves $W(\gamma)$ invariant. 
$\Fix (h)$ cannot be disjoint from $T_0$ or completely lie in $T_0$, 
otherwise $h$ would give rise to an involution on $S^3$ whose fixed set 
was $k_0$ or a satellite of $k_0$---contradicting the $\zed_2$--Smith 
Conjecture \cite{Wa}. 
In particular the quotient of $W$ under $h$, the exterior of the crossing 
ball corresponding to $(a,\partial a)$, is a tangle  $\T(*,*)$ in 
$S^2\times I$ satisfying the hypotheses of \fullref{thm2.2}; 
hence $\T(*,*)$ is as described there. 
This implies $\P_0 \cong\widehat W(\gamma) = \T(*,\beta) \cong \A_\ep (\ell,m)$ 
for some $\ep,\ell,m$. 
Since $\bT(M) = \bT(X) \cup T$, $N_1$ and $M_2$ of $W(\gamma)$ are 
components of $\seif (M)$. 
Furthermore, $\bT^+ (M) = \bT(M)$ since $h$ preserves the sides of $T$ 
and since $\bT^+(X) = \bT(X)$ by the argument for (1). 
Thus $\partial \P_0 \in \bT(K)$ and $\widehat N_1,\widehat M_2$ are 
components of $\seif (K)$.
\vspace{2pt}

(4)\qua Here $W$ is a cable space and $W(\gamma)$ is a solid torus.
Let $X_0 = \overline{X-W}$, and $\partial X_0 = T_0$. 
Then $M= X_0 \cup_{T_0} W(\gamma)$, and $\bT(M) = \bT(X) -T_0 = \bT(X_0)$.
Since $h$ leaves $X$ invariant, we can isotop $\bT(X)$ in $X$ to be 
$h$--invariant.
Then $h(T_0) = T_0$, $\bT(X_0)$ is $h$--invariant, and $\bT^+ (M) =\bT^+(X_0) 
= \bT(X_0) = \bT(M)$. 
Therefore $\bT(K)$ is the quotient $\widehat \bT(X_0)$.
\vspace{2pt}

Let $\mu_0$ be the meridian of $K_0$ on $T_0$, and let $\gamma_0$ be the 
meridian of the solid torus $W(\gamma)$. 
Then $\Delta (\gamma_0,\mu_0) >2$, and $M= X(\gamma) = X_0(\gamma_0)$.
Let $W_0$ be the component of $X_0$ cut along $\bT(X_0)$ that contains 
$\partial X_0$.
By case~(4)(b) of Addendum~\ref{add:JSJ}, (1) holds for $X_0(\gamma_0)$; 
thus (1)(a), (1)(b) and (1)(c) of Addendum~\ref{add:JSJ} hold for $X_0,W_0$.
Conclusions 1(a) and 1(b) now follow from the argument in case~(1) above 
applied to $X_0,W_0,\gamma_0$.
\end{proof}
\vspace{2pt}

\begin{thm}\label{thm:Main}		
Let $K$ be a knot with unknotting number 1. 
Then one of the following three possibilities holds.
\begin{itemize}
\item[(1)] {\rm(a)}\qua 
Any unknotting arc $(a,\partial a)$ for $K$ can be isotoped 
in $(S^3,K)$ so that $a\cap \bT(K) = \emptyset$.

\item[{}]  {\rm(b)}\qua  
If $\bT(K)\ne\emptyset$ and $K$ has an unknotting arc $a$ in $\seif (K)$ then 
$a$ is isotopic to an $(r,s)$--cable of an exceptional fiber of $\seif (K)$, 
for some $s\ge1$.

\item[(2)] {\rm(a)}\qua  $K$ is an EM--knot $K(\ell,m,n,p)$.

\item[{}]  {\rm(b)}\qua  
$O(K)$ has a unique connected incompressible 2--sided toric 2--suborbifold $S$, 
a Conway sphere, $K$ has an unknotting arc $a$ with $|a\cap S|=1$ (the 
standard unknotting arc for $K(\ell,m,n,p)$) , and 
$K$ has no unknotting arc disjoint from $S$.

\item[(3)] 
$K$ is the union of essential tangles $\P\cup \P_0$, where 
$\P_0$ is an EM--tangle $\A_\ep (\ell,m)$ and $\partial\P_0$ is in $\bT(K)$.
Any unknotting arc for $K$ can be isotoped into $\P_0$.
The standard unknotting arc for $\A_\ep (\ell,m)$ is an unknotting arc 
for $K$.
\end{itemize}
\end{thm}

\begin{remarks}
(A)\qua If $K= K(\ell,m,n,p)$, then (2) must hold. 
If $K= \P\cup \P_0$ where $\P_0$ is an EM--tangle then (1) or (3) may hold. 
If (1) holds then, still,
by \fullref{lem9.6} (and \fullref{def:EMtangle}), 
any unknotting arc of $K$ can be isotoped 
into $\P_0$ (hence into an exceptional fiber of $\P_0$).

(B)\qua  In conclusion (3), to say that $K$ is 
unknotted by the standard unknotting 
arc for $\A_\ep (\ell,m)$ (as described in Section~4)
we mean that there is a tangle homeomorphism 
from $\P_0$ to $\A_\ep(\ell,m)$ which makes this identification. 
Any two such will differ by an isotopy of the tangle ball fixed on 
the boundary, which will isotop the two unknotting arcs.
Indeed, any homeomorphism of tangles, $h\co \A_\ep (\ell,m) \to \A_\ep (\ell,m)$, 
preserves the markings, hence is isotopic to the identity.
\end{remarks}

\begin{Quest}
Is any unknotting arc for $K(\ell,m,n,p)$ or for $\A_\ep(\ell,m)$ 
isotopic to its standard unknotting arc?
\end{Quest}

One approach to the above question would be to prove an analog of 
\fullref{thm:doublebranched}  for the exteriors of $k(\ell,m,n,p)$,
$J_\ep (\ell,m)$ (resp.). 
That is, show that there is a unique tangle quotient arising from 
involutions on any such knot exterior.

\begin{proof}[Proof of \fullref{thm:Main}]
\fullref{thm:Main}(2) is the same as \fullref{lem6.1}(2). 
So we assume $K$ is not an EM--knot. 
Furthermore we may assume that $K$ can be written as the union of 
essential tangles $\P\cup \P_0$ with $\P_0$ homeomorphic to $\A_\ep(\ell,m)$ 
and $\partial \P_0\subset \bT(K)$ (otherwise (1) holds for $K$ 
by \fullref{lem6.1}).

If $(a,\partial a)$ is an unknotting arc for $K$ which cannot be isotoped 
into $\P_0$, then \fullref{lem6.1}(1) applies. 
But this says the unknot can be written as $\T\cup \A_\ep (\ell,m)$ 
for some tangle $\T$ and some $\ep,\ell,m$. 
This contradicts \fullref{lem9.6} and \fullref{def:EMtangle}. 
Thus any unknotting arc for $K$ can be isotoped into $\P_0$.

If conclusion (1) of \fullref{thm:Main} does not hold then there is 
an unknotting arc $(a,\partial a)$ satisfying \fullref{lem6.1}(3). 
Let $\T(*,*)$ be the exterior of the crossing ball corresponding to 
$(a,\partial a)$ and $h_1 \co \T(*,\beta) \to \T_\ep (\ell,m)(\frac12)$, 
$h_2 \co \T(*,\alpha) \to \T_\ep (\ell,m)(\frac10)$ be the tangle 
homeomorphisms provided by \fullref{thm2.2}. 
The standard unknotting arc for $\A_\ep (\ell,m)$ corresponds to the 
crossing move $\T_\ep (\ell,m)(\frac12)\allowbreak\to \T_\ep (\ell,m)(\frac10)$.
Thus $h_1$ identifies the standard unknotting arc of $\P_0 \cong \A_\ep(\ell,m)$
for which we are looking. 
That is, performing the crossing move on $K$ dictated by the standard 
unknotting arc gives the knot gotten by gluing $\T_\ep (\ell,m)(\frac10)$ 
to $\P$ via $h_1^{-1}|\partial$. 
Since $(h_2|\partial)(h_1|\partial)^{-1}$ is the identity, this is the 
same as gluing $\T_\ep (\ell,m)(\frac10)$ to $\P$ via $h_2^{-1}|\partial$, 
which is the unknot by assumption.\end{proof}

The following is a generalization of the result of \cite{ST1,ST2} 
(see also \cite{Kob}) that 
an unknotting arc for a satellite knot can always be taken to be 
disjoint from the companion 2--torus.

\begin{cor}\label{cor4.2}
Let $K$ be a knot with unknotting number 1, that is neither an EM--knot 
nor a knot containing an EM--tangle with essential boundary. 
Let $F$ be an incompressible 2--sided toric 2--suborbifold of $O(K)$. 
Then any unknotting arc $(a,\partial a)$ for $K$ can be isotoped in 
$(S^3,K)$ so that $a\cap F=\emptyset$.
\end{cor} 

This is an immediate consequence of \fullref{thm:Main} and the following 
lemma. 

\begin{lem}\label{lem4.3} 
Let $F$ be an incompressible 2--sided toric 2--suborbifold of $O(K)$. 
Then $F$ is orbifold isotopic to a vertical suborbifold of $\seif (K)$. 
\end{lem}
\vspace{-2pt}

\begin{proof} 
The fact that $F$ is isotopic into $\seif (K)$ follows from the discussion 
in \cite{BS} beginning at the paragraph immediately preceding Lemma~7 on 
page~456, and ending at the statement ``and therefore that 
$F\cap F'=\emptyset$'' near the bottom of page~457.
\vspace{-2pt}

So assume $F\subset \seif (K)$. 
Any component of $F$ that is boundary parallel in $\seif (K)$ can be 
isotoped to be vertical, so by \cite[Verticalization Theorem 4]{BS} it is 
enough to show that $F$ cannot be isotoped to be horizontal. 
Let $p\co M\to O(K)$ be the double branched covering projection. 
If $F$ were a horizontal 2--suborbifold of $O(K)$, then $p^{-1}(F)$ would 
be a horizontal surface in $M$. 
But $H_2 (M;\que)=0$, so $M$ contains no horizontal surface. 
\end{proof} 
\vspace{-2pt}

\section{Mutation}	
\vspace{-2pt}

Let $\T = (B,A)$ be a knot in $S^3$ or a tangle.
Let $\T_0 = (B_0,A_0)$ be a  subtangle of $\T$ such that $B_0$ is a 3--ball.
Let $S_0 = \partial B_0$. 
Let $h\co B_0\to B^3$ be a homeomorphism such that $h(S_0\cap A) = Q\subset S^2$.
Let $\Gamma_0 \cong \zed_2\times\zed_2$ be the group of automorphisms of 
$(S^2,Q)$ consisting of rotations through $\pi$ about any one of the three 
co-ordinate axes in $\real^3$ together with the identity.
Let $\T_1 = \T - \Int \T_0$ and regard $\T$ as $\T_0 \cup \T_1$, where $\T_0$ 
is glued to $\T_1$ by the identity map on $S_0$. 
Now for $g\in\Gamma_0 - \{1\}$, let $\mu = h^{-1} gh\co S_0 \to S_0$ and 
define $\T' = \T_0 \cup_\mu \T_1$. 
We say $\T'$ is a {\em mutant} of $\T$.
The operation of replacing $\T$ by $\T'$ is {\em mutation} of $\T$ along 
$S_0$, by the {\em mutation involution} $\mu$.
\vspace{-4pt}

Boileau has asked \cite[Problem 1.69(c)]{Ki} if the unknotting number of a 
link is a mutation invariant. 
We prove that this is at least true for knots with unknotting number~1.
\vspace{-4pt}

\begin{thm}\label{thm:M}
Having unknotting number 1 is invariant under mutation.
\end{thm}
\vspace{-4pt}

\begin{proof}
Let $K$ be a knot with unknotting number 1, and let $K'$ be a mutant of $K$.
Then there is a Conway sphere $S$ which decomposes $K$ into two tangles 
$\T_1$ and $\T_2$, such that $K' = \T_1 \cup \rho(\T_2)$, where $\rho$ is 
rotation of the ball $B_2$ containing $\T_2$ through $\pi$ about one of 
the co-ordinate axes. 
Note that we also have $K' = \rho(\T_1) \cup \T_2$. 
If either $\T_1$ or $\T_2$ is trivial then $K= K'$. 
Also, since $u(K)=1$, $K$ is prime \cite{S1}.
Hence we may assume that the tangles $\T_i$ are prime, and therefore that 
$S$ is essential.
\vspace{-4pt}

First suppose that $K$ is not an EM--knot nor the union of two essential 
tangles, one of which is an EM--tangle. 
Then, taking $F=S$ in \fullref{cor4.2} we get that 
there is an unknotting arc $a$ for $K$ disjoint from $S$,
and therefore, without loss of generality, contained in $B_2$. 
As marked tangles, the crossing move  determined by $a$ transforms $\T_2$ 
to a rational tangle $\R$. 
Then the crossing move on $K'$ determined by $\rho(a)$ transforms 
$\rho (\T_2)$ to the rational tangle $\rho (\R) = \R$, and hence transforms 
$K'$ to $\T_1\cup \R=$ unknot.
\vspace{-4pt}

If $K$ is an EM--knot, then $K'=K$ \cite{E1}.
\vspace{-4pt}

Finally, suppose 
$K$ is the union of essential tangles $\P \cup \P_0$, where $\P_0$ is 
an EM--tangle. 
Then $\P_0$ is of the form $\S(\alpha,\beta;\gamma,*)$. 
If $S$ is not isotopic in $O(K)$ to $S_0 = \partial \M(\alpha,\beta)$, then 
the argument above shows that $K'$ has unknotting number~1.
So assume that $S= S_0$. 
Since rotating $\M(\alpha,\beta)$ about the horizontal axis leaves it 
invariant, we may assume that $\rho$ is rotation about the axis 
perpendicular to the plane of the paper.
Thus $K' = \P\cup \P'_0$, where $\P'_0 = \S(\beta,\alpha;\gamma,*)$.
\vspace{-4pt}

By part (3) of \fullref{thm:Main}, $K$ has an unknotting arc $a$ 
that lies in the 3--ball $B_0$ containing the tangle $\P_0$. 
By \fullref{lem:rotate}, there is a homeomorphism $h\co \S(\alpha,\beta;\gamma,*)
\to \S(\beta,\alpha;\gamma,*)$ such that $h|\partial D$ is rotation through
$\pi$ about the horizontal axis. 
Note that $h$ is isotopic to the corresponding rotation of the ball 
$\overline{S^3-D}$.
Hence there is a rotation $g$ of the ball $B_0$ which takes $\P_0$ to $\P'_0$.
The crossing change of $K$ determined by $a$ converts $\P_0$ to a rational
tangle $\R$, where $\P\cup \R$ is the unknot.
Therefore the crossing change of $K'$ determined by $g(a)$ converts 
$K' = \P\cup \P'_0 = \P \cup g(\P_0)$ to 
$\P \cup g(\R) = \P\cup \R = \text{unknot}$.
\end{proof}
\vspace{-4pt}

\section{Algebraic knots}	
\vspace{-4pt}

For the definition of an {\em algebraic\/} knot or link (in the sense 
of Conway) see Section~2 of \cite{Th1}. 
We briefly summarize this in a form suitable for our present purposes.
\vspace{-4pt}

\begin{figure}[ht!]
\centerline{\includegraphics[height=1.15truein]{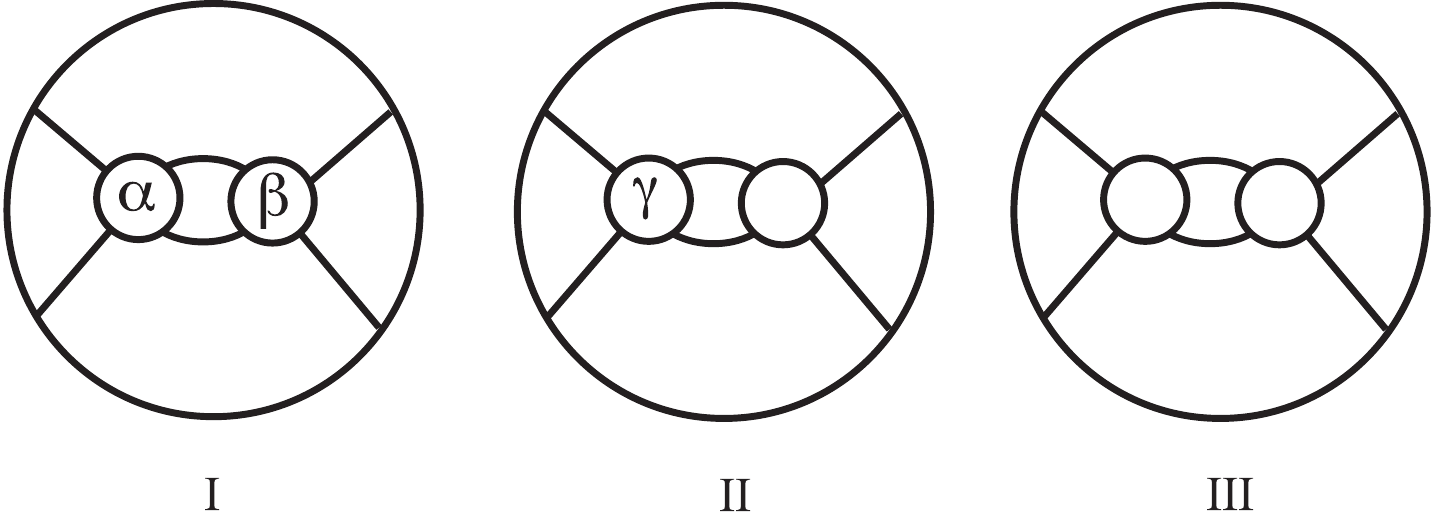}}
\caption{}		
\label{alg-tangle}
\end{figure}

An {\em elementary (algebraic) tangle\/} is a tangle of the form 
$\M(\alpha,\beta)$, $\M(\gamma,*)$, or $\M(*,*)$, where $\alpha,\beta,\gamma
\in \que-\zed$. 
We shall refer to these as elementary tangles of type~I, II, or III 
respectively. 
See \fullref{alg-tangle}.
\vspace{-5pt}

\begin{figure}[ht!]
\centerline{\includegraphics[height=0.9truein]{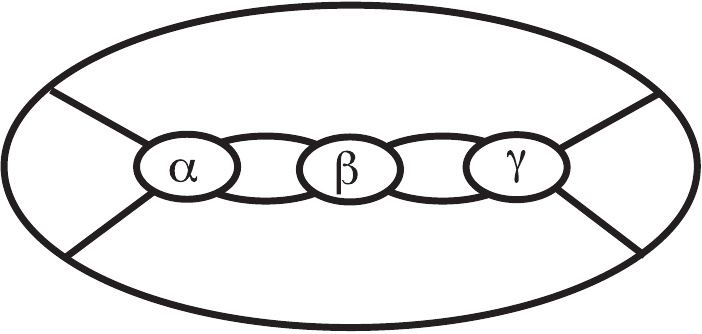}}
\caption{}		
\label{alg-tangle2}
\end{figure}

A {\em Montesinos tangle of length~3}, $\M(\alpha,\beta,\gamma)$, 
$\alpha,\beta,\gamma\in \que-\zed$, is defined in the obvious way; 
see \fullref{alg-tangle2}.
\vspace{-5pt}

Recall \cite{C} that if $\T$ is a marked tangle in $B^3$ then $1^*\T$ 
is the knot or link obtained by capping off $\T$ with the rational 
tangle $\R(0)$ (ie, $1^*\T$ is the numerator closure of $\T$). 
\vspace{-5pt}

Then an {\em algebraic\/} knot is a knot of one of the following forms:
\begin{itemize}
\item[{(a)}] $1^* \R(\alpha)$;
\item[{(b)}] $1^* \M(\alpha,\beta,\gamma)$;
\item[{(c)}] a union along boundary components of elementary tangles.
\end{itemize}
\vspace{-5pt}

The knots of type (a) are the 2--bridge or rational knots. 
(Any knot $1^*\M(\alpha,\beta)$ can also be expressed as $1^*\R(\gamma)$.) 
Those of type~(b) are the {\em Montesinos knots of length 3}. 
We will call a knot of type~(c) a {\em large} algebraic knot. 
Thus an algebraic knot is large if and only if it has an essential Conway 
sphere. 
If $K$ is an algebraic knot, of type~(a), (b), or (c), then the double 
branched cover of $K$ is a lens space, a Seifert fiber space over $S^2$ 
with three exceptional fibers, or a toroidal graph-manifold, respectively.
\vspace{-5pt}

An elementary tangle comes equipped with a marking, given by \fullref{alg-tangle}. 
In constructing a large algebraic knot $K$, the gluing homeomorphisms 
between the boundary components of the elementary tangles will not in 
general preserve the markings.
To describe $K$ as a marking-preserving union of marked tangles we need to 
interpolate marked tangles of 4--string braids in $S^2\times I$ between 
the boundary components. 
\vspace{-5pt}

This can also be described in terms of diagrams. 
\fullref{alg-tangle}~III is a diagram in a pair of pants of an elementary tangle 
of type~III. 
A diagram of an elementary tangle of type~I or II, in a disk or annulus 
respectively, may be obtained by inserting diagrams of the appropriate 
rational tangles into the diagrams in \fullref{alg-tangle}, I or II. 
Also, a 4--string braid in $S^2\times I$ has a diagram in an annulus. 
Then a knot is large algebraic if and only if it has a diagram in $S^3$ 
that is a union along boundary components of such elementary tangle 
diagrams and 4--string braid diagrams.

\begin{lem}\label{lemK}
$\R(p/q)$ can be transformed to $\R(1/0)$ by a crossing move if and only
if there exist coprime integers $r,s$ such that $p/q = \frac{2rs\pm1}{2s^2}$.
\end{lem}

\begin{proof}
If $p/q= \frac{2rs\pm 1}{2s^2}$, $q\ne0$, then $p/q$ has a continued 
fraction expansion of the form 
$[a_1,a_2,\ldots,a_k,\pm2,-a_k,\ldots,,-a_2,-a_1,a]$, 
(see \cite{KM} or \cite{Koh}), and hence can be transformed to $\R(1/0)$ 
by a crossing move. 

Conversely, suppose $\R(p/q)$ can be transformed to $\R(1/0)$ by a 
crossing move. 
Let the double branched covers of $\R(1/0)$ and $\R(p/q)$ be $V$ and $V'$ 
respectively. 
Then there is a knot $K$ in $V$ such that (with respect to some 
framing of $K$), 
$m/2$--Dehn surgery on $K$ gives $V'$. 
Note that, with respect to the basis of $H_1(\partial V)$ corresponding to the 
standard marking, the meridian of $V'$ has slope $p/q$.

If $K$ is unknotted in $V$, then $p/q = 1/0$, while if $K$ is a
core of $V$, then 
$p/q = m/2$; in both cases $p/q$ is of the stated form.

Otherwise, it follows from \cite[Theorem 2.4.4]{CGLS} 
that $K$ is an $(r,s)$--cable of the core of $V$. 
With respect to the usual framing on $K$, the Seifert fiber of the cable 
space $V- \Int N(K)$ has slope $rs$ on $\partial N(K)$. 
Hence, $K(m/2)$ will be a solid torus $V'$ if and only if 
$\Delta (m/2,\, rs/1) =1$, ie $m= 2rs\pm1$. 
The meridian of $V'$ then has slope $m/2s^2 = \frac{2rs\pm1}{2s^2}$.
\end{proof}

\begin{thm}\label{thm8.3}
Let $K$ be a large algebraic knot with unknotting number $1$. 
Then either
\begin{itemize}
\item[(1)] any unknotting arc for $K$ can be isotoped into either
\item[]{\rm(a)}\qua one of the rational tangles $\R(p/q)$ in an elementary tangle of 
type I; or 
\item[]{\rm(b)}\qua  the rational tangle $\R(p/q)$ in an elementary tangle of type II.
\end{itemize}
In case (a), the crossing move transforms $\R(p/q)$ to $\R(k/1)$ 
for some integer $k$, 
and $p/q = \frac{2s^2}{2rs\pm1} +k$, where $s\ge1$ and $(r,s)=1$.

\noindent 
In case (b), the crossing move transforms $\R(p/q)$ to $\R(1/0)$, and 
$p/q = \frac{2rs\pm1}{2s^2}$, where $s\ge1$ and $(r,s)=1$.
\begin{itemize}

\item[(2)] {\rm (a)}\qua $K$ is an EM--knot $K(\ell,m,n,p)$.

\item[{}] {\rm (b)}\qua
$O(K)$ has a unique connected incompressible 2--sided toric 2--suborbi\-fold $S$,
a Conway sphere, $K$ has an unknotting arc $a$ with $|a\cap S|=1$ (the
standard unknotting arc for $K(\ell,m,n,p)$) , and
$K$ has no unknotting arc disjoint from $S$.

\item[(3)] 
$K$ is the union of essential tangles $\P\cup \P_0$, where
$\P_0$ is an EM--tangle $\A_\ep (\ell,m)$ and $\partial\P_0$ is in $\bT(K)$.
Any unknotting arc for $K$ can be isotoped into $\P_0$.
The standard unknotting arc for $\A_\ep (\ell,m)$ is an unknotting arc
for $K$.
\end{itemize}
\end{thm}

\begin{remark}
The remarks (A), (B) following \fullref{thm:Main} also apply here.
\end{remark}
\vspace{-5pt}

\begin{proof} 
Note that the characteristic orbifold decomposition of $O(K)$ is gotten by 
amalgamating subcollections of the constituent elementary tangles. 
Applying \fullref{thm:Main}, we are left to check that 
\fullref{thm:Main}(1) implies \fullref{thm8.3}(1). 
\fullref{thm:Main}(1) implies that the unknotting move replaces  
$\R(p/q)$ in some elementary tangle of type~I or II with another 
tangle $\T$. 
$\T$ must be an integer tangle, $\R(k/1)$, if the elementary tangle is of 
type~I, and $\R(1/0)$ if it is of type~II. 
\fullref{lemK} gives  the desired result (using $\R(1/(\frac{p}{q}-k))$
for $\R(p/q)$ in type~I). 
\end{proof}
\vspace{-5pt}

\section{Some algebraic tangle calculations} 		
\vspace{-5pt}

In this section we do some calculations concerning crossing moves on certain 
algebraic tangles. 
These will be used in Sections~10 and 11.
\vspace{-5pt}

\begin{lem}\label{lem9.1}
Suppose $q>1$, and that $\M(p/q,\chi)$ is a rational tangle, where 
$\chi \in \que \cup\{\infty\}$. 
Then 
\begin{itemize}
\item[(1)] $\chi = k\in\zed$;
\item[(2)] $\M(p/q,k) = \R \left(\frac{kq+p}q\right)$
\item[(3)] if $\M (p/q,k) = \R (1/x)$, $x\in\zed$, then there exists 
$\ep = \pm1$ such that $x = \ep q$ and $kq+p=\ep$.
\end{itemize}
\end{lem}
\vspace{-5pt}

\begin{proof}
(1)\qua Since $\M(p/q,\chi)$ is a disk sum of $\R (p/q)$ and $\R(\chi)$, and 
$q>1$, we must have $\chi = k\in\zed$. 
\vspace{-5pt}

(2)\qua Incorporating the $k$ horizontal twists into $\R(p/q)$, we see 
that $\M(p/q,\chi) = \R \left( \frac{kq+p}q\right)$.
\vspace{-5pt}

(3)\qua This follows immediately from (2).
\end{proof}
\vspace{-5pt}

\begin{lem}\label{lem9.2}
Suppose $q_1,q_2>1$, and that $\M(p_1/q_1,\, p_2/q_2)$ can be transformed 
to a rational tangle $\R$ by a crossing move. 
Then 
\begin{itemize}
\item[(1)] the crossing arc is isotopic to an $(r,s)$--cable, $s\ge1$, of one 
of the two exceptional fibers of $\M(p_1/q_1,\, p_2/q_2)$; 
\vspace{-5pt}
\item[(2)] the crossing move transforms the corresponding rational tangle, 
$\R (\frac{p_1}{q_1})$, say, of $\M$, to an integral tangle $\R(k)$;
\vspace{-5pt}
\item[(3)] $\R = \R \left( \frac{kq_2 +p_2}{q_2}\right)$;
\vspace{-5pt}
\item[(4)] $p_1/q_1 = k+ \frac{2s^2}{2rs\pm 1}$.
\end{itemize}
\end{lem}

\begin{proof}
(1)\qua 
The argument is very similar to the proofs of \fullref{lem6.1} and 
\fullref{thm:JSJ}. 
Let $M$ be the double branched cover of $B^3$ along 
$\M(\frac{p_1}{q_1}, \frac{p_2}{q_2})$. 
$M$ is a Seifert fiber space over the disk with two exceptional fibers. 
The crossing move corresponds to replacing a marked tangle $\T_0$ in 
$\M(\frac{p_1}{q_1}, \frac{p_2}{q_2})$  
with a marked tangle $\T'_0$, resulting in a rational tangle. 
Let $V_0$ be the solid torus preimage of $\T_0$ in $M$ and 
$X = M-\Int (V_0)$ be its exterior. 
Let $\gamma$ be the meridian of $V_0$. 
Then $M= X(\gamma)$ and $X(\mu)$ is a solid torus for some $\mu$ with 
$\Delta (\gamma,\mu)=2$. 
By \cite[Proposition 9]{MZ}, $X$ must be either Seifert fibered or toroidal. 
In the first case, $X$ is the exterior of an exceptional fiber in $M$ 
and we are done. 
So $X$ is toroidal and we let $W$ be the component of $X$ cut along 
$\bT(X)$ (canonical torus decomposition) that contains $\partial X$. 
Then $\partial W -\partial X$ is compressible in $W(\gamma),W(\mu)$. 
Furthermore $W(\gamma),W(\mu)$ are irreducible. 
Then $W(\gamma),W(\mu)$ are solid tori. 
By \cite[Theorem 2.4.4]{CGLS}, $W$ is a cable space. 
In particular, say that $W$ is the exterior of the $(p,q)$ curve in the 
solid torus $W(\mu$), $q\ge2$. 
If $\gamma_0,\mu_0$ are the slopes of the meridian disks on 
$\partial W(\gamma)$, $\partial W(\mu)$ respectively, 
then $\Delta (\gamma_0,\mu) = q^2 (\Delta \gamma,\mu)\ge8$. 
Thus, if $X_0 = X-W$, we may argue as above to conclude that $X_0$ is 
the exterior of a $(p,q)$--cable on some knot $k_0$ in $X(\mu)$. 
But then a coordinate calculation, as in the proof of \fullref{thm:JSJ} 
(Case~IIB), says that $\gamma_0 = \frac{n_1p_1q_1\pm1}{n_1}$, $n_1\ge2$ 
and $\gamma_0 = \frac{npq\pm1}{nq^2}$ for $n\ge2$. 
This contradiction finishes the proof of \fullref{lem9.2}(1).

\fullref{lem9.2}(2) and (3) 
now follow from \fullref{lem9.1}, (1) and (2). 
Finally, (4) follows from \fullref{lemK} applied to 
$\R(1/(\frac{p_1}{q_1} -k))$. 
\end{proof}

\begin{cor}\label{newcor:9.3}
Suppose $q_1,q_2>1$ and that $\M (p_1/q_1,\, p_2/q_2)$ can be transformed 
to a vertical twist tangle $\R(1/x)$ by a crossing move.
Then there exist $\ep = \pm1$ and $k\in\zed$ such that, after possibly 
interchanging $p_1/q_1$ and $p_2/q_2$,
\begin{itemize}
\item[(1)] $x= \ep q_2$;
\item[(2)] $kq_2 + p_2 =\ep$;
\item[(3)] $p_1/q_1 = k+\frac{2s^2}{2rs\pm1}$, for some $s\ge1$, 
$(r,s)=1$.
\end{itemize}
\end{cor}

\begin{proof} 
This follows from Lemmas~\ref{lem9.1} and \ref{lem9.2}.
\end{proof}

\begin{lem}\label{lem9.4}
$1^*\M(\alpha,\beta)$ is the unknot if and only if $\Delta (\alpha,-\beta)=1$.
\end{lem}

\begin{proof}
$1^* \M(\alpha,\beta)$ is the unknot if and only if its double branched 
cover $M$ is $S^3$. 
But $M$ is the union of the two solid tori $\tilde{\R} (\alpha)$ 
and $\tilde{\R}(\beta)$, 
whose meridians have slopes $\alpha$ and $-\beta$ respectively on the torus
$T = \partial \tilde{\R} (\alpha)$ with respect to the basis of $H_1(T)$ 
corresponding to the lifts of the slopes $1/0$ and $0/1$ on 
$\partial \R(\alpha)$. 
Hence $M$ is $S^3$ if and only if $\Delta (\alpha,-\beta)=1$.
\end{proof}

\begin{lem}\label{newlem:9.5}
Let $K$ be the knot shown in \fullref{fig-new8-1}, where 
$q_1,q_2>1$ and $\T$ is some marked tangle. 
Then 
\begin{itemize}
\item[(1)] $K$ is the unknot if and only if 
$\T = \R(x)$, where $x\in\zed$ satisfies
$$xq_1 q_2 + p_1 q_2 + p_2 q_1 = \pm1\ ;$$
\item[(2)] if $|p_1/q_1|,|p_2/q_2| <1$, and $x$ is as in {\rm (1)},
then $|x| \le 1$.
\end{itemize}
\begin{figure}[ht!]
\labellist\small
\pinlabel $\left(\frac{p_1}{q_1},\frac{p_2}{q_2}\right)$ <-2pt,0pt> [l] at 174 620
\endlabellist
\centerline{\includegraphics[height=1.3truein]{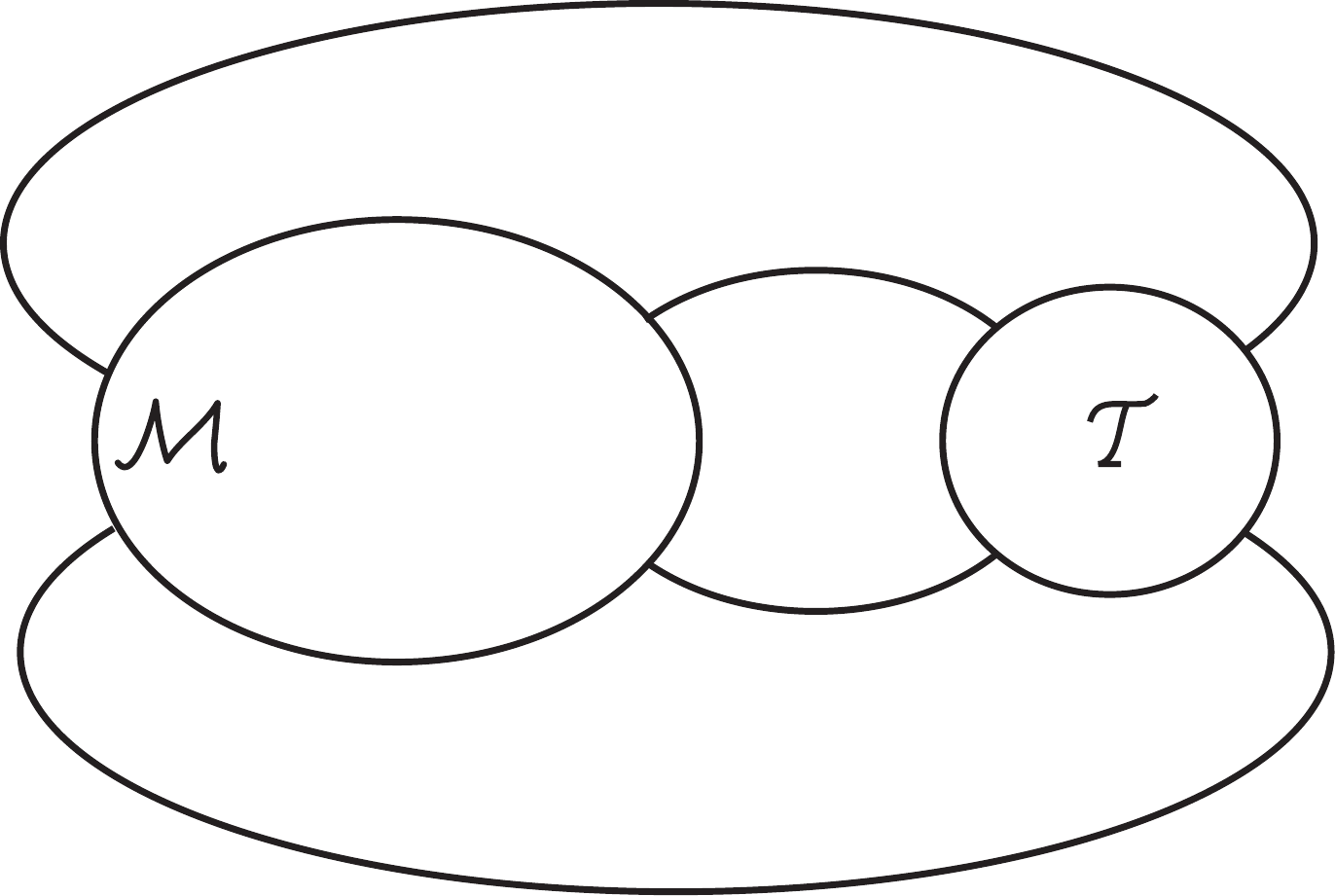}}
\caption{}
\label{fig-new8-1}
\end{figure}
\end{lem}
\vspace{5pt}

\begin{proof}
Suppose $K$ is the unknot. 
Passing to double branched covers, we see that $\tilde{\T}$ must be 
a solid torus, implying that $\T$ is a rational tangle.
Moreover, the meridian of $\tilde \T$ must be distance~1 from the 
Seifert fiber of $\tilde \M (p_1/q_1,\, p_2/q_2)$. 
Since the latter projects to the slope $1/0$ on $\partial \M(p_1/q_1,\,
p_2/q_2)$, this implies that $\T= \R(x)$ for some $x\in\zed$. 
Incorporating this twist tangle $\R(x)$ with $\R(p_1/q_1)$, we see 
that the unknot is the union of the rational tangles 
$\R(x+ \frac{p_1}{q_1}) = \R\left(\frac{xq_1 +p_1}{q_1}\right)$ and 
$\R (p_2/q_2)$. 
Then  $\Delta \left(\frac{xq_1+p_1}{q_1}, \frac{-p_2}{q_2}\right)=1$ 
by \fullref{lem9.4}, giving the equation in (1).

Conversely, if $\T = \R(x)$ where $x$ satisfies the given equation then 
$\Delta (\frac{xq_1+p_1}{q_1},  -\frac{p_2}{q_2}) =1$ and 
$K$ is the unknot by \fullref{lem9.4}. 

To prove (2), suppose $|p_1/q_1|, |p_2/q_2| <1$. 
Then from (1) we have 
$$x+ \frac{p_1}{q_1} + \frac{p_2}{q_2} = \pm \frac1{q_1q_2}\ ,$$
giving $|x| \le \frac{(q_1-1)}{q_1} + \frac{(q_2-1)}{q_2} +
\frac1{q_1q_2} <2$.
\end{proof}

\begin{lem}\label{lem9.6}
Suppose $\alpha,\beta,\gamma\in\que -\zed$, and 
$|\alpha|,|\beta| <1$. 
Then $\S(\alpha,\beta;\gamma,\T)$ is not the unknot, for any tangle $\T$.
\end{lem}

\begin{proof} 
Suppose $\S (\alpha,\beta;\gamma,\T)$ is the unknot. 
Then passing to double branched coverings as usual we see that $\T$ 
must be a rational tangle $\R(\chi)$.
Furthermore, by \fullref{newlem:9.5}, $\M(\gamma,\chi)$ must be a 
vertical twist tangle $\R(-1/x)$, where $|x|\le1$ since $|\alpha|,|\beta|<1$.
But by \fullref{lem9.1}(3), $|x|\ge 2$.
\end{proof}
\vspace{-5pt}

\begin{cor}\label{cor9.7}
Let $K = \S(\alpha,\beta;\gamma,\delta)$, where $\alpha,\beta,\gamma,\delta
\in\que-\zed$, and $|\alpha|,|\beta|<1$. 
Then $K$ cannot be unknotted by a crossing move in $\M(\gamma,\delta)$.
\end{cor}
\vspace{-5pt}

\begin{proof}
By \fullref{lem9.2}(2), the crossing move has the effect of 
replacing  one of the rational tangles $\R (\gamma)$ or $\R(\delta)$ 
in $\M(\gamma,\delta)$ with some other (rational) tangle. 
But this contradicts Lemmas~\ref{lem9.6} and \ref{D8}(2). 
\end{proof}
\vspace{-6pt}

Recall (\fullref{cor:EMknot}) that the EM--knots are all of the form 
$\S(\alpha,\beta;\gamma,\delta)$ with $\alpha,\beta,\gamma,\delta$ $\in 
\que-\zed$ and $|\alpha|,|\beta|,|\gamma|,|\delta| <1$.
\vspace{-6pt}

\begin{thm}\label{thm9.8}
Let $K\! = \S(\alpha,\beta;\gamma,\delta)$ where $\alpha,\beta,\gamma,\delta 
\in \que-\zed$ and $|\alpha|,|\beta|,|\gamma|,|\delta|\allowbreak <1$.
Then $K$ has unknotting number~1 if and only if $K$ is an EM--knot 
$K(\ell,m,n,p)$.
\end{thm}
\vspace{-6pt}

\begin{proof}
$K$ has a unique essential Conway sphere $S= \partial\M (\alpha,\beta) 
= \partial \M(\gamma,\delta)$. 
Therefore, by \fullref{cor4.2}, if $K$ has unknotting number~1 then 
either $K$ is an EM--knot or $K$ can be unknotted by a crossing move
disjoint from $S$.
But the latter is impossible by \fullref{cor9.7}.
\end{proof}
\vspace{-6pt}

\section{Examples}		
\vspace{-6pt}

In this section we apply our results to certain families of knots 
defined in terms of the notation of Conway \cite{C}; since that notation 
naturally encodes the characteristic toric orbifold decomposition of a 
knot it is eminently suited to our techniques. 
In particular, we consider all the knots up to 11 crossings in Conway's 
tables \cite{C} with the property that their description in the tables 
makes it clear that they contain an essential Conway sphere.
It turns our that these are all large algebraic knots. 
Specifically, they are the knots that are listed in \cite{C} as either 
$.a.b$, $.a.b.c$, $.a.b.c.d$, $.a.(b,c)$, $.(a,b).c$, $(a,b)(c,d)$, 
$(a,b)1(c,d)$, or $a,b,c,d$.
We determine exactly which of them have unknotting number~1.
Note that the knots $a,b,c,d$ are Montesinos knots of length~4, so they 
have unknotting number greater than 1 by \cite{Mot}. 
\vspace{-6pt}

Throughout this section, $a,b,c$ and $d$ will denote rational numbers.
\vspace{-6pt}

We start with the knots $.a.b.c.d$.
Recall \cite[page 335]{C} that in Conway's tables the form $.a.b.c.d$ is used 
only when $a,b,c$ and $d$ are $>0$. 
It turns out that all the EM--knots $K(\ell,m,n,p)$ are of this form 
(up to mirror image), and that a knot $.a.b.c.d$ with $a,b,c,d>0$ has 
unknotting number~1 if and only if it is an EM--knot.
Recall also that in \cite{C} $.a.b.c.1$ is abbreviated to $.a.b.c$, and 
$a.b.1$ to $.a.b$~.

\begin{figure}[ht!]

\centerline{\includegraphics[height=1.0truein]{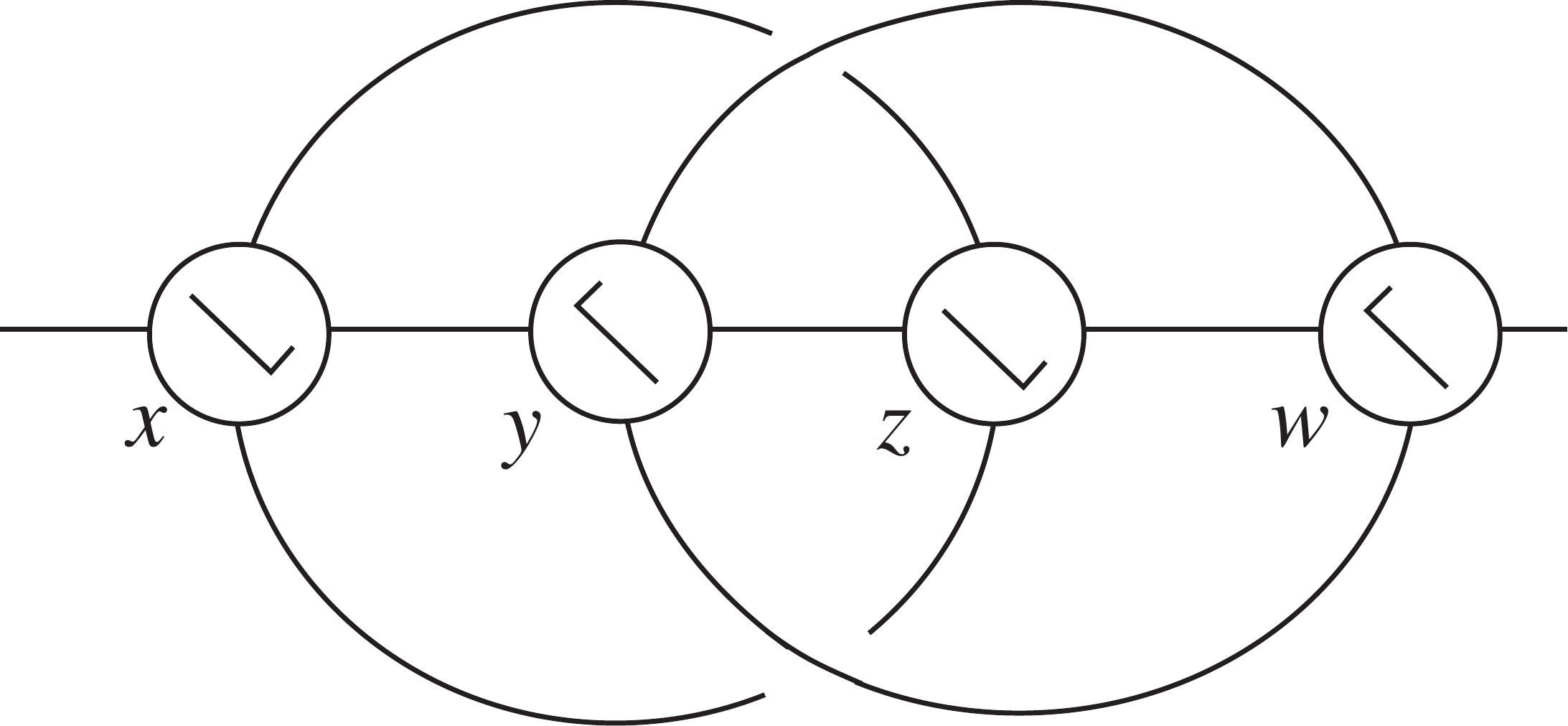}}
\caption{} 
\label{maskedtangles}
\end{figure}

First  we describe some symmetries of Conway's $.x.y.z.w$ notation. 
Recall that if $x,y,z$ and $w$ are arbitrary marked tangles, then 
$.x.y.z.w$ is the knot shown in \fullref{maskedtangles} (where the 
leftmost and rightmost horizontal arcs are understood to meet at the 
point at infinity in the projection $S^2$).

Let $D_8$ be the dihedral group of order 8, the group of symmetries 
of the square. 
By cyclically numbering the vertices of the square 1,2,3,4, we regard $D_8$ 
as a subgroup of $S_4$. 
Then $D_8$ acts on the set of expressions of the form $.x.y.z.w$ by 
permuting the substituent tangles; thus, by a slight abuse of notation, 
we write 
$\pi (.x_1.x_2.x_3.x_4) = .x_{\pi (1)}.x_{\pi(2)} .x_{\pi(3)}. x_{\pi(4)}$, 
for $\pi \in D_8$.

Recall also \cite[pages 330--331]{C} that $\zed_2 \times\zed_2$ acts on the 
set of marked tangles as follows. 
If $x$ is a marked tangle, then $x_h,x_v$ and $x_r$ are the marked 
tangles obtained by rotating $x$ through $180^\circ$ about the horizontal 
axis, the vertical axis, and the axis perpendicular to the plane of the paper,
respectively. 
Let $\Mut \cong \zed_2\times\zed_2$ be the group $\{h,v,r,\id\}$. 
(Then mutation is the equivalence relation on knots generated by 
replacing a tangle $x$ in (some diagram of) $K$ by $x_f$ for some 
$f\in\Mut$.) 
Again by a slight abuse of notation, we write  $(.x.y.z.w)_f  = 
.x_f.y_f.z_f.w_f$, for $f\in\Mut$.

Let $\mu\co D_8\to\Mut$ be the epimorphism defined by $\mu((1\,2\,3\,4))=h$, 
and $\mu((14)(23)) =v$. 
Finally, let $\ep \co \Mut \to\zed_2 = \{\pm1\}$ be the homomorphism defined 
by $\ep (h) = \ep (v) =-1$, and recall that -- denotes mirror-image.

\begin{thm}\label{thm9.1}
If $\pi \in D_8$ then:
$$\pi (.x.y.z.w) = \ep (\mu (\pi)) (.x.y.z.w)_{\mu (\pi)}$$
\end{thm}

\begin{cor}\label{cor9.2}
Up to mutation and mirror-image, $.x.y.z.w$ is invariant under the 
action of $D_8$.
\end{cor}

Since $t_f =t$ for rational tangles $t$, for all $f\in\Mut$, we have
\vspace{-3pt}

\begin{cor}\label{cor9.3} 
If $\pi \in D_8$ then: 
$$\pi (.a.b.c.d) = \ep (\mu (\pi)) .a.b.c.d$$
\end{cor}
\vspace{-3pt}

\begin{proof}[Proof of \fullref{thm9.1}] 
Consider $.x.y.z.w$, as shown in \fullref{maskedtangles}.
By sliding the tangle $x$ around the point at infinity we get 
\fullref{slidingtangle-B}.
Changing all crossings, we see that $.x.y.z.w = -.y_h.z_h.w_h.x_h$. 
This shows that the theorem holds for $\pi = (1\,2\,3\,4)$.
\vspace{-3pt}

\begin{figure}[ht!]

\centerline{\includegraphics[height=1.0truein]{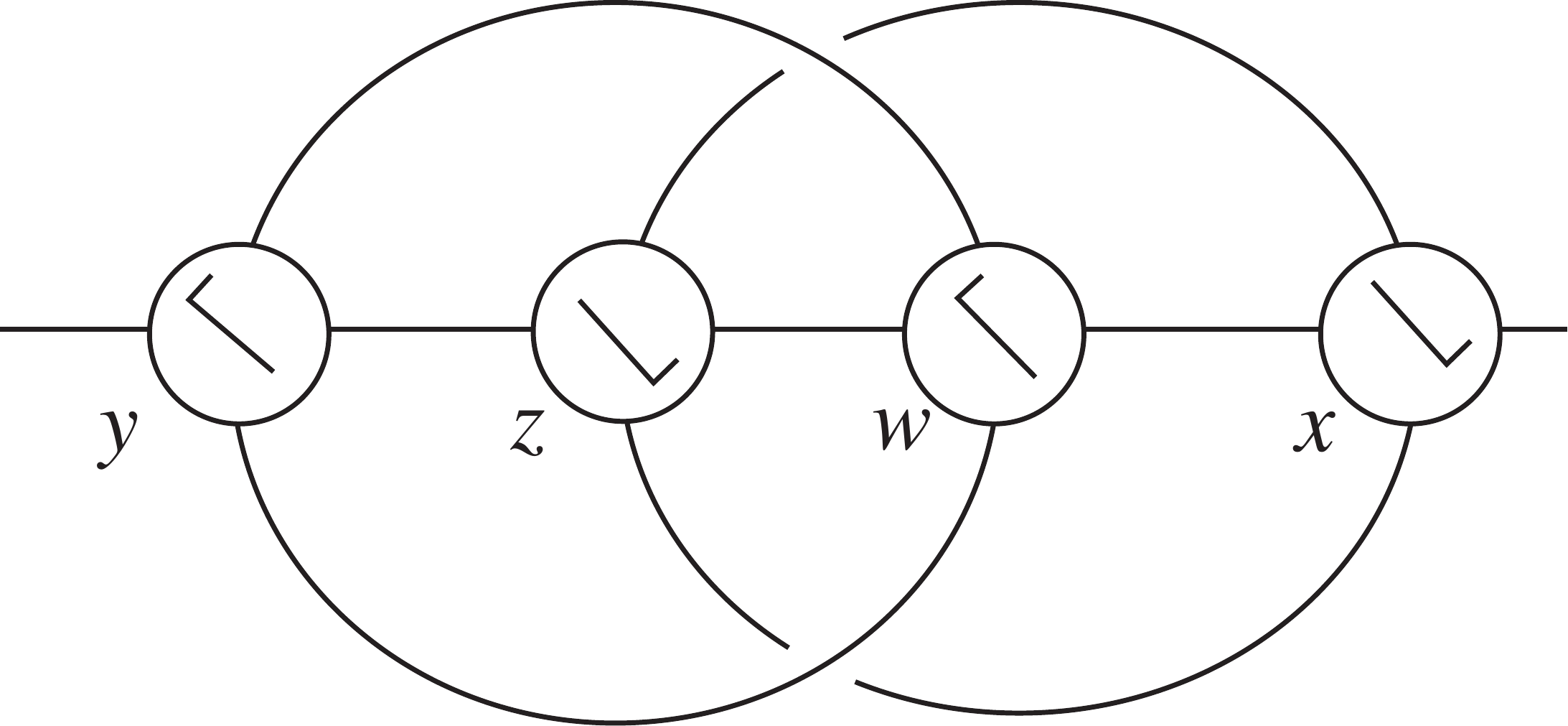}}
\caption{} 
\label{slidingtangle-B}
\end{figure}

Rotating \fullref{maskedtangles} through $180^\circ$ about the 
central axis perpendicular to the  plane of the paper, we get 
\fullref{rotatingtangles}.
Changing all crossings shows that $.x.y.z.w = - .w_v.z_v.y_v.x_v$, 
in other words, the theorem holds for $\pi = (14)(23)$.
\vspace{-3pt}

\begin{figure}[ht!]

\centerline{\includegraphics[height=1.0truein]{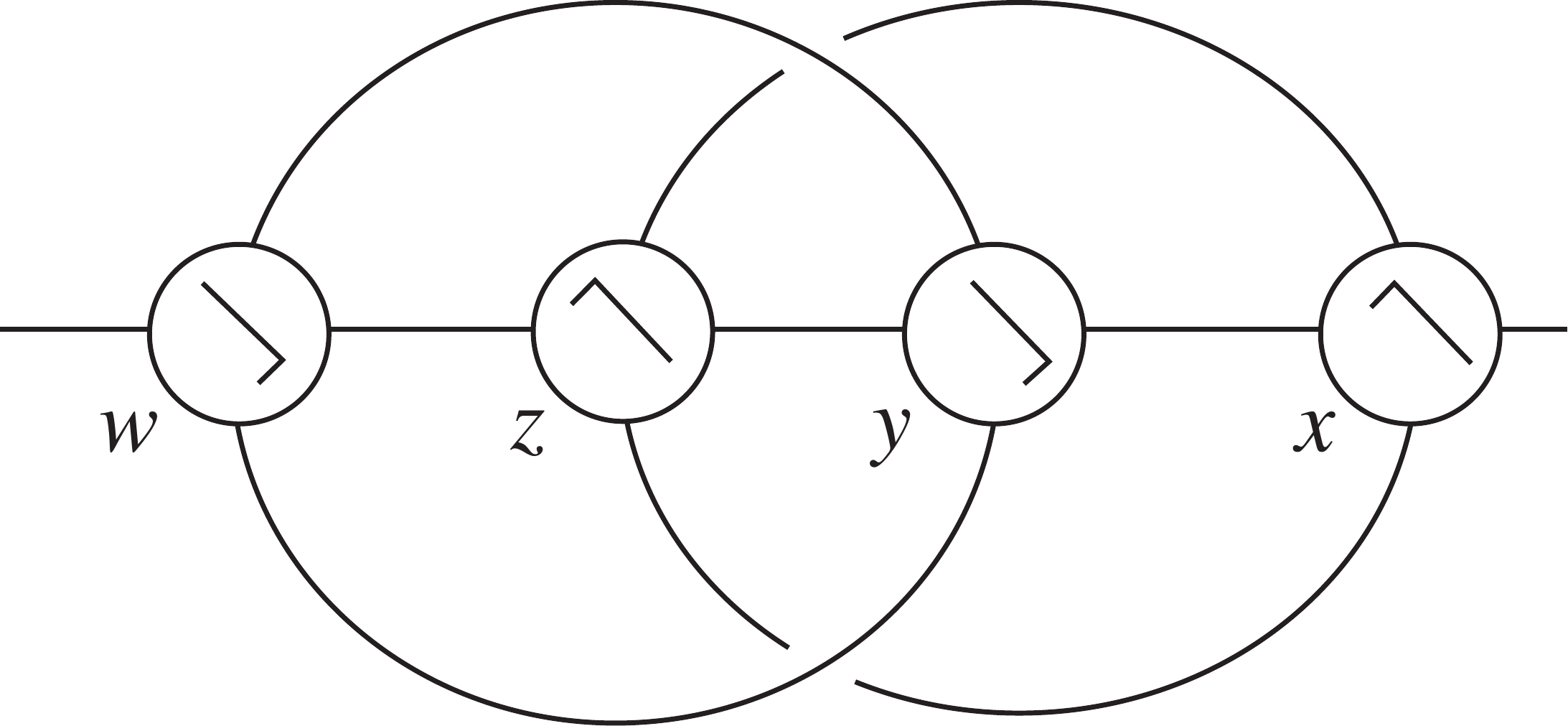}}
\caption{} 
\label{rotatingtangles}
\end{figure}

Since $(1\,2\,3\,4)$ and $(14)(23)$ generate $D_8$, the result follows.
\end{proof}
\vspace{-3pt}

We consider a further symmetry. 
Recall \cite[page 331]{C} that $t0$ denotes the reflection of the tangle $t$ 
in a plane perpendicular to the paper through the NW/SE--diagonal.
Note that $(t0)_r$ is then the reflection of $t$ in a plane through the 
NE/SW--diagonal.

\begin{thm}\label{thm:9.4}
$.x.y.z.w = - .x0.(y0)_r.z0.(w0)_r$
\end{thm}

In particular, for rational tangles we have 

\begin{cor}\label{cor:9.5}
$.a0.b0.c0.d0 = -.a.b.c.d$
\end{cor}

\begin{proof}
$t00 =t$ and $a_r =a$ when $a$ is rational.
\end{proof}

\begin{cor}\label{cor:9.6}\
\begin{itemize}
\item[(1)] $.a.b.c0 = -(.a0.b0.c)$
\item[(2)] $.a.b = -(.a0.b0)$
\item[(3)] $.a0.b = -(.a.b0)$
\end{itemize}
\end{cor}

\begin{proof}
These follow immediately from \fullref{cor:9.5}, along with the facts
that $10=1$ and $t00=t$.
\end{proof}

\begin{figure}[ht!]

\centerline{\includegraphics[height=1.0truein]{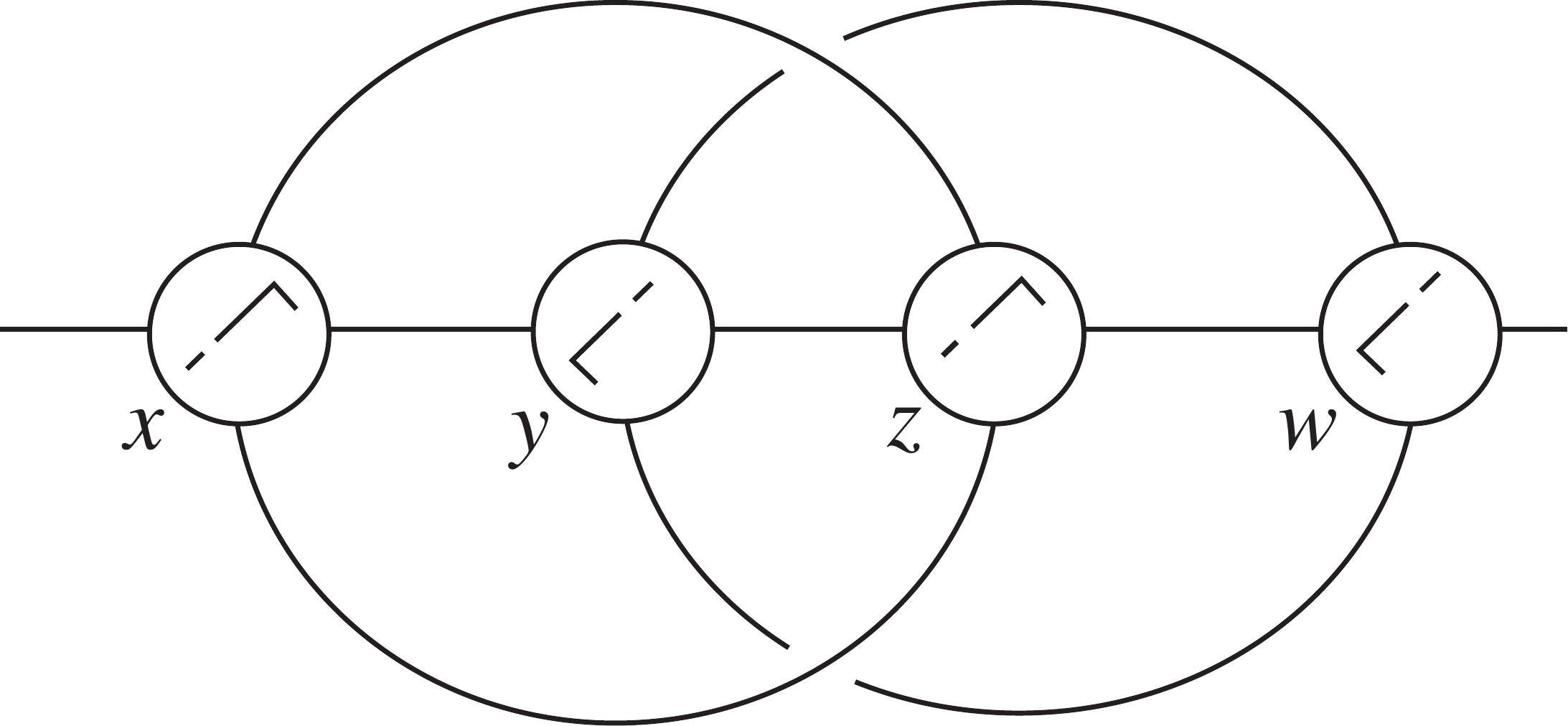}}
\caption{} 
\label{rotating180-B}
\end{figure}

\begin{proof}[Proof of \fullref{thm:9.4}]
Rotating \fullref{maskedtangles} through $180^\circ$ about the 
horizontal axis gives \fullref{rotating180-B}.
Changing all crossings, we now get $.x0.(y0)_r.z0.$ $(w0)_r$.
\end{proof}

The following lemma describes $.a.b.c.d$ in terms of the square tangle $\S$.

\begin{lem}\label{lem:9.7} 
$.a.b.c.d = \S (\frac{-1}{c+1},\frac{a}{a+1};\frac1{b+1},\frac{-d}{d+1})$
\end{lem}

\begin{proof}
This follows from the second deformation shown in \cite[Figure 10]{C}, 
possibly together with \fullref{D8}. 
\end{proof}

Recall that if $K$ is an EM--knot $K(\ell,m,n,p)$, then, by taking the 
mirror image of $K$ if necessary, we may assume that $\ell>1$.

\begin{lem}\label{lem:9.8}
Assume $\ell>1$. 
Then $K(\ell,m,n,p) = .a.b.c$, where $a,b,c$ are $>0$ and are given by:
\begin{equation*}
\begin{split}
p=0\ :\quad & a= \frac{m}{(\ell-1)m-1}\ ,\quad 
b= \frac{2mn - m+n}{2mn-m-n+1}\ ,\quad c=\ell-1\\
\noalign{\vskip6pt}
n=0\ :\quad & a= \frac{2mp-m-p}{(\ell-1)(2mp-m-p)-2p+1}\ ,\quad 
b = \frac{m}{m-1}\ ,\quad c= \ell-1
\end{split}
\end{equation*}
\end{lem}

\begin{proof}
This follows from Lemmas~\ref{lem:EMknot}  and \ref{lem:9.7}.
\end{proof}

\begin{thm}\label{thm:9.9}
Let $K= .a.b.c.d$, with $a,b,c,d>0$. 
Then $K$ has unknotting number 1 if and only if it is an EM--knot
$K(\ell,m,n,p)$.
\end{thm}
\vspace{-5pt}

\begin{proof}
By \fullref{lem:9.7}, $.a.b.c.d = \S (\alpha,\beta;\gamma,\delta)$ 
where $\alpha,\beta,\gamma,\delta \in \que -\zed$ and 
$|\alpha|,|\beta|,|\gamma|,|\delta| <1$. 
The result is now an immediate consequence of \fullref{thm9.8}.
\end{proof}
\vspace{-5pt}

We now determine the EM--knots up to 11 crossings. 
\vspace{-5pt}

\begin{thm}\label{thm:9.10} 
The EM--knots $K(\ell,m,n,p)$ with at most 11 crossings are listed 
in Table~EM 11, up to mirror image.
\end{thm}
\vspace{-5pt}

\vbox{
$$\vbox{\offinterlineskip\halign{
\quad\strut\hfil$#$\hfil\quad
&\vrule#&\quad$#$\hfil\quad
&\vrule#&\quad\hfil$#$\quad\cr
\text{Rolfsen}&&\text{Conway}&&K(\ell,m,n,p)\cr
\noalign{\hrule}
8_{17}&&.2.2&&K(3,2,0,1)\cr
9_{33}&&.21.2&&-K(2,3,0,1)\cr
10_{82}&&.4.2&&K(2,2,0,2)\cr
10_{84}&&.22.2&&K(2,2,0,-1)\cr
10_{88}&&.21.21&&K(2,3,0,0)\cr
10_{91}&&.3.2.20&&-K(4,1,1,0)\cr
10_{95}&&.210.2.2&&K(3,2,0,0)\cr
&&.311.2&&-K(2,2,-1,0)\cr
&&.23.2&&K(2,2,0,3)\cr
&&.212.2&&K(2,2,0,-2)\cr
&&.2111.2&&-K(2,2,2,0)\cr
&&.31.21&&-K(2,-3,1,0)\cr
&&.22.2.20&&-K(3,-2,1,0)\cr
&&.210.21.2&&K(3,2,1,0)\cr
\noalign{\hrule} }}$$
\centerline{$K(\ell,m,n,p)$ with at most 11 crossings}
\medskip
\centerline{{Table} EM 11}
}
\vspace{-5pt}

\begin{remark}
The third column of Table EM 11 represents the knot in $K(\ell,m,n,p)$ form.
These representations are not unique.
\end{remark}
\vspace{-5pt}

\begin{proof}[Proof of \fullref{thm:9.10}]
By possibly taking mirror images, Proposition~1.4 of \cite{E1} allows us 
to assume that $\ell>1$. 
By \fullref{lem:9.8}, $K(\ell,m,n,p) = .a.b.c$ with $a,b,c$ positive 
rational numbers. 
Inserting alternating diagrams of the corresponding rational tangles 
$\R(a)$, $\R(b)$, $\R(c)$ into 
Conway's $6^{**}$ polyhedron gives an alternating diagram of $.a.b.c$.
This alternating diagram will be a minimal crossing diagram from which 
we can compute the crossing number of $.a.b.c$.
To get  alternating diagrams of the rational tangles we can use 
the positive continued fraction expansions of $a,b$, and $c$. 
That is, if $r$ is a positive rational number and 
$r= [a_1,a_2,\ldots,a_n]$ (see Section~2) with $a_i>0$ if $i\ne n$ 
and $a_n\ge 0$, then there is an alternating diagram  of the rational 
tangle $\R(r)$ that has exactly $\sum_{i=1}^n a_i$ crossings. 
If $r\ne1$, we will take $a_1>1$.

In Tables E1--E3 below, we list nonnegative continued fraction expansions 
of $a,b,c$.
\medskip

\vbox{
$$\vbox{\offinterlineskip\halign{
\quad\strut\hfil$#$\hfil\quad
&\vrule#&\quad\hfil$#$\hfil\quad
&\vrule#&\quad\hfil$#$\hfil\quad
&\vrule#&\quad\hfil$#$\hfil\cr
m,p&&a&&b&&c\cr
\noalign{\hrule}
m>0,\ p>0&&[p-1,2,m-2,1,\ell-2,0]&&[m-1,1]&&\ell-1\cr
m>0,\ p<0&&[|p|,1,1,m-2,1,\ell-2,0]&&[m-1,1]&&\ell-1\cr
m<0,\ p>0&&[p-1,1,1,|m|,\ell-1,0]&&[|m|,1,0]&&\ell-1\cr
m<0,\ p<0&&[|p|,2,|m|,\ell-1,0]&&[|m|,1,0]&&\ell-1\cr}}$$
$$n=0\quad [\ell >1,\quad m\ne0,1,\quad (\ell,m,p)\ne (2,2,1)]$$
\centerline{{Table} E1}
}

\bigskip
\vbox{
$$\vbox{\offinterlineskip\halign{
\quad\strut\hfil$#$\hfil\quad
&\vrule#&\quad\hfil$#$\hfil\quad
&\vrule#&\quad\hfil$#$\hfil\quad
&\vrule#&\quad\hfil$#$\hfil\cr
m,n&&a&&b&&c\cr
\noalign{\hrule}
m>0,\ n>0&&[m-1,1,\ell-2,0] &&[n-1,1,1,m-1,1]&&\ell-1\cr
m>0,\ n<0&&[m-1,1,\ell-2,0] &&[|n|,2,m-1,1]&&\ell-1\cr
m<0,\ n>0&&[|m|,\ell-1,0]&&[n-1,2,|m|-1,1,0]&&\ell-1\cr
m<0,\ n<0&&[|m|,\ell-1,0] &&[|n|,1,1,|m|-1,1,0]&&\ell-1\cr}}$$
$$p=0\quad [\ell>1,\quad m\ne0,\quad (\ell,m)\ne (2,1);
\ (m,n)\ne (1,0),(-1,1)]$$
\centerline{{Table} E2}
}

\bigskip
\vbox{
$$\vbox{\offinterlineskip\halign{
\quad\strut\hfil$#$\hfil\quad
&\vrule#&\quad\hfil$#$\hfil\quad
&\vrule#&\quad\hfil$#$\hfil\quad
&\vrule#&\quad\hfil$#$\hfil\cr
m&&a&&b&&c\cr
\noalign{\hrule}
m>0&&[m-1,1,\ell-2,0]&&[m-1,1]&&\ell-1\cr
m<0&&[|m|,\ell-1,0]&&[|m|,1,0]&&\ell-1\cr}}$$
$$n=0=p\qquad (\ell>1;\ m\ne 0,1)$$
\centerline{{Table} E3}
}

If $0$ appears in any but the last entry of one of these expansions 
(for certain $\ell,m,n,p$), we may use one of the following rules to 
eliminate it:
\begin{itemize}
\item[(1)] $[0,a,b,c,\cdots]\longrightarrow [b,c,\cdots]$
\item[(2)] $[\cdots a,b,0,c,d,\cdots]\longrightarrow [\cdots a,b+c,d,\cdots]$
\end{itemize}
Note that the sum of the entries is changed by (1) only and that amounts 
to deleting the second entry from the sum.
In this way we enumerate those $K(\ell,m,n,p)$ with crossing number at
most~11.

As an example, consider the case when $n=0$, $m>0$, $p>0$. 
Assuming $p>1$, Table~E1 gives the crossing number of $K(\ell,m,0,p)$ as 
$2m+2\ell+p$ (the crossings from $\R(a)$, $\R(b)$, $\R(c)$ 
plus 3 from the Conway polyhedron 
$6^{**}$). 
Thus if $K(\ell,m,n,p)$ has at most 11 crossings, 
$(\ell,m,p) \in \{	
(2,2,2), (2,2,3)  	
\}$ (noting that $m\ne 1$ when $n=0$). 
If $p=1$ and $m>2$, then the crossing number of $K(\ell,m,0,1)$ 
is $2m+2\ell-1$.
Thus $(\ell,m) \in \{(2,3), (2,4), (3,3)\}$.
Finally, if $p=1$ and $m=2$, then $\ell>2$ and the crossing number 
is $2\ell+2$. 
That is, $\ell\in \{3,4\}$.

Once having enumerated the $K(\ell,m,n,p)$ with at most 11 crossings and 
written them in the form $.a.b.c$, we can now locate them in the tables. 
To do this we use the symmetries given by Corollaries~\ref{cor9.3} 
and \ref{cor:9.6} (note that if $t= [a_1,\ldots,a_n]$ is rational 
then $t0 = [a_1,a_2,\ldots,a_n,0]$).

This completes the proof of \fullref{thm:9.10}.
\end{proof}
\medskip

Of the knots listed in Conway's tables in the form $.a.b.c.d$ 
(or $.a.b$, or $.a.b.c$) that are not EM--knots, there are: 
1 with 8 crossings, 
3 with 9 crossings, 
13 with 10 crossings, and 
45 with 11 crossings. 
By \fullref{thm:9.9} these all have unknotting number greater than~1.

Next we consider the knots of the form $(a,b)(c,d)$ in Conway's notation. 
For these we first  have the following result.

\begin{thm}\label{thm:9.11} 
Let $K= (a,b)(c,d)$, where $|a|,|b|,|c|$ and $|d|$ are $>1$ and 
either $ab>0$ or $cd>0$. 
Then $K$ does not have unknotting number~1.
\end{thm}

\begin{proof}
It is easy to see (possibly using \fullref{D8}) that 
$(a,b)(c,d)\! =\!\S(-\frac1a,-\frac1b;$ $ \frac1c,\frac1d)$. 
It follows from \fullref{thm9.8} that $K$ has unknotting number~1 
if and only if $K$ is an EM--knot. 
Now by \fullref{cor:EMknot}, the EM--knots are all of the form 
$\S(\alpha,\beta;\gamma,\delta)$ with $|\alpha|,|\beta|,|\gamma|,|\delta|<1$
and $\alpha\beta<0$, $\gamma\delta <0$. 
Moreover, it is easy to verify, using \fullref{D8}, that if 
$\S(\alpha,\beta;\gamma,\delta) = \pm\S(\alpha',\beta';\gamma',\delta')$, 
where $|\alpha|,|\alpha'|$, etc. are all $<1$, and $\alpha\beta<0$ 
and $\gamma\delta<0$, then $\alpha'\beta' <0$ and $\gamma'\delta'<0$. 
It follows that a knot $K$ of the form described in the theorem is never 
an EM--knot.
\end{proof}

In Conway's tables, there are 48 knots listed in the form $(a,b)(c,d)$:
3 10--crossing alternating knots, 
7 10--crossing non-alternating knots, 
10 11--crossing alternating knots, and 
28 11--crossing non-alternating knots. 

\begin{thm}\label{thm:9.12} 
Up to 11 crossings, of the 48 knots listed in Conway's tables as 
$(a,b)(c,d)$, the only ones with unknotting number~1 are the four 
non-alternating 11--crossing knots 
$(3,2+)(21,2-)$, $(21,2+)(21,2-)$, $(3,2+)\!\!-\!\!(21,2)$, and  
$(21,2+)\!\!-\!\!(21,2)$.
\end{thm}

\begin{proof}
All the knots in question satisfy the hypotheses of \fullref{thm:9.11} 
except those with $(a,b) = (3,2+)$ or $(21,2+)$. 
Of these, the four listed in the theorem are easily seen to have 
unknotting number~1. 
The others have $(c,d) = (3,2)$, $(21,2)$, $(3,2-)$ or $-(3,2)$. 
{From} now on, let $K$ be a knot of the form $(a,b)(c,d)$ where $(a,b)$ 
(resp.\ $(c,d)$) is one of the two (resp.\ four) possibilities listed. 
Recall (see proof of \fullref{thm:9.11}) that 
$K= \S (-\frac1a,-\frac1b;\frac1c,\frac1d)$.

First note that $K$ is not an EM--knot. 
For $\{\frac1a,\frac1b\} = \{\frac13,\frac32\}$ or $\{\frac23,\frac32\}$, 
and it follows easily from \fullref{D8} that $K$ is  not of the form 
$\S(\alpha,\beta;\gamma,\delta)$ with $|\alpha|,|\beta|,|\gamma|,|\delta|<1$, 
and hence not an EM--knot by \fullref{cor:EMknot}.

By \fullref{thm:Main}, if $K$ has unknotting number~1 then it can be 
unknotted by a crossing move in $\M(-\frac1a,-\frac1b)$ or 
$\M (\frac1c,\frac1d)$. 
Since $|c|,|d|>1$, the former is impossible, by \fullref{cor9.7}. 
So the unknotting move is contained in $\M(\frac1c,\frac1d)$, transforming 
$\M(\frac1c,\frac1d)$ to a tangle $\T$ that unknots 
$\M(-\frac1a,-\frac1b)= \M(-\frac13,-\frac32)$ or 
$\M(-\frac23,-\frac32)$. 
By \fullref{newlem:9.5}, we see that in both cases $\T = \R (-1/2)$.

Consider the case $(c,d) = (3,2)$. 
Then $\M(\gamma,\delta) =\M(\frac13,\frac12)$. 
By \fullref{newcor:9.3},  
this can be transformed to $\R(-1/2)$ by 
a crossing move if and only if there is an integer $k$ such that 
$k\cdot 2 +1 = -1$,
and $1/3 = k+\frac{2s^2}{2rs\pm1}$ for some $s\ge1$, $(r,s)=1$. 
The first equation gives $k=-1$, and the second now gives 
$\pm 2s^2 = 1-3k=4$, which is impossible.

The other three cases $(c,d) = (21,2)$, $(3,2-)$ and $-(3,2)$ are similar. 
We omit the details.
\end{proof}

We now consider the knots of the form $.a.(b,c)$ or $.(b,c).a$.
First we have the following, which is an immediate consequence of 
\fullref{thm9.1}.

\begin{lem}\label{lem9.13}
$.a.(b,c)$ and $-.(b,c).a$ are mutants.
\end{lem}

\fullref{lem9.13} and \fullref{thm:M} 
imply  that $.a.(b,c)$ has unknotting number~1 if and only if 
$.(b,c).a$ does, so we restrict attention to knots of the first type.

\begin{thm}\label{thm9.14}
Suppose that $a>0$ and $|b|,|c|>1$. 
Then $.a.(b,c)$ has unknotting number~1 if and only if 
$\Delta (\frac{a}{a+1},\frac12) = \Delta (-\frac1b,\frac{c+1}c)=1$.
\end{thm}

\begin{proof}
The knot $.x.y$ has the form shown in \fullref{knotxy}.
\begin{figure}[ht!]
\centerline{\includegraphics[height=1.5truein]{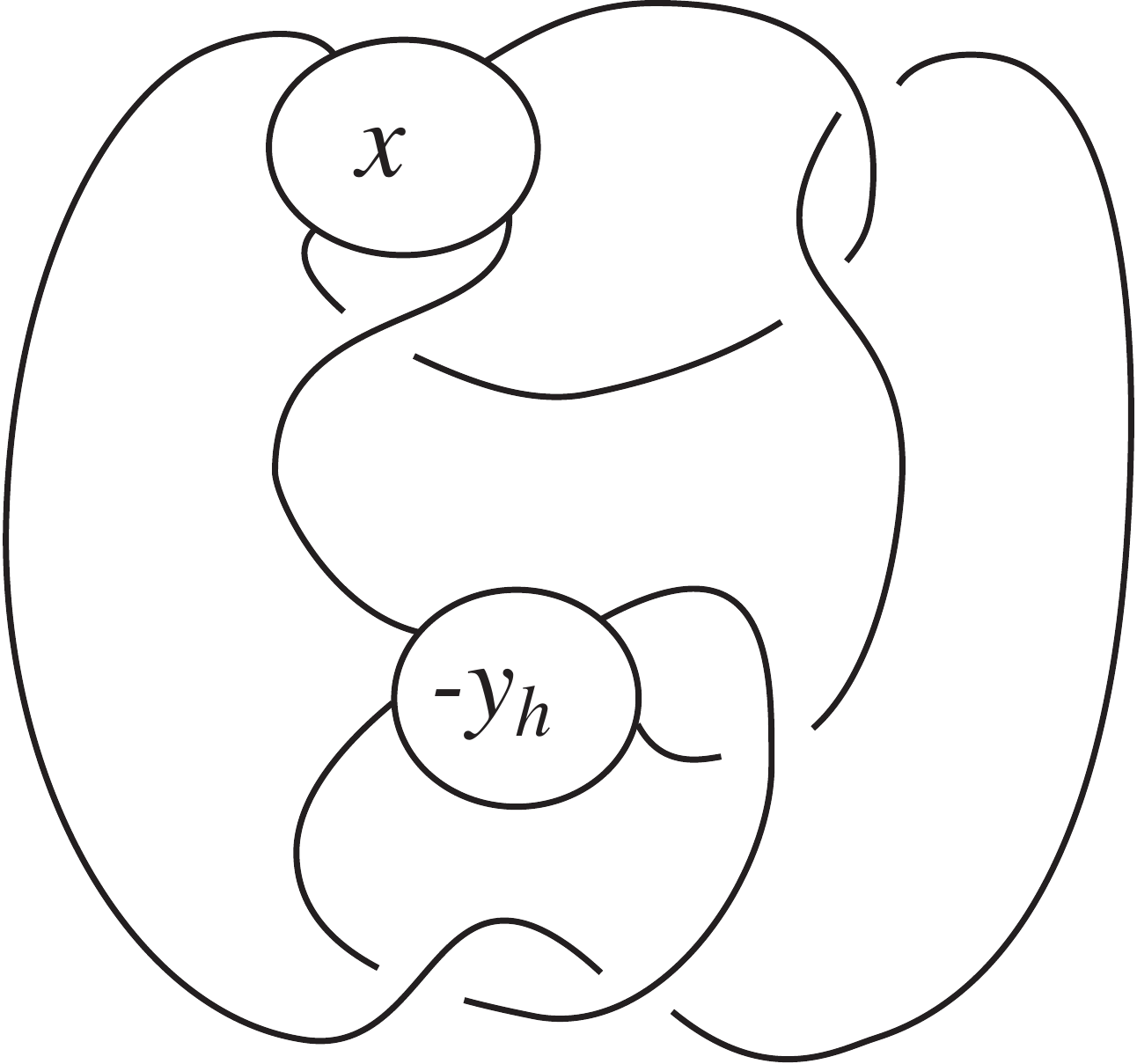}}

\caption{} 
\label{knotxy}
\end{figure}
Therefore $K= .a.(b,c)$ is of the form $\M(\alpha,\beta) \cup 
\M (\gamma,\delta) \cup \M(-1/2,*)$, where ${\alpha = a/(a+1)}$, 
$\beta = -1/2$, $\gamma = -1/b$, and $\delta = -(c+1)/c$; 
see \fullref{knotxyB}.
Note that the conditions on $a,b$ and $c$ guarantee that 
$\alpha,\beta,\gamma,\delta \in \que -\zed$.
This is a decomposition of $K$ into elementary tangles.

\begin{figure}[ht!]
\centerline{\includegraphics[height=1.75truein]{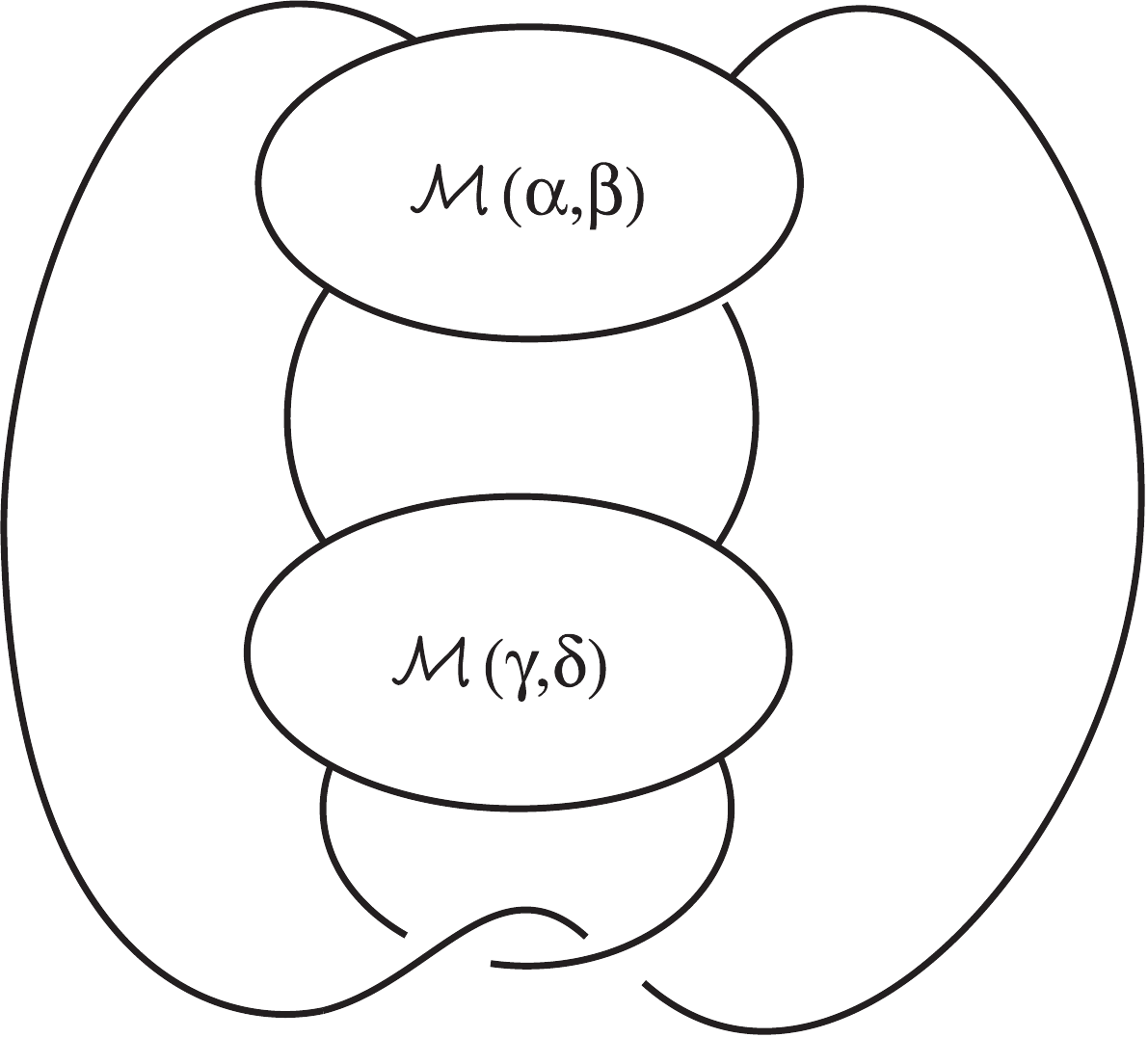}}

\caption{} 
\label{knotxyB}
\end{figure}

Since $K$ does not have a unique essential Conway sphere, $K$ is not 
an EM--knot. 
Hence, if $u(K)=1$, then either conclusion (1) or (3) of 
\fullref{thm8.3} holds. 

First, suppose (1)(a) of \fullref{thm8.3} holds for $\M(\alpha,\beta)$. 
By \fullref{lem9.2}(3),  
$\M(\alpha,\beta)$ is transformed to a 
rational tangle $\R(p'/q)$ with $q\ge 2$. 
If this unknots $K$, then $\R(p'/q)\cup \M(-1/2,*)$ must be an integral 
tangle; see \fullref{unknotsK}. 
Hence $p' =\pm1$, and $\R(p'/q) \cup \M(-1/2,*) = \R(\frac2{1\pm 2q})$.
Therefore $q= 0$ or $\pm1$, a contradiction.
Exactly the same argument shows that (1)(a) of \fullref{thm8.3} does 
not hold for $\M(\gamma,\delta)$.

\begin{figure}[ht!]
\centerline{\includegraphics[height=0.8truein]{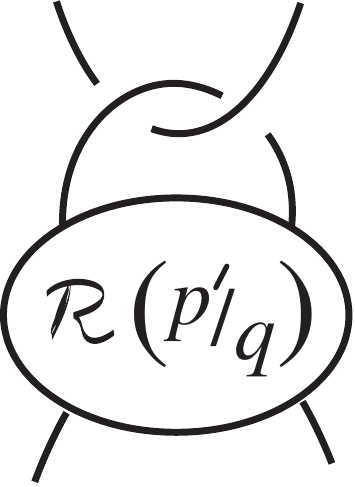}}

\caption{} 
\label{unknotsK}
\end{figure}

Next suppose that we have conclusion (3) of \fullref{thm8.3}
and that (1) of that theorem does not hold. 
By the remark after \fullref{thm8.3},
for $\M = \M(\alpha,\beta)$ or $\M(\gamma,\delta)$, and 
$\M' = \M(\gamma,\delta)$ or $\M(\alpha,\beta)$, respectively, we have 
$\M\cup \M(-1/2,*) =\A$, an EM--tangle,
and the rational tangle $\R$ resulting from the 
standard unknotting move in $\A$ must unknot $\M'$. 
Hence $\R$ must be an integral tangle, $\R(r/1)$.
Here $r= \pm\Delta (\alpha,\beta)$ where  $\alpha$ is the tangle slope 
corresponding to $\R$ and $\beta$ is the tangle slope of the Conway disk 
of $\A$. 
{From} \fullref{Decomp}, we determine that $r= \pm 4m$ when 
$\A = \A_1(\ell,m)$ (the rational tangle is $\R(\frac{1-2m}{4m})$) and 
that $r= \pm (\ell (1-2m)+2)$ (the rational tangle is 
$\R(\frac{2m-1}{\ell(1-2m)+2})$) when $\A = \A_1 (\ell,m)$.

Since $|\alpha|,|\beta| <1$, it follows from \fullref{newlem:9.5}(2)  
that $\M = \M(\alpha,\beta)$, $\M' = \M(\gamma,\delta)$. 
Let $-\frac1b = \frac{r_1}{t_1}$, $-\frac1c = \frac{r_2}{t_2}$ where 
$r_i,t_i \in\zed$ and $t_i >1$. 
Applying \fullref{newlem:9.5}(1) with $x=r$ we get that 
$r  {-\frac1b}  {-\frac1c}  -1 = \frac{\pm1}{t_1t_2}$. 
As $|b|,\, |c|,\, t_1,\, t_2 >1$, $|r| <4$. 
{From} the preceding paragraph, this implies that when $\A = \A_1 (\ell,m)$, 
$|4m| <4$. 
But $m\ne0$, a contradiction.  
When $\A = \A_2 (\ell,m)$, $|\ell(1-2m)+2| <4$. 
Thus $|\ell|\, |1-2m| <6$. 
Since $|\ell| \ge2$ and $m\notin \{0,1\}$, this is again a contradiction. 
Thus conclusion~(3) of \fullref{thm8.3} cannot hold.

We conclude that (1)(b) of \fullref{thm8.3} must hold. 
The crossing move transforms $\M (-\frac12,*)$ to $\M(\frac10,*)$, 
which transforms 
$K$ to the connected sum of $1^* \M(\alpha,\beta)$ and $1^*\M(\gamma,\delta)$.
The result now follows from \fullref{lem9.4}.		
\end{proof}

In Conway's tables, there are 22 knots listed in the form $.a.(b,c)$ or 
$.(b,c).a$: 
8 11--crossing alternating knots and 
14 11--crossing non-alternating knots. 
(They come in mutant pairs $.a.(b,c)$ and $.(b,c).a$.) 
They all have $a= 2$ or $20 = 1/2$, the alternating knots have $(b,c) = 
(3,2)$ or $(21,2)$, while the non-alternating knots have $(b,c)=(3,2-)$, 
$(21,2-)$, $-(3,2)$ or $-(21,2)$.

\begin{thm}\label{thm9.15} 
Of the knots listed in Conway's tables of the form $.a.(b,c)$ or 
$.(b,c).a$, those with unknotting number~1 are precisely the 6 
11--crossing non-alternating knots with $(b,c) = -(3,2)$ or $-(21,2)$.
\end{thm}

\begin{proof} 
This follows easily from \fullref{thm9.14}.
\end{proof}

Finally, we consider the three knots in Conway's tables of the form 
$(a,b) 1 (c,d)$. 
First we have the following lemma. 

\begin{lem}\label{lem:notEMknot}
If $a,b,c,d>1$ then $(a,b)1(c,d)$ is not an EM--knot.
\end{lem}

\begin{proof}
The knot $(a,b)1(c,d)$ has the form shown in \fullref{knotK}. 
Consider the arc that joins the SE--corners of the two Montesinos tangles. 
By swinging this arc over the right-hand tangle, one sees that 
$(a,b)1(c,d)=\S (1/a,1/{b+1};1/c,\allowbreak 1/{d+1})$. 
By Lemma~2.2, it is easy to see that this is never of the form 
$\S(\alpha,\beta;\gamma,\delta)$ with $|\alpha|,|\beta|,|\gamma|,|\delta|<1$.
But the EM--knots are all of this form, by \fullref{lem:EMknot}.
\end{proof}

\begin{thm}\label{thm9.16} 
The 11--crossing alternating knots $(3,2)1(3,2)$, $(3,2)1(21,2)$ and 
$(21,2)1(21,2)$ do not have unknotting number~1.
\end{thm}

\begin{proof}
Let $K= (a,b)1(c,d)$.  		
Note that for the knots under consideration, $\frac1a,\frac1b,\frac1c,
\frac1d \in \que -\zed$. 
Clearly $K$ does not contain an EM--tangle. 
Also, it is not an EM--knot by \fullref{lem:notEMknot}.
\begin{figure}[ht!]
\centerline{\includegraphics[height=1.3truein]{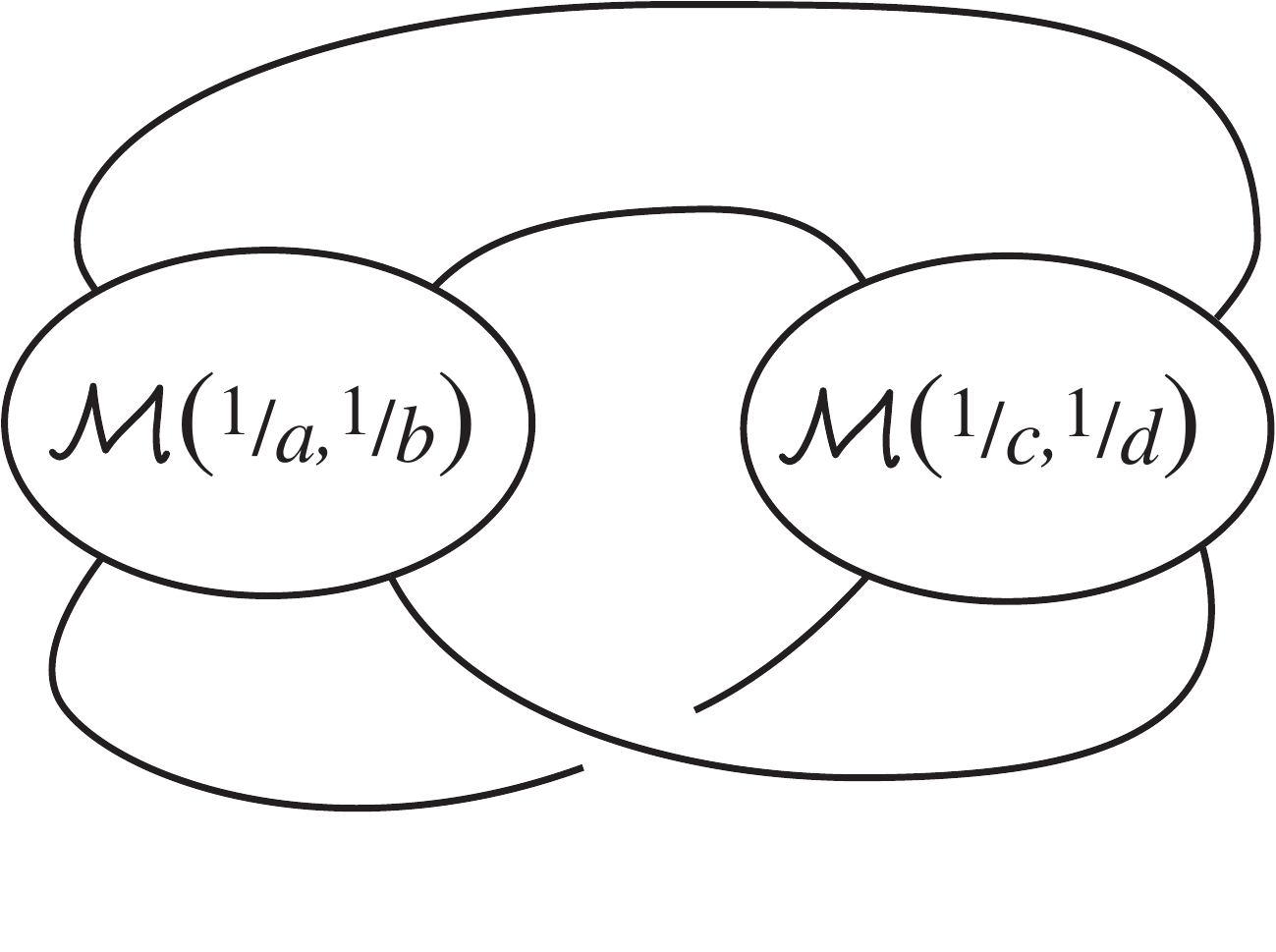}}
\caption{} 
\label{knotK}
\end{figure}

Hence, by \fullref{thm8.3}, if $u(K)=1$ then the unknotting move takes 
place in a rational substituent of $\M(\frac1a,\frac1b)$ or 
$\M(\frac1c,\frac1d)$, transforming it to a rational tangle $\R(r)$. 
By symmetry we may assume that the unknotting move takes place in 
$\M(\frac1a,\frac1b)$. 
By \fullref{newlem:9.5}(1), the tangle $\R(r)$
must be an integral tangle $\R(m)$. 
It follows that $r= [m,-1,0] = \frac{m}{1-m}$. 
By \fullref{newlem:9.5}(2), since $|\frac1c|,|\frac1d|<1$, we must have 
$|m|\le 1$. 
On the other hand, by \fullref{lem9.2}(3),  
there exists an integer 
$k$ such that $\frac{kq+p}q = r = \frac{m}{1-m}$, where $p/q = 1/a$ or $1/b$.
Since $q\ge2$ and $(p,q)=1$, $m=0$ or 1 is impossible. 
Hence $m= -1$ and $q= 2$. 
Therefore $p/q = 1/2 = 1/a$ (since $(a,b) = (3,2)$ or $(21,2)$), 
and $2k+1 = m=-1$, giving $k=-1$.
Also, by \fullref{lem9.2}(4), the other rational substituent of 
$\M(\frac1a,\frac1b)$, namely $1/b$, must satisfy $\frac1b = -1+\frac{2s^2}
{2rs\pm1}$, $s\ge1$, $(r,s)=1$. 
But in both cases, $b=3$ and $b=21=3/2$, this is easily seen to be 
impossible.
\end{proof}

\begin{remark}
For knots of up to 10 crossings, until the work described here and the work
of Ozsv\'ath and Szab\'o \cite{OS}, there were 41 knots for which it was not 
known whether or not they had unknotting number~1.
We ruled out 14 of them (see above): 
the two 9--crossing knots $9_{29} = .2.20.2$ and $9_{32} =.21.20$, and 
the twelve 10--crossing knots $10_{79} = (3,2)(3,2)$, 
$10_{81} = (21,2)(21,2)$, 
$10_{83} = .31.2$, 
$10_{86} = .31.20$, 
$10_{87} = .22.20$, 
$10_{90} = .3.2.2$, 
$10_{93} = .3.20.2$, 
$10_{94} = .30.2.2$, 
$10_{96} = .2.21.2$, 
$10_{148} = (3,2)(3,2-)$, 
$10_{151} = (21,2)(21,2-)$, and 
$10_{153} = (3,2)\!\! -\!\! (21,2)$. 
Meanwhile, Ozsv\'ath and Szab\'o \cite{OS}, using their remarkable Heegaard 
Floer homology theory, ruled out all 41 knots except $10_{153}$. 
So the knots with 10 or fewer crossings and unknotting number~1 are 
completely determined. 
\end{remark}

\section{Deciding if a large algebraic knot has unknotting number~1}

It is unknown if the unknotting number of a knot is a computable invariant. 
Even the following special case is open:

\begin{quest}\label{quest11.1}
Is there an algorithm to decide whether or not a given knot has unknotting 
number~1?
\end{quest}

Note that by Haken \cite{Hak} there is an algorithm to decide if a knot 
has unknotting number~0.

\fullref{thm8.3} allows us to answer \fullref{quest11.1} affirmatively
for large algebraic knots.

\begin{thm}\label{thm11.1}
There is an algorithm  to decide whether or not a given large algebraic 
knot $K$, described as a union of elementary marked tangles 
(\fullref{alg-tangle}) and 4--braids in $S^2\times [0,1]$, has 
unknotting number~1, and, if so, to identify an unknotting crossing move.
\end{thm}

\proof
Note that none of the rational tangles in a constituent elementary 
tangle of $K$ is integral or $\R(1/0)$ (that is, the distance between the 
slope of such a rational tangle and the orbifold $S^1$--fiber of the elementary 
tangle is at least two).

\begin{Defin}
Let $\T$ be a marked tangle in the 3--ball. 
The number $p/q \in \que \cup \{\infty\}$ is called an {\em unknotting slope}
for $\T$ if $\T\cup \R(-p/q)$ (ie $1^*  (\T+\R (-p/q))$) 
is the unknot.
\end{Defin}

\begin{remark} 
As long as $\T$ is not rational, an unknotting slope for $\T$, if there is one,
is unique. 
This follows, for example, from the fact that knots are determined by their 
complements and the $\zed_2$--Smith Conjecture (applied to the double 
branched covers). 
\end{remark}

\begin{lem}\label{lem-alg}
There is an algorithm to decide whether or not a given 
algebraic tangle (ie a union of elementary tangles and 4--braids) 
has an unknotting slope, and, if so, to find it.
\end{lem}

\begin{proof} 
Again we assume the tangle is given as a union of elementary marked 
tangles and 4--braids in $S^2\times [0,1]$. 
Begin with the innermost, constituent, elementary tangles, necessarily of 
type~I, and compute their unknotting slopes (if they exist)
via \fullref{newlem:9.5}(1).
Then work outward along elementary tangles of type~II determining unknotting
slopes at each step.
This is equivalent to solving equations of the following form:
given $\frac{r_1}{s_1},\frac{r_2}{s_2}$ find a rational $\frac{x}y$ 
such that $\frac{x}y + \frac{r_1}{s_1} = \frac{r_2}{s_2}$.
In working outward, if one comes to a constituent tangle of type~III, 
the corresponding unknotting slope (if it exists) must be the slope of 
the orbifold $S^1$--fiber.

Working outward in this  way, either one finds at some point that 
there is no unknotting
slope, in which case there is none for the algebraic tangle, $\T$, or one
determines the unknotting slope for $\T$.
\end{proof}

By \fullref{thm8.3}, $K$ has unknotting number 1 if and only if 
one of the following options holds: 
(A)~$K$ is an EM--knot; 
(B)~$K$ contains an EM--tangle for which the standard crossing change 
unknots $K$; or 
(C)~$K$ unknots by replacing a rational tangle in a constituent elementary 
tangle of type~I or II with another rational tangle as described in 
\fullref{thm8.3}(1).

We show that these options can be checked algorithmically.

\noindent {\bf (A)}\quad
To see if $K$ is an EM--knot, first check that there are two elementary 
tangles of type~I whose union is $K$. 
If so, compute the orders of the exceptional orbifold $S^1$--fibers of 
each elementary tangle and list the finite number of EM--knots having 
exceptional fibers of the same order.
Check if $K$ is equivalent to one of these (eg, check that the orbifold 
$S^1$--fibers of the two elementary tangles intersect twice at the unique 
Conway sphere [ie, that the distance between the slopes of these fibers 
is 1 on the Conway sphere]. 
If so then $K$ can be written in the form $\S(\alpha,\beta; \gamma,\delta)$
and \fullref{D8} may be applied).

\noindent {\bf (B)}\quad
To check the unknotting of $K$ by a standard crossing change in 
an EM--tangle, we list all pairs 
$\{X_1,X_2\}$ of a type~I and type~II elementary tangle of $K$ that 
share a common Conway sphere. 
The union $X_1\cup X_2$ is a candidate for an EM--tangle.
Let $-\frac{p}q$ be the unknotting slope of the complementary tangle 
$(S^3,K) - (X_1\cup X_2)$. 
List the finitely many EM--tangles that have exceptional orbifold 
$S^1$--fibers of the same order as $X_1\cup X_2$.
Then check if there is a homeomorphism from the candidate $X_1\cup X_2$
to one of these EM--tangles taking the slope $\frac{p}q$ to the 
slope of the rational tangle that results from the standard crossing change.
(For example, check that the orbifold $S^1$--fibers of $X_1$ and $X_2$ intersect
twice along the common Conway sphere. 
If so, $X_1\cup X_2$ may be rewritten in the form $\S(*,\frac{p_1}{q_1};
\frac{p_2}{q_2},\frac{p_3}{q_3})$.
Then check that 
$$\frac{p_1}{q_1} +n_1= \frac{v_1}{w_1}\ ,\quad 
\frac{p_2}{q_2} +n_2 = \frac{v_2}{w_2}\ ,\quad
\frac{p_3}{q_3} -n_2 = \frac{v_3}{w_3}$$
for some $n_1,n_2\in\zed$, where $\A_\ep (\ell,m) = \S (*,\frac{v_1}{w_1};
\frac{v_2}{w_2}, \frac{v_3}{w_3})$.
If so, there is a unique homeomorphism up to isotopy taking $X_1\cup X_2$ 
to $\A_\ep (\ell,m)$. 
One then checks that $\frac{p}q$ is identified with the slope of the 
rational tangle gotten by the standard crossing move on $\A_\ep (\ell,m)$.)

\noindent {\bf (C)}\quad
To check the condition of \fullref{thm8.3}(1), we check each 
constituent elementary tangle of type~I or II as follows. 
For any elementary tangle of type~II, we replace the rational tangle 
$\R (p/q)$ with $\R(1/0)$. 
If this yields the unknot (which can be checked algorithmically), then we 
check that $\frac{p}q = \frac{2rs\pm1}{2s^2}$ for some $(r,s) =1$. 
If so one identifies the unknotting crossing move of $K$ as described, 
for example, in \fullref{lemK} (or \fullref{lem:crossingchange}).

Consider a constituent elementary tangle, $X$, of type~I, and let 
$-\frac{v}{w}$ be the unknotting slope for the complementary tangle 
$(S^3-K)-X$ (if there is none, $K$ cannot be unknotted in this way).
Decide if replacing some rational tangle, $\R(p/q)$, of $X$ by an integer 
tangle, $\R(k)$, changes $X$ to the rational tangle $\R(v/w)$.
If so, determine if $\frac{p}q = \frac{2s^2}{2rs\pm 1} +k$ for some $(r,s)=1$.
If so, then a variation of \fullref{lemK} will determine an unknotting 
crossing move for $K$ (see also \fullref{lemH}).
\qed
\vspace{-4pt}

\section{Unknotting in a minimal diagram}
\vspace{-4pt}

In \cite{Koh}, Kohn made the following conjecture, which he showed was 
true for 2--bridge knots and links. 
\vspace{-4pt}

\begin{conjecture}[Kohn]\label{kohn}
Let $K$ be a knot or link with $u(K)=1$. 
There is a crossing in a minimal diagram of $K$ which, when changed, 
unknots $K$.
\end{conjecture}
\vspace{-4pt}

Note that the analog for knots with $u(K)>1$ is false \cite{B}, \cite{N}.
\vspace{-4pt}

We shall show that the conjecture is true for alternating large 
algebraic knots.
\vspace{-4pt}

{From} \cite{SuT}, we take the following 
\vspace{-4pt}

\begin{Defin}
A tangle diagram $D$ in a disk $\Delta$ is {\em prime} iff
\begin{itemize}
\item[{(i)}] the underlying projection of $D$ is a connected subset 
of the disk $\Delta$
\item[{(ii)}] if $C$ is a circle in $\Delta$ meeting $D$ transversely 
in two points, then these points belong to the same edge of $D$
(ie, diagrammatic connected sums are not allowed).
\end{itemize}
\end{Defin}

\begin{prop}[Flyping conjecture for alternating tangles in a 3--ball]
\label{prop:flyping}
If $D,D'$ are prime, alternating diagrams of the same marked tangle 
in a 3--ball, then $D$ and $D'$ differ by a sequence of flypes.
\end{prop}

\begin{proof}
As outlined in \cite[page 333]{SuT} or \cite[page 998]{T2}, 
this follows from the rigid vertex version of the Flyping Conjecture
proved in \cite{MT}.
\end{proof}

\begin{Defin}
A tangle $(B,T)$ in a 3--ball $B$ is {\em prime} iff 
\begin{itemize}
\item[{(i)}] $(B,T)$ is not rational;
\item[{(ii)}] if $S$ is a 2--sphere in $B-T$, then the 3--ball bounded by 
$S$ does not meet $T$;
\item[{(iii)}] if a 2--sphere $S$ in $B$ meets $T$ transversely in two points,
then the 3--ball in $B$ bounded by $S$ meets $T$ in an unknotted arc.
\end{itemize}
\end{Defin}

\begin{cor}\label{cor11.2}
If $D,D'$ are alternating diagrams of the same marked tangle in a 3--ball
which is either prime or rational but not $\R(1/0)$, $\R(0/1)$, then any 
(marked) tangle obtained by changing a crossing in $D$ can be gotten by 
changing a crossing in $D'$.
\end{cor}

\begin{proof}[Proof of \fullref{cor11.2}] 
The hypotheses guarantee that if $D$ is not a prime diagram, then it has 
a nugatory crossing. 
(A circle violating (ii) of primeness of the diagram 
must encircle an alternating diagram of the unknot. 
This must have a nugatory crossing by the minimality of crossing number for 
reduced, alternating diagrams of links.)
Thus by reducing all nugatory crossings in $D$, we leave a prime, 
alternating diagram. 
Similarly for $D'$. 
Thus $D,D'$ are related by a sequence of flypes and nugatory crossing 
reductions or creations. 
One checks that changing a crossing commutes with these operations.
\end{proof}

Recall (\fullref{lemK}) that $\R(p/q)$ can be transformed to $\R(1/0)$ 
by a crossing move if and only if there are coprime integers $r,s$ 
such that $p/q = \frac{2rs\pm1}{2s^2}$.

\begin{lem}\label{lem:crossingchange}
If $p/q = \frac{2rs\pm1}{2s^2}$, where $(r,s) =1$ and $s\ne0$, then $\R(p/q)$ 
can be transformed to $\R(1/0)$ by a crossing change in any alternating 
diagram for $\R(p/q)$.
\end{lem}
\vspace{2pt}

\begin{proof}
By \cite{Koh} the condition on $p/q$ is equivalent to the 
condition that $\pm p/q =$
\begin{equation*}
\begin{split}
&[c_1,\ldots,c_{\ell-1},c_\ell,1,1,c_\ell-1,c_{\ell-1},\ldots,c_1,c_0]
\ \text{ or}\\
 &[c_1,\ldots,c_{\ell-1},c_\ell-1,1,1,c_\ell,c_{\ell-1},\ldots,c_1,c_0]
 \end{split}
 \end{equation*}
where $c_i\ge1$, $0\le i\le \ell-1$, $c_\ell\ge 2$. 
The corresponding diagram of $\R(p/q)$ is alternating, and the crossing 
change is visible  in that diagram.
By \fullref{cor11.2}, 
the transformation $\R(p/q)\mapsto \R(1/0)$ can be effected 
by a crossing change in any alternating diagram. 
\end{proof}
\vspace{2pt}

\begin{lem}\label{lemH}
If $p/q\in \que-\zed$, and  
$\R(p/q)$ can be transformed to $\R(k)$ by a crossing move, then it can 
be transformed to $\R(k)$ by a crossing change in any alternating diagram of 
$\R(p/q)$, unless 
$p/q = \pm [\ell,2,m]$, where $\ell>0$, $m\ge0$, and  $k=\pm  (\ell+m+2)$.
\end{lem}
\vspace{2pt}

\proof 
By \fullref{cor11.2}, we only need exhibit a crossing change in 
some alternating diagram of $\R(p/q)$.
If $k=0$ then the result holds by rotating and applying Lemms~\ref{lemK}
and \ref{lem:crossingchange}. 
So assume $k\ne0$. 
We may also suppose, without loss of generality, that $p/q>0$. 
Write $p/q = p'/q+m$, $0<p'/q <1$, $m\ge 0$. 
Then $\R(p/q-k) = \R(p'/q +(m-k))$ can be transformed to $\R(0)$ by a 
crossing move.

\begin{figure}
\centerline{\includegraphics[height=0.75truein]{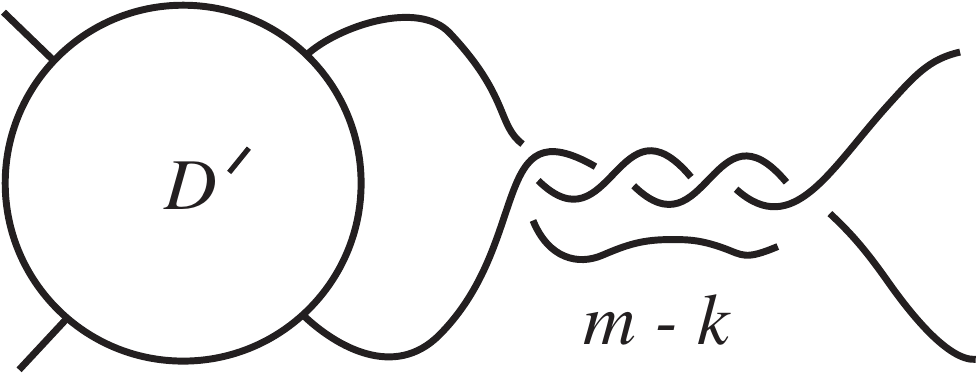}}

\caption{}		
\label{case1}
\end{figure}

Let $D'$ be a positive alternating diagram of $\R(p'/q)$
(see \fullref{case1}).  
There is such a diagram since $\frac{p'}q >0$.
\medskip

\noindent {\bf Case (1)\qua \boldmath$m-k\ge0$}

In this case $\R(p/q-k)$ has the alternating diagram  shown in 
\fullref{case1}

\noindent By the case $k=0$ above, 
$\R(p/q-k)$ can be transformed to $\R(0)$ by changing 
a crossing $c$ in this diagram. 
Clearly $c$ must be a crossing of the diagram $D'$ 
(eg, by computing the associated rational number after the crossing change). 
Let $D$ be the alternating diagram of $\R(p/q) = \R(p'/q +m)$ obtained by 
putting $m$ horizontal $\frac12$--twists on the right of $D'$. 
Then changing the crossing $c$ in $D$ transforms $\R(p/q)$ to $\R(k)$.
\medskip

\noindent {\bf Case (2)\qua \boldmath$m-k<0$}

\begin{figure}[ht!]
\centerline{\includegraphics[height=1.5truein]{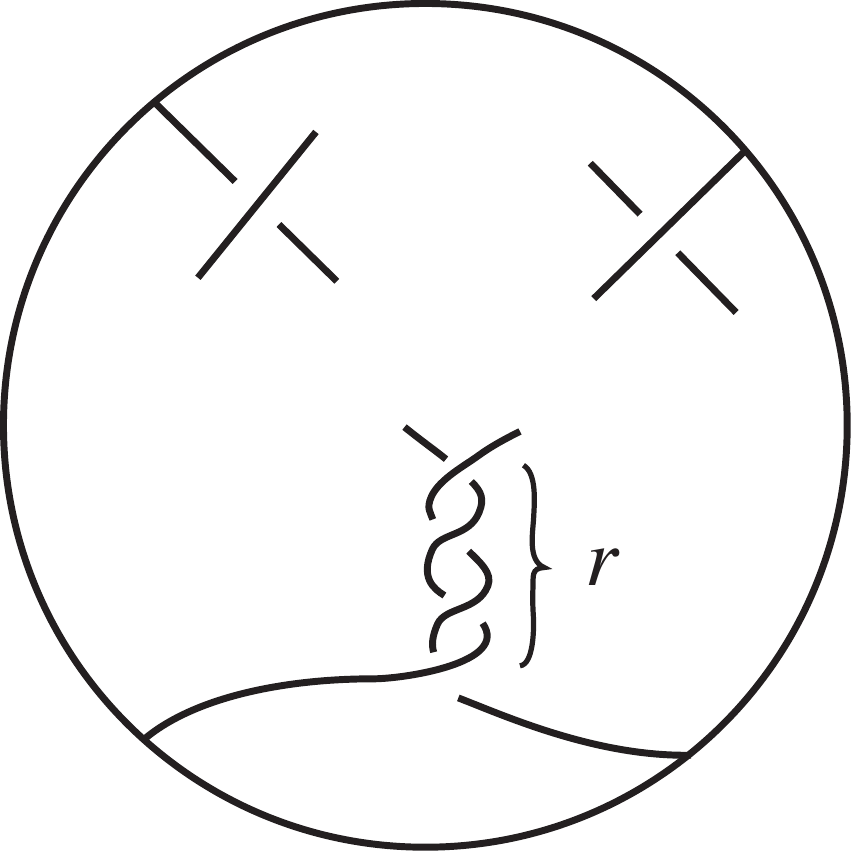}}

\caption{}		
\label{case2a}
\end{figure}

\begin{figure}[ht!]
\centerline{\includegraphics[height=1.25truein]{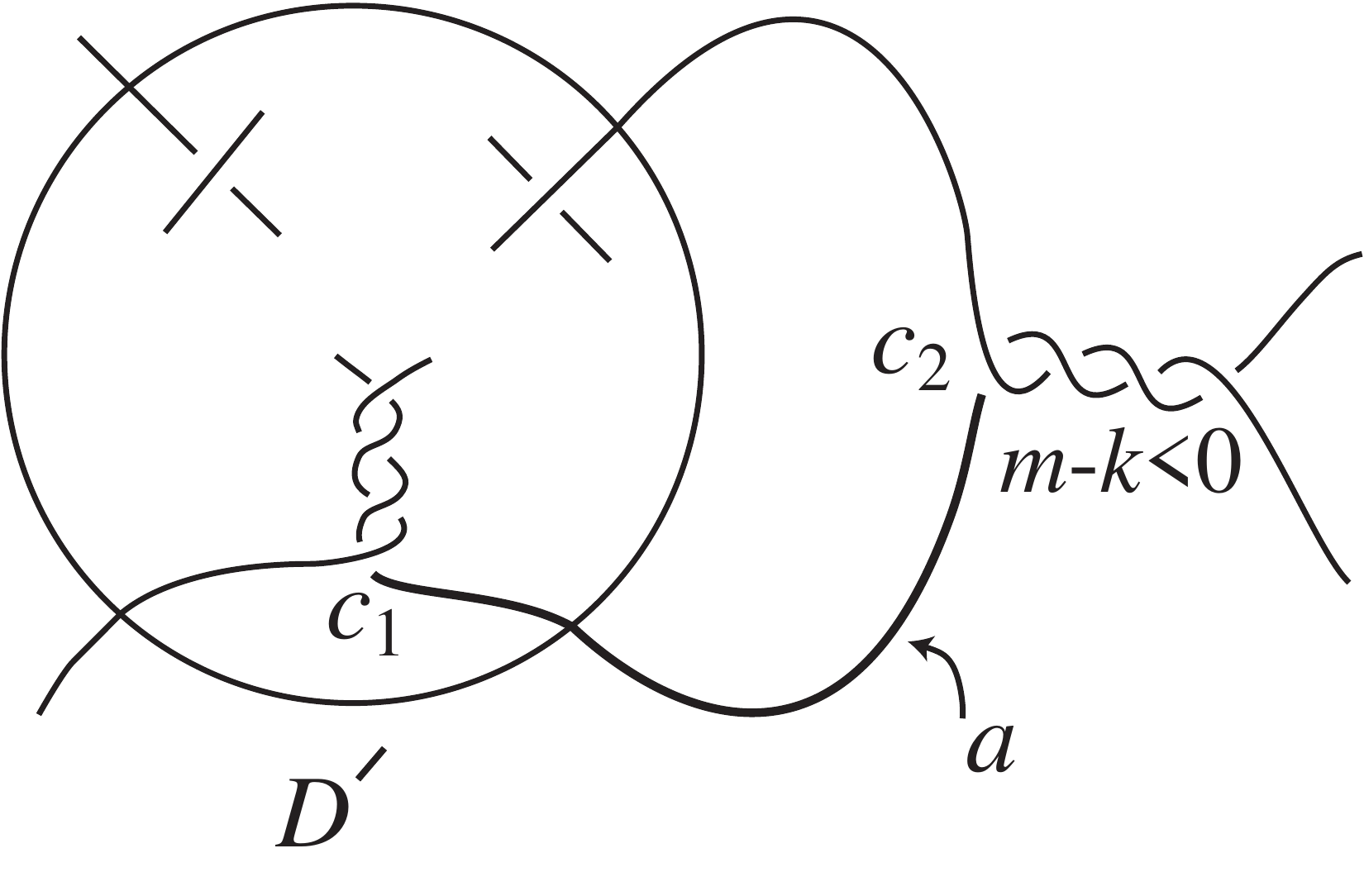}}

\caption{}		
\label{case2b}
\end{figure}

Since $0<p'/q <1$, the diagram $D'$ ends up with $r\ge1$ vertical 
$\frac12$--twists; see \fullref{case2a}.
Hence $\R(p/q -k) = \R(p'/q + (m-k))$ has the diagram $D_0$ shown in 
\fullref{case2b}.


\begin{figure}
\centerline{\includegraphics[height=1.2truein]{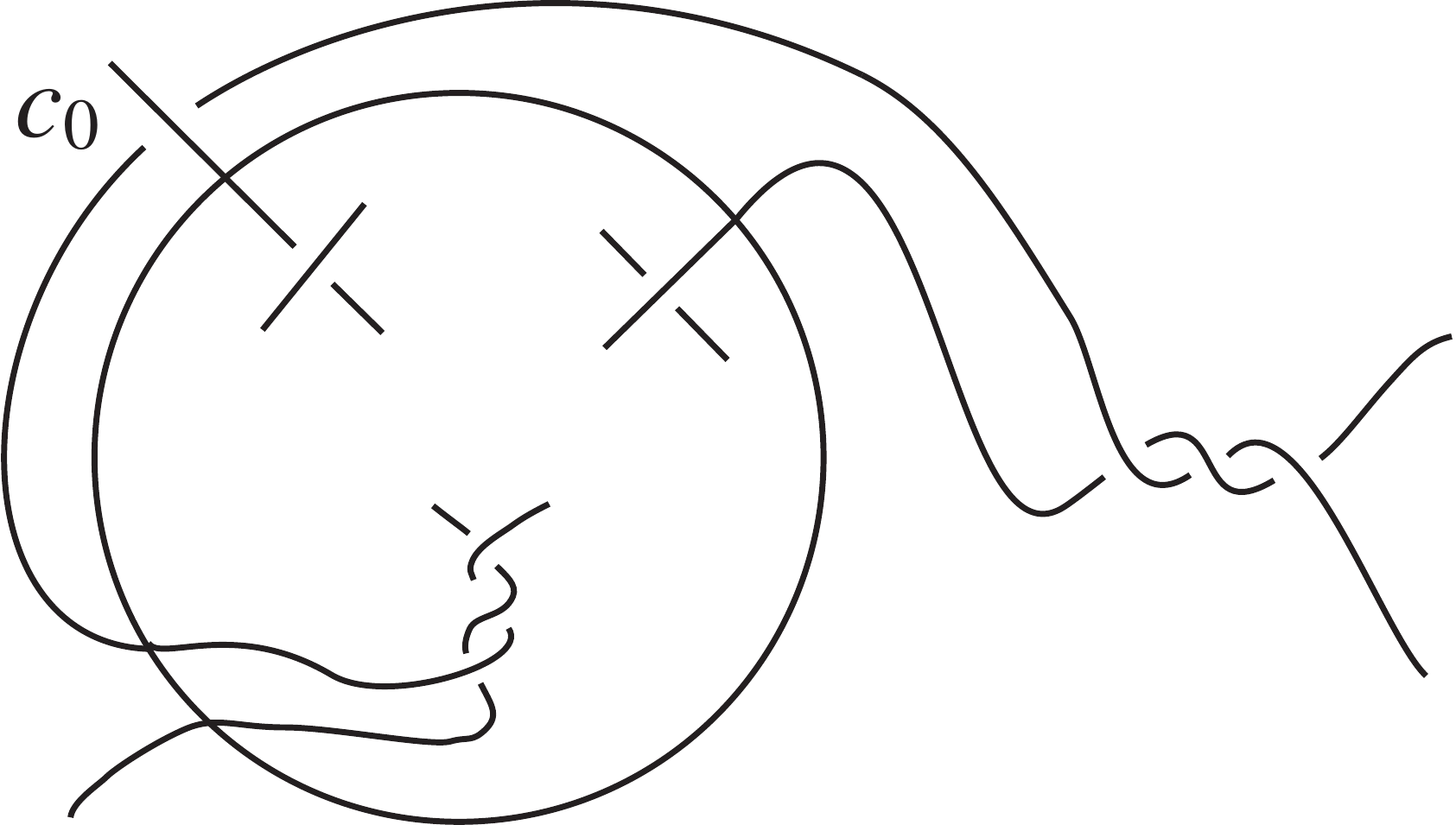}}

\caption{}		
\label{case2c}
\end{figure}

Let $a$ be the arc indicated by the bold line in \fullref{case2b}. 
Swinging $a$ underneath $D'$ gives the diagram $D_1$ shown in 
\fullref{case2c}. 
Note that this is a (negative) alternating diagram. 
Therefore, by the case $k=0$ above,
$\R(p/q-k)$ can be transformed to $\R(0)$ 
by changing some crossing $c$ in $D_1$. 
Clearly $c$ is either a crossing of $D'$ (other than $c_1$) or the 
new crossing $c_0$.

In the first case, $\R(p'/q +m) = \R(p/q)$ is transformed to $\R(k)$ by 
changing the same crossing $c$  in the alternating diagram $D$ of 
$\R(p'/q+m)$ defined in Case~(1). 

In the second case, changing the crossing $c_0$ in $D_1$ and swinging the 
arc $a$ over $D'$ clearly gives the 
same tangle as changing the crossings $c_1$ and $c_2$ in the diagram $D_0$.
By hypothesis, this is $\R(0)$. 
Hence changing only $c_1$ in $D_0$ gives $\R(-2)$. 
Therefore $r\le 2$ (else the crossing change would not yield an 
integral tangle). 

\begin{Claim}
If $r=1$, then $\frac{p}q = [\ell-1,1,1,m]$ where $\ell\ge 2$, $k=m+2-\ell$.

If $r=2$, then $\frac{p}q = [\ell,2,m]$ where $\ell>0$, $k=m+2+\ell$.
\end{Claim}

\begin{proof}[Proof of Claim]
Assume $r=1$. 
The positive continued fraction expansion for the rational tangle  of 
\fullref{case2b} with crossing $c_1$ changed gives 
$$-2 = (m-k) + \frac1{-1+\frac{s}{\ell}} 
= (m-k) + \frac{\ell}{s-\ell}$$
where $0<\frac{s}{\ell} <1$, $(s,\ell) =1$. 
Integrality gives $s=\ell-1$ and 
$$\frac{p}q = m+\frac1{1+\frac{s}{\ell}} 
= m+ \frac1{1+\frac1{1+\frac1{\ell-1}}}$$
where $\ell-1>0$. 
Thus $\frac{p}q = [\ell-1,1,1,m]$. 
Furthermore 
$$-2 = (m-k) + \frac1{-1+\frac1{1+\frac1{\ell-1}}} 
= m-k-\ell$$
as required.

Similarly, if $r=2$ we have 
$$-2 = (m-k) + \frac1{0+\frac{s}{\ell}}$$
This implies that $s=1$, $\ell>0$. 
Then $-2= (m-k) +\ell$. 
Finally $\frac{p}q = [\ell,2,m]$.\qed

By the Claim, if $r=1$ then we are in Case~(1). 
When $r=2$, we obtain the list of tangles stated in the Lemma.
\end{proof}

\begin{thm}\label{thm8.4}		
Let $K$ be an alternating large algebraic knot with unknotting number~$1$. 
Then $K$ can be unknotted by a crossing change in any alternating diagram 
of $K$.
\end{thm}

\proof
By \cite{MT}, any two reduced alternating diagrams of $K$ are related 
by flype moves. 
It follows easily that if $K$ can be unknotted by a crossing change in 
some alternating diagram then it can be unknotted by a crossing change 
in any alternating diagram. 
\vspace{-4pt}

In what follows we use the notion of the ``visibility'' of a Conway 
sphere or disk in an alternating diagram as discussed in \cite{Th1}. 
In particular, \cite{Me} shows that in an alternating diagram a Conway 
sphere is either visible, or {\em hidden} in a very specific way 
(see Figures~3(i), (ii) of \cite{Th1}). 
In the latter case, there  is a standard move on the diagram to make 
the sphere visible (see Figure~3(iii) of \cite{Th1}). 
For a reduced alternating diagram of an elementary tangle, \cite[page 326]{Th1} 
shows that the arguments of \cite{Me} can also be used to say that the 
Conway disk must be visible.
\vspace{-4pt}

First suppose that we are in Case (1) of \fullref{thm8.3}.  
\vspace{-4pt}

Let $D$ be a reduced alternating diagram of $K$. 
Suppose we are in subcase~(a),  so that the unknotting crossing move takes
place in a rational subtangle of an elementary tangle $\T$ of type~I.
Let $S$ be the boundary of $\T$, and suppose that $S$ is visible in $D$.
Then after flyping if necessary (see \cite[page 326]{Th1} for the visibility of 
the Conway disk), 
we may assume that $D$ contains a 
subdiagram of the form shown in \fullref{alg-tangle}(I). 
By \fullref{thm8.3}, the crossing move transforms $\R(p/q)$ to $\R(k)$. 
Since $p/q\notin \zed$, it follows from \fullref{lemH} that this can be 
achieved by a crossing change in the diagram $D_1$, unless (without loss 
of generality) $p/q = [\ell,2,m]$, $\ell>0$, $m\ge0$, and $k=\ell+m+2$. 
Since $k>0$ and $D$ is alternating, we see that replacing $D_1$ with the 
standard diagram of $\R(k)$ gives an alternating diagram $D'$. 
Also, since $D$ is reduced and $S$ is essential, it is easy to see that $D'$ 
is reduced. 
Hence $D'$ is a diagram of a non-trivial knot, a contradiction.
\vspace{-4pt}

Next suppose that we are in subcase (b) of \fullref{thm8.3}, Case~(1), and 
that the boundary components of the corresponding elementary tangle of 
type~II are both visible in $D$ (\fullref{alg-tangle}(II)).
The crossing move transforms $\R(p/q)$ to $\R(1/0)$, and,  
by \fullref{prop:flyping} and (the proof of) 
\fullref{lemK}, 
this can be achieved by a crossing change in the diagram $D_3$ ($q\ne0$).
\vspace{-4pt}

\begin{figure}[ht!]
\centerline{\includegraphics[height=2.2truein]{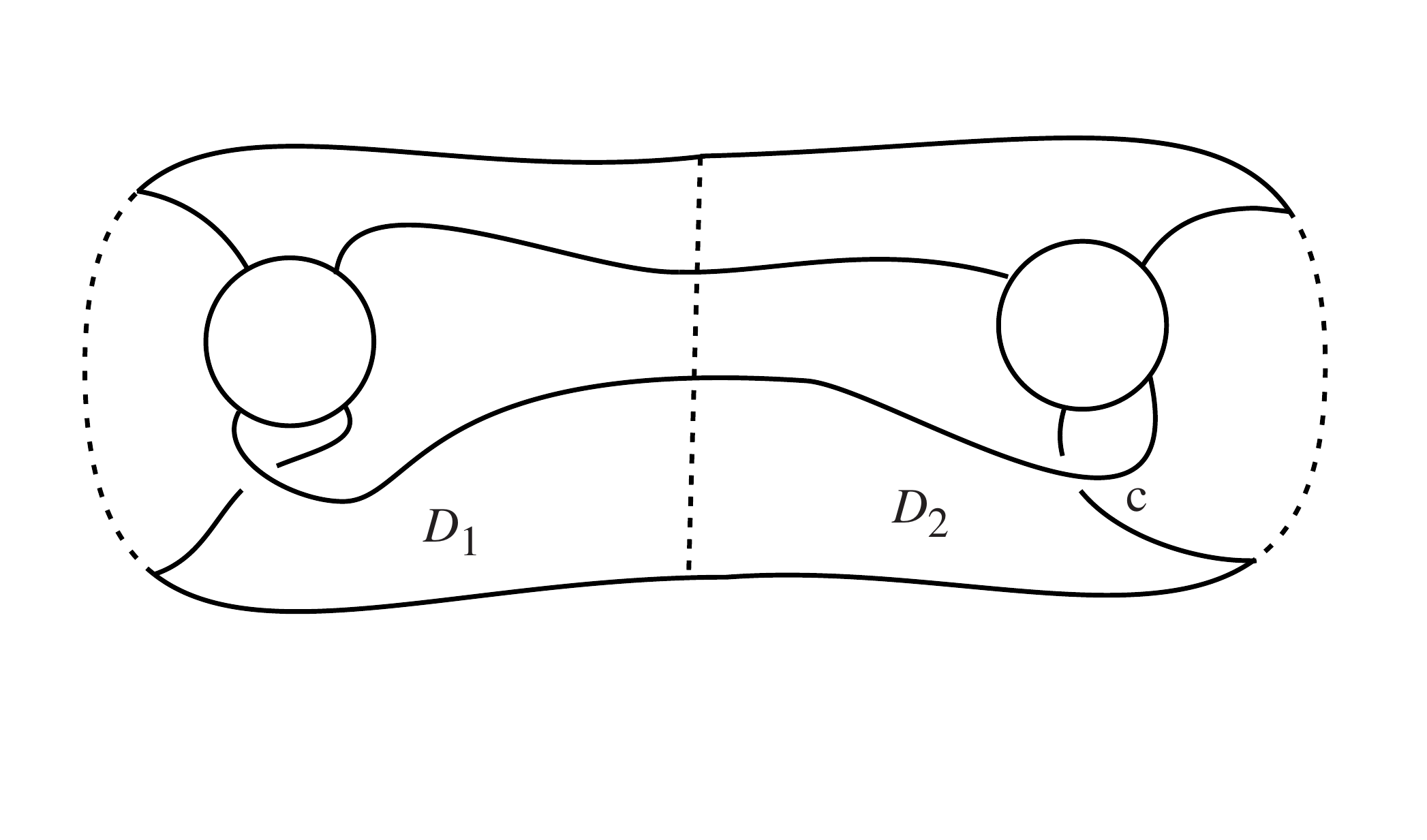}}
\vspace{-15mm}
\caption{$\T = \R(r/s) + \R(p/q)$}		
\label{hid}
\end{figure}

It remains to consider (a) and (b) when the relevant Conway spheres are 
hidden in $D$.
So suppose we are in subcase~(a), and the boundary $S$ of the corresponding 
elementary tangle $\T$ of type~I is hidden in $D$.
Making $S$ visible as described in \cite{Th1},
we get a diagram in 
which the tangle $\T$ appears as in \fullref{hid}. 
Note that the diagrams $D_1,D_2$ of $\R(r/s)$, $\R(p/q)$ in 
\fullref{hid} will be alternating.
Suppose, without loss of generality, that the unknotting crossing move takes 
place in the right-hand rational tangle $\R(p/q)$, transforming it to 
$\R(k)$ for some integer $k$.
\vspace{-4pt}

We first argue that $k=0$ or 1. 
Assume not.
Let $\R(c/d)$ be the subtangle of $\R(r/s)$ 
encapsulated in the circle of $D_1$ of \fullref{hid}. 
We then have the equation $\frac{d}c = 1+\frac{s}r$. 
Note that since the diagram of $\R(r/s)$ is alternating, $c/d <0$. 
Now the crossing move we are considering turns $\T$ into $\R(\frac{r}s +k)$. 
This has to be of the form $\R(1/x)$, $x\in\zed$, or $\R(0/1)$. 
(Write $K = \T\cup \T'$, as in Figure~3(iii) of \cite{Th1}, where $\T'$ 
is also of the form of \fullref{hid}. 
Since the corresponding subdiagrams, $D'_i$, of $\T'$ are alternating, 
if either $D'_i$ is a diagram of a rational tangle $\R(p/q)$, then $q>0$. 
Arguing on the level of double branched covers, as in 
\fullref{newlem:9.5}(1), we see that $\Delta (\frac{r}s +k,\frac01)\le1$. 
Hence $\frac{r}s +k = \frac1x$ or $\frac01$.)
Since $D_1$ in \fullref{hid} for $\T$ is alternating, $|\frac{r}s | <1$. 
Thus $\frac{r}s+k\ne \frac01$. 
We assume $\frac{r}s +k=\frac1x$, $x\in \zed$. 
Then $\frac{r}s = \frac{1-kx}{x}$, hence $\frac{s}r = \frac{x}{1-kx}$. 
Therefore $0>\frac{d}c = 1+\frac{x}{1-kx} = (1-(k-1)x)/(1-kx)$, implying 
that $k=0,1$.
\vspace{-4pt}

By \fullref{lemH}, since $\frac{p}q \notin\zed$ 
and $k=0,1$, we see that $\R(p/q)$ can be 
transformed to $\R(k)$ by a crossing change in the minimal diagram $D_2$.
Now $D_2$ can be obtained by adding a vertical right-handed twist, given 
by the crossing marked $c$, to the diagram of the tangle $\R(a/b)$ which 
was visible in the original alternating diagram of $K$. 
If the crossing changed in $D_2$ is not $c$, then this is a crossing change 
in the original diagram. 
So suppose the crossing changed is $c$.
Since this gives $\R(k)$, we see that $\R(p/q)$ ($\R(a/b)$, resp.) is gotten
by adding two (one, resp.) right-handed vertical twists to $\R(k)$. 
Since $\frac{p}q \ne \frac01$, we see that $k\ne0$. 
Thus $k=1$ and $\frac{a}b = \frac12$, $\frac{p}q = \frac13$. 
But then one sees that the crossing change at $c$ can be accomplished by a 
crossing change in the diagram of  $\R(a/b)$ by \fullref{cor11.2}. 
Thus $K$ can be unknotted by a crossing change in $D$.

\begin{figure}[ht!]
\centerline{\includegraphics[height=1.35truein]{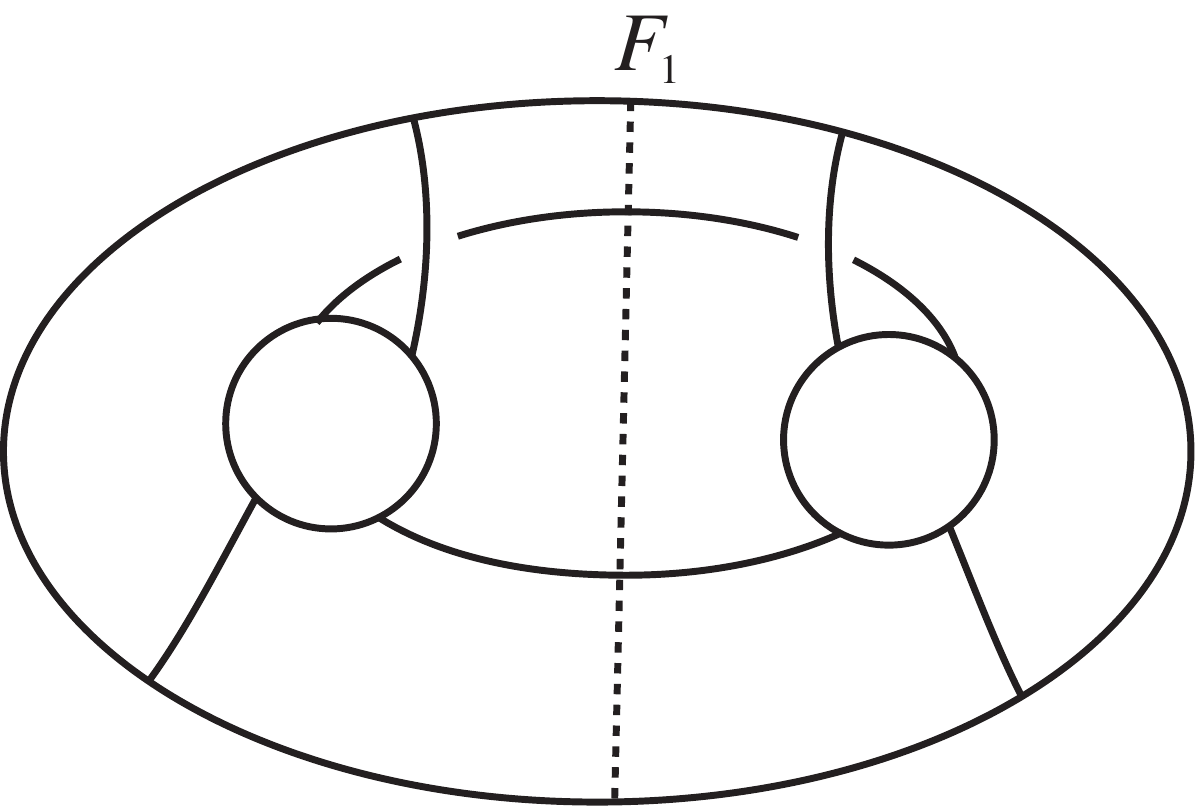}}
\caption{} 
\label{algknots2.8}
\end{figure}

To finish the proof of \fullref{thm8.4} in Case~I, we consider the case
when the crossing move is in an elementary tangle of type~II where one 
of the boundary components, $S_1$, is hidden in $D$. 
Making $S_1$ visible gives a diagram containing a subdiagram as shown in 
\fullref{algknots2.8}.

\begin{lem}\label{lem:visibleII}
Let $F_1$ be the disk pictured in \fullref{algknots2.8}. 
Then $F_1$ is a Conway disk for $S_1$. 
Furthermore, any Conway disk for $S_1$ is parallel to $F_1$.
\end{lem}

\begin{proof}
By Corollary 3.3 of \cite{Th1}, $F_1$ is a Conway disk for the prime tangle 
bounded by $S_1$ in \fullref{algknots2.8}. 
We assume for contradiction that there is a Conway disk, $F_2$, of $S_1$ 
that is not parallel to $F_1$. 
Then we may take $F_2$ to be disjoint from $F_1$. 
In particular, the slope of $F_2$ on $S_1$ is $\frac10$ (in the diagram 
coordinates). 
Let $\T$ be the tangle containing $F_2$ after cutting \fullref{algknots2.8}
along $F_1$. 
Without loss of generality assume this is the right-hand side of $F_1$. 
Then $\T$ is an alternating tangle for which $F_2$ is an essential Conway 
disk with slope $\frac10$ ($F_2$ is not parallel to $F_1$). 
After possibly flyping, we can write $\T$ as the union of a positive braid 
in $S^2\times I$, with at least one vertical twist, and a reduced alternating 
tangle $\T'$. 
See \fullref{algknots2.9}. 

\begin{figure}[ht!]
\centerline{\includegraphics[height=1.35truein]{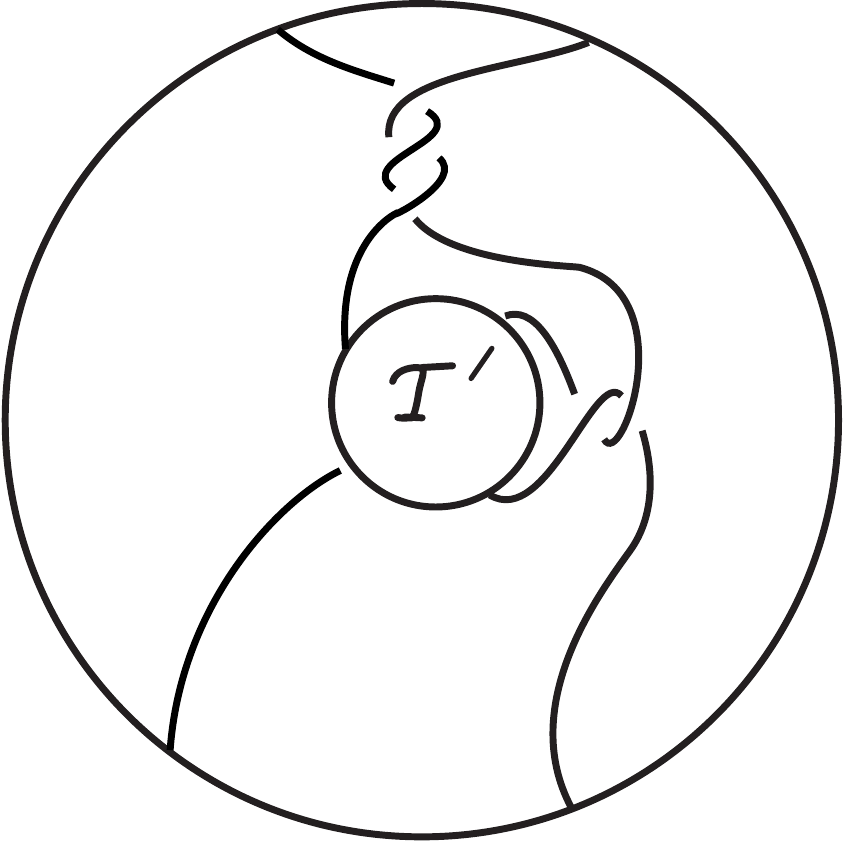}}

\caption{} 
\label{algknots2.9}
\end{figure}

After an isotopy we may assume that $F_2$ intersects the boundary of $\T'$ 
in a single circle, thereby writing $F_2$ as the union of an annulus in 
$S^2\times I$ and an essential Conway disk, $F'_2$, in $\T'$. 
By \cite[page 326]{Th1}, $F'_2$ can be taken to be visible, or hidden in a 
very special way.
If visible, then its slope on $\partial \T'$ (with coordinates from the 
diagram) is either $\frac01$ or $\frac10$. 
Since the braiding in $S^2\times I$ is positive with at least one vertical 
twist, this means the slope of $F_2$ on $\T$ cannot be $\frac10$, a 
contradiction (note that the flyping did not change the slope of $F_2$).

\begin{figure}[ht!]
\centerline{\includegraphics[height=1.0truein]{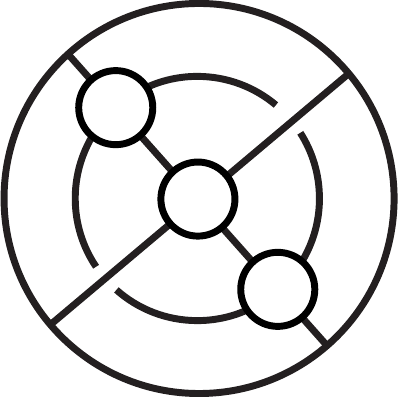}}

\caption{} 
\label{algknots4}
\end{figure}

So we assume $F'_2$ is hidden in $\T'$. 
But then \cite{Th1} shows that $\T'$ has a subdiagram as in 
\fullref{algknots4} and shows that the slope of $F'_2$ on 
$\partial \T'$ is either $\frac01,\frac10,\frac11$ ($-\frac11$ does not 
occur since the braiding  is positive). 
Again, the fact that the braiding in $S^2\times I$ is positive with at least
one vertical twist guarantees that the slope of $F_2$ on $\partial\T$ is 
not $\frac10$ as we have assumed.
\end{proof}

\begin{figure}[ht!]
\centerline{\includegraphics[height=1.55truein]{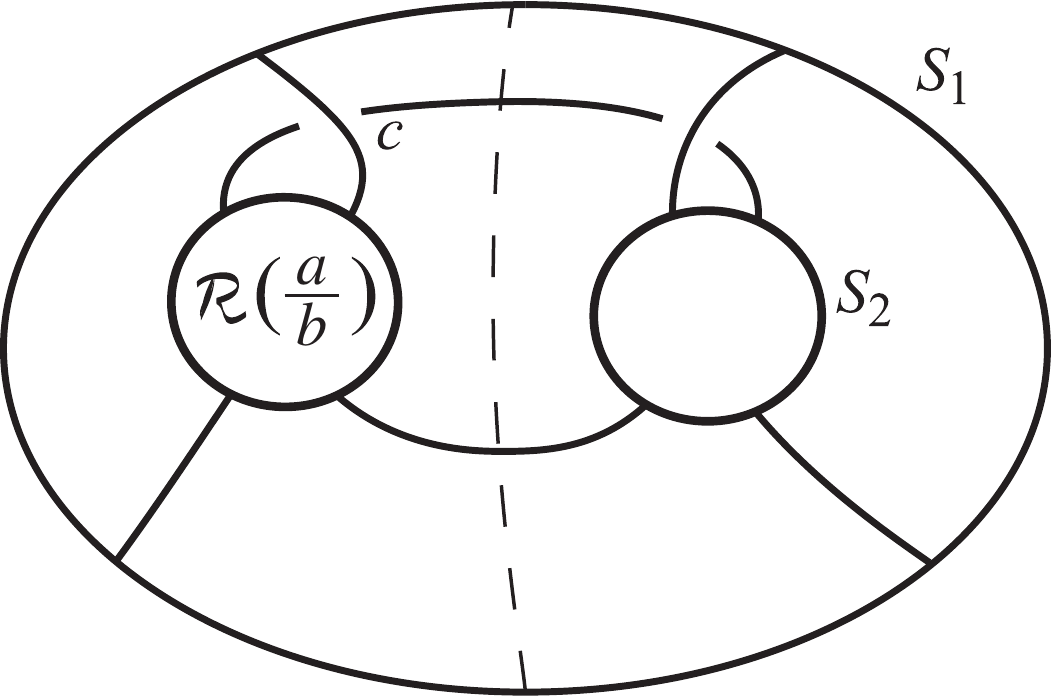}}
\caption{}		
\label{algknots3}
\end{figure}

\fullref{lem:visibleII} allows us to say that, after isotoping to  make 
$S_1$ visible, the elementary tangle of type~II becomes visible as in 
\fullref{algknots3}. 
In particular, we may assume the unknotting arc lies to the left of $F_1$ 
in this figure. 
By \fullref{thm8.3} and \fullref{lem:crossingchange}, 
this can be achieved by a crossing change in the 
(alternating) diagram on the left of \fullref{algknots3}.
If the crossing changed is not $c$, then this crossing change corresponds 
to a crossing change in the original diagram $D$.
If the crossing changed is $c$, then $a/b = -1$. 
As argued above, changing the crossing $c$ is equivalent to changing the 
single crossing in $\R(-1)$, which corresponds to a crossing in $D$.

This finishes the proof of \fullref{thm8.4}   	
in Case (1).
Cases (2) and (3) of the theorem are proved in Theorems~\ref{thm8.5} and 
\ref{thm8.6}.
\qed

\begin{thm}\label{thm8.5}
Let $K$ be an EM--knot. 
Then $K$ can be unknotted by a crossing change in any alternating 
diagram of $K$.
\end{thm}

\begin{proof} 
Again by the Flyping Conjecture \cite{MT}, we need only show that $K$ 
can be unknotted in some alternating diagram.
Let $K = K(\ell,m,n,p)$. 
By \cite[Proposition 1.4]{E1}, we may assume, by taking the mirror-image 
of $K$ if necessary, that $\ell>1$.
Lemmas~\ref{lem:9.7} and \ref{lem:9.8} imply that 
$K\!=\! .a.b.c = \S(\frac{-1}{c+1},\frac{a}{a+1};\frac1{b+1},\frac{-1}2)$, 
where $a,b,c>0$. 
Thus $K$ has a diagram of the form shown in \fullref{EMknot}.
\begin{figure}[ht!]
\centerline{\includegraphics[height=1.8truein]{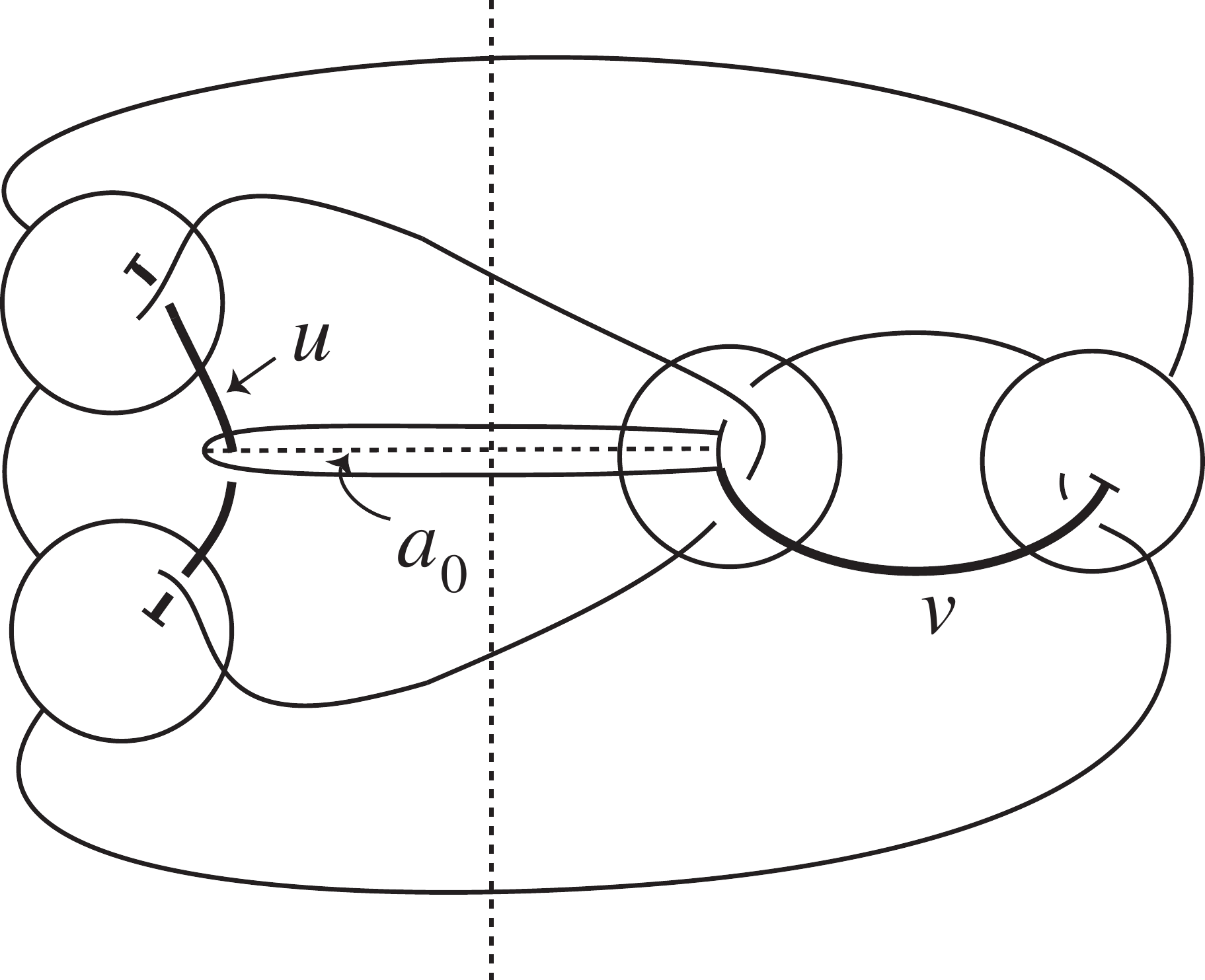}}
\caption{}		
\label{EMknot}
\end{figure}
The unknotting arc $a_0$ as well as 
the corresponding crossing move that unknots $K$ are shown. 
Let $u$ and $v$ be the arcs of the diagram indicated by the bold lines. 
Swinging $u$ ``under'' and $v$ ``over'' gives the diagram shown in 
\fullref{EMunknot}.
\begin{figure}[ht!]
\centerline{\includegraphics[height=1.6truein]{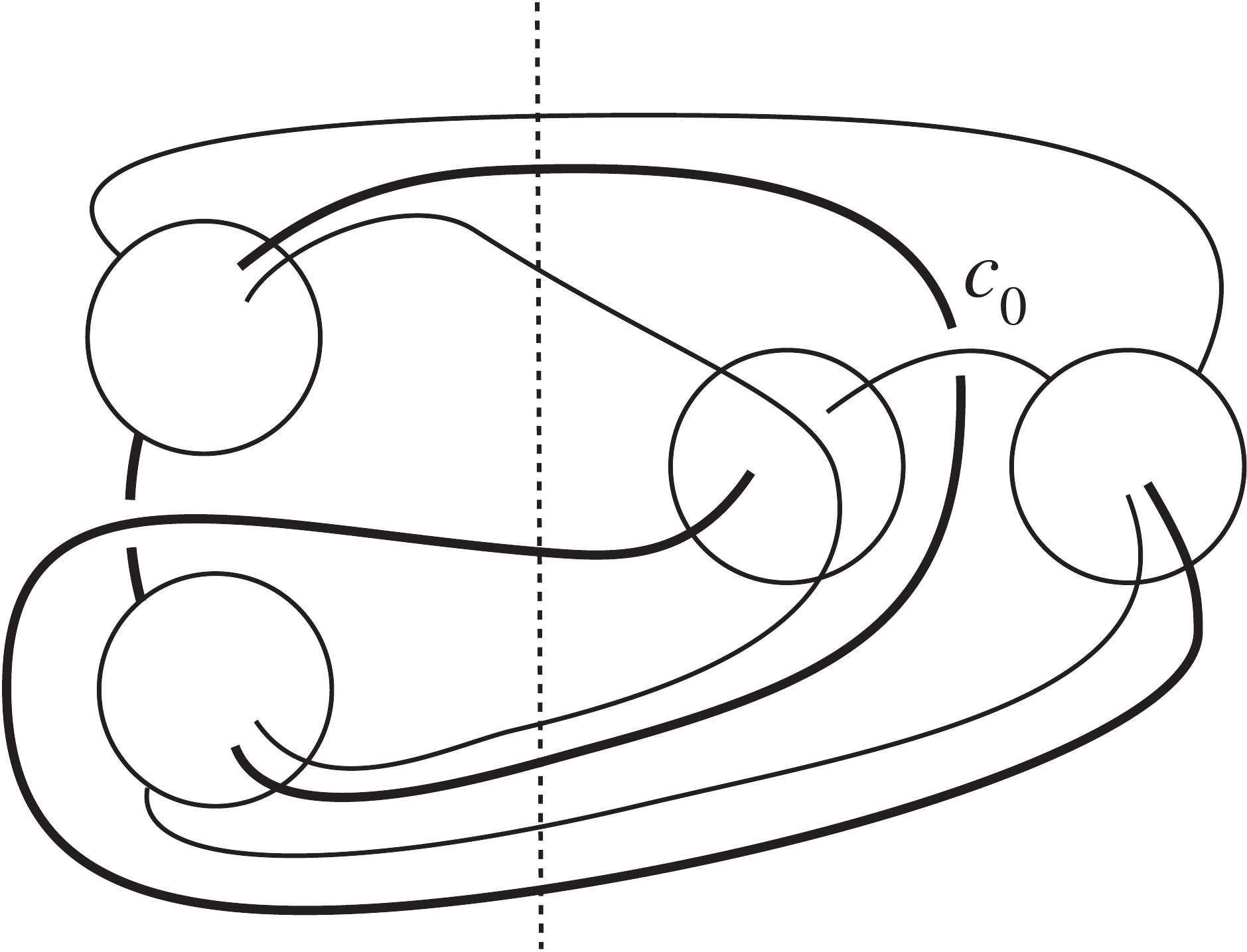}}
\caption{}		
\label{EMunknot}
\end{figure}
This diagram is alternating since $a,b,c >0$. 
Also, changing the crossing $c_0$ shown in \fullref{EMunknot} has the 
same effect as performing the crossing move shown in \fullref{EMknot}.
\end{proof}

\begin{thm}\label{thm8.6} 
Let $K$ be an alternating algebraic knot containing an EM--tangle whose
boundary is an essential Conway sphere. 
If $u(K)=1$, then $K$ can be unknotted by a crossing change in any 
alternating diagram of $K$.
\end{thm}

\proof
Again, by \cite{MT} we need only verify this for some alternating diagram.
Let $K$ be the union of essential tangles $\P\cup \P_0$ where $\P_0$ is an 
EM--tangle $\A_\ep(\ell,m)$. 
Applying \fullref{thm:Main}, we are either in Case~(1) or (3).  
In Case~(1), the proof of \fullref{thm8.4} guarantees the existence 
of an unknotting crossing change in an alternating diagram. 
Thus we assume we are in Case~(3).

\begin{lem}\label{lem:partialP}
In an alternating diagram of $K$, $\partial \P_0$ must be visible.
\end{lem}

\begin{proof} 
Assume not. 
Then there is a diagram of $\P_0$ as in \fullref{algknots3}. 
By \cite{Th1}, the disk $F_1$ in that figure is a Conway disk for $\P_0$. 
Thus it is the unique Conway disk for $\P_0$. 
Each side of $F_1$ is an alternating tangle. 
But this contradicts the fact that capping $\A_\ep (\ell,m)$ along slope 
$\frac10$ (ie, taking the denominator closure) gives either the unknot 
or Hopf link (by inspection of \fullref{Decomp}). 
\end{proof}

Thus $\partial\P_0$ is visible and let $D_0$ be the corresponding subdiagram.
This allows us to regard $\P_0$ as a marked tangle.
By \fullref{cor11.2} we need to find some alternating diagram of
this marked tangle  that exhibits a crossing change which unknots $K$.

\fullref{thm:Main}(3) gives an (unmarked) 
tangle homeomorphism $h\co \P_0 \to \A_\ep 
(\ell,m)\allowbreak
= \T_\ep (\ell,m) (1/2)$ which identifies the standard crossing 
move on $\A_\ep (\ell,m)$ as an unknotting, crossing move for $K$. 
After possibly rotating,  reflecting, or applying a mutation involution 
to $\A_\ep (\ell,m)$, there is an alternating braided tangle $\C$ in 
$S^2\times I$ such that extending $\A_\ep (\ell,m)$  by adjoining $\C$ 
gives a marked tangle 
$\A_\ep (\ell,m) \cup \C$ which is equivalent to $\P_0$ 
as a marked tangle --- via an extension of $h$.

\begin{Defin}
Let $E$ be the diagram in a disk $\Delta$ of a tangle in a 3--ball. 
A crossing $c$ of $E$ is said to be {\em inessential\/} iff there is a 
properly embedded arc in $\Delta$ that intersects $E$ only in $c$, dividing 
the four arcs of $E$ at $c$ into pairs. 
A diagram is {\em reduced\/} iff it contains no inessential crossings.
\end{Defin}

\begin{figure}[ht!]\small
\centerline{\includegraphics[height=1.3truein]{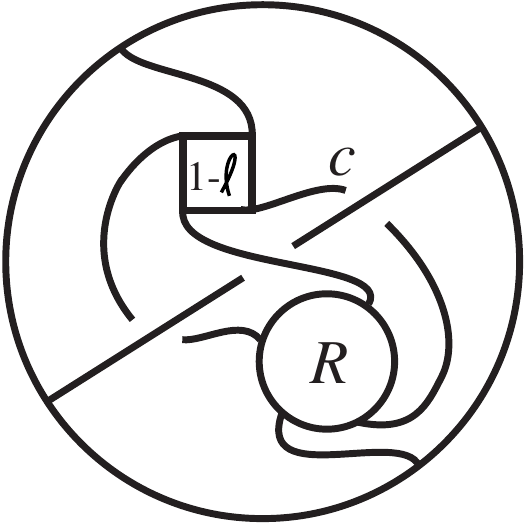}\hskip1.0truein 
\includegraphics[height=1.3truein]{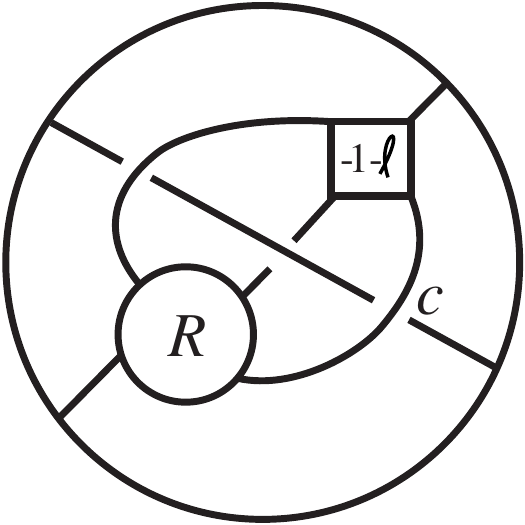}}
\centerline{$\ell>0$\hskip2.3truein $\ell <0$}
\vskip.2in
\centerline{\hskip1.0truein $R = \R \Big(\frac{(1-\ell)m+1}{m}\Big)$\hfill
$R = \R \Big(\frac{(-1-\ell)m+1}{m}\Big)$\hskip1.0truein}
\centerline{crossing change at $c\to\R\Big(\frac{1-2m}{1+2m}\Big)$
\hskip.3truein 
crossing change at $c\to \R\Big(\frac{1-2m}{1+2m}\Big)$}
\vskip.1in
\centerline{$\A_1 (\ell,m)$}
\vskip.2in
\centerline{\includegraphics[height=1.3truein]{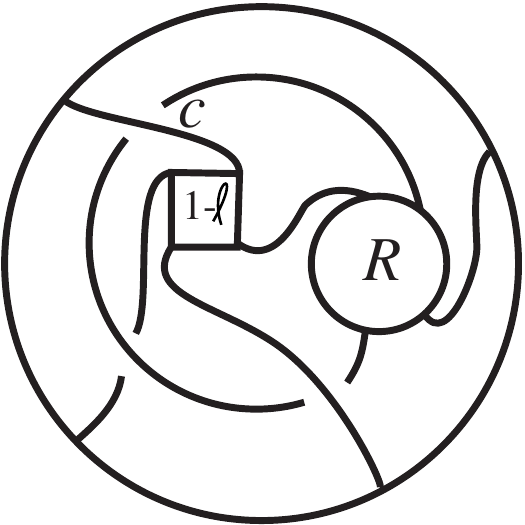}\hskip1truein 
\includegraphics[height=1.3truein]{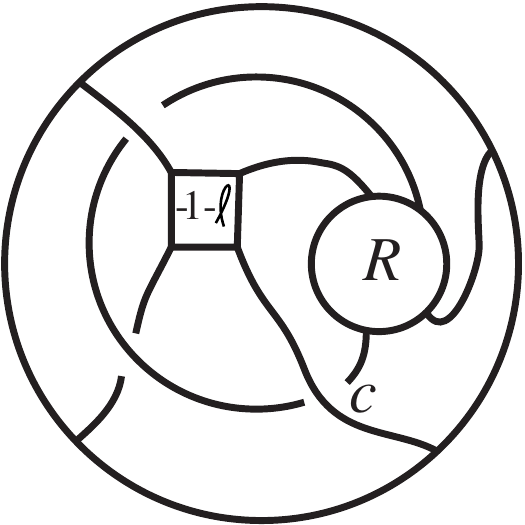}}
\centerline{$\ell>0$\hskip2.3truein $\ell <0$}
\vskip.2in
\centerline{\hskip1truein $R = \R \Big(\frac{m-1}{m}\Big)$\hfill
$R = \R \Big(\frac{m-1}{m}\Big)$\hskip1truein}
\vskip.05in
\centerline{\hskip1truein crossing change at\hfill 
crossing change at\hskip0.93truein} 
\vskip.05in
\centerline{\hskip1truein$c\to\R\Big(\frac{(\ell-1)(1-2m)+2}{1-2m}\Big)$
\hfill
$c\to \R\Big(\frac{2m-1}{(-1-\ell)(2m-1)+2}\Big)$\hskip0.5truein}
\vskip.1in
\centerline{$\A_2 (\ell,m)$}
\caption{Reduced, alternating diagrams of $\A_\ep (\ell,m)$}
\label{Crch}
\end{figure}

\fullref{Crch} displays a reduced, prime, alternating diagram, $E$, 
for $\A_\ep (\ell,m)$ which exhibits a crossing change sending 
$\A_\ell(\ell,m)$ to the rational tangle \break $\T_\ep (\ell,m)(\frac10)$.
Letting $C$ denote an alternating braided diagram of  the 
braided tangle $\C$, we have that $E\cup C$ is a prime diagram for the 
marked tangle $\P_0$, which also has the diagram $D_0$.
If $E$ and $C$ have the same sign (defined in Section~2) 
as alternating tangles, then $E\cup C$ is the desired alternating 
representative of $\P_0$, exhibiting the appropriate crossing change.
If $C$ and $D_0$ have opposite sign, then 
adjoining to $D_0$ the diagram $\bar C$ of the reverse braiding of $\C$ 
gives an alternating diagram $D_0\cup \bar C$ for the marked tangle 
$\A_\ep (\ell,m)$, which also has the diagram $E$.
Since $E$ is reduced, so is $D_0\cup \bar C$.
(For, we may assume $D_0$ has no 
nugatory crossings, hence neither does $D_0\cup \bar C$. 
Then $E$ and $D_0\cup \bar C$ have the same crossing number. 
But then a reducing arc for $D_0\cup \bar C$ suggests a capping of 
$D_0\cup \bar C$ and of $E$ giving two alternating diagrams of the same 
link with the same number of crossings --- one of which contains a 
nugatory crossing, the other does not). 
That is, $C$ is empty and $E$, $D_0$ represent the same marked tangle.
Thus $E$ is the sought after diagram for $\P_0$.

Thus we assume $E$ and $C$ have opposite sign, and $C$ and $D_0$ have the
same sign. 
Let $F$ be the diagram outside of $D_0$ in the alternating diagram of $K$. 
\fullref{Crch} shows that after the crossing change in $E$ the sign of the 
rational number corresponding to $\T_\ep (\ell,m)(\frac10)$ is opposite 
to the sign of $E$. 
That is, the result of the crossing change has an alternating diagram $R$ 
such that $R\cup C$ is an alternating diagram with the same sign as $D_0$.
Then  the diagram $R\cup C\cup F$ is an alternating diagram of the unknot. 
Hence it must have a nugatory crossing. 
Since we may assume neither $F$, $C$ or $R$ contain nugatory crossings, 
this means that either $F$ or $R$ is a split diagram (ie, there is a 
properly embedded arc that separates the arcs of the tangle). 
Since $K$ cannot be a connected sum (by \cite{S1}), and since $F$ is 
not rational, $R$ must represent $\R(0/1)$ or $\R(1/0)$.
But \fullref{Crch} shows this is not true.
\qed

\noindent{\bf Remarks about \fullref{Crch}}

(1)\qua  The markings have been changed between 
\fullref{Decomp} and \fullref{Crch}. 
In particular, to get from the markings in \fullref{Crch} to those in 
\fullref{Decomp} add the following braiding outside of \fullref{Crch}:
\begin{equation*}
\begin{split}
\A_1(\ell,m) ,\ \ell>0\ :\ &-1\ \text{ vertical twist below}\\
\A_1(\ell,m) ,\ \ell<0\ :\ &-1\ \text{ vertical twist above}\\
\A_2(\ell,m) ,\ \ell>0\ :\ &+1\ \text{ vertical twist below, then}\\
&-1\ \text{horizontal twist  to left}\\
\A_2(\ell,m)  ,\ \ell<0\ :\ &+1\ \text{ vertical twist below}
\end{split}
\end{equation*}

(2)\qua The diagrams  for $\A_1(\ell,m)$ were obtained by the moves in 
Figures~\ref{EMknot} 
and \ref{EMunknot}   
as well as twistings that change the markings.

\bibliographystyle{gtart}
\bibliography{link}

\end{document}